\documentclass[11pt]{article}
\usepackage{bbm}
\usepackage{amsmath,amsthm,amssymb}
\usepackage{graphicx}
\usepackage{tikz}
\usepackage[utf8]{inputenc}
\usepackage{pgfplots}
\usepackage{hyperref}
\usepackage{setspace}
\usepackage{pict2e,color}
\usepackage{multirow}
\usepackage{booktabs}
\usepackage{amsfonts}
\usepackage{verbatim} %for adding comments
\usepackage{sectsty}
\usepackage{url}
\usepackage{empheq}
\usepackage[normalem]{ulem}
\usepackage[titletoc,title]{appendix}
\usepackage[flushleft]{threeparttable}
\usepackage{soul} % use this (many fancier options)
\usepackage{natbib}
\usepackage{algorithm,algorithmic}
\usepackage{microtype}
\usepackage{subcaption}
\usepackage{pbox}
\usepackage{etoolbox}
\preto\subequations{\ifhmode\unskip\fi}
\usepackage{booktabs}
\usepackage{multirow}
\theoremstyle{definition}
\newtheorem{definition}{Definition}
\newtheorem{theorem}{Theorem}

\newtheorem{proposition}{Proposition}
\newtheorem{corollary}{Corollary}

\newtheorem{remark}{Remark}

\newtheorem{example}{Example}

\allowdisplaybreaks

\pgfplotstableread{
x         y      y-min y-max
2	0.081269906 0.009200852	0.133721914
3 0.120356703 0.014524797	0.178521315
4 0.138398513 0.018497724	0.218103121
5 0.141935673 0.017284335	0.23073966
}{\RVMSC}

\pgfplotstableread{
x         y      y-min y-max
3	0.081269906 0.009200852	0.133721914
4 0.121580874 0.023632245	0.160094506
5 0.161268633 0.033261422	0.178087752
6 0.181477098 0.041630059	0.210996841
7 0.206564113 0.053631882	0.223229032
8 0.228041379 0.060038851	0.233105214
9 0.232445938 0.068962122	0.23254455
10 0.249707189 0.075523453	0.283562069
}{\RVMST}

\pgfplotstableread{
x         y      y-min y-max
0	0.121421711 0.016129058	0.118468225
0.2 0.103277735 0.014756792	0.137113002
0.4 0.090164587 0.013856631	0.123657441
0.6 0.077791966 0.01400838	0.12892815
0.8  0.058526046 0.008533759	0.131525212
1 0.042120239 0.007702227	0.14111562
}{\RVMSL}

\pgfplotstableread{
x         y      y-min y-max
0.2	0.033305674 0.005678905	0.07834461
0.4 0.058298967 0.00857199	0.103492785
0.6 0.079716841 0.011289554	0.11776404
0.8 0.081269906 0.009200852	0.133721914
}{\RVMSS}

\pgfplotstableread{
x         y      y-min y-max
2	0.047249868	0.006051757	0.102019723
3 0.059624061	0.007823007	0.148719391
4 0.079003529	0.008269176	0.16059111
5 0.082601612	0.008026563	0.187532019
}{\RVMSCSAA}

\pgfplotstableread{
x         y      y-min y-max
3	0.047249868	0.006051757	0.102019723
4 0.043323422	0.007551923	0.11485011
5 0.042672586	0.008571207	0.101717359
6 0.042325944	0.011168122	0.10782674
7 0.044753596   0.011634766	0.106264134
8 0.045210809   0.013142589	0.104368796
9 0.04005351 0.004550389	0.064694648
10 0.062513911  0.014406163	0.0789229
}{\RVMSTSAA}

\pgfplotstableread{
x         y      y-min y-max
0	0.081540354	0.010495555	0.111406822
0.2 0.072638059	0.009560715	0.101594619
0.4 0.051665765	0.005578743	0.104718884
0.6 0.036500582	0.004295405	0.096225383
0.8  0.017707293	0.00195222	0.112466415
1 0.00196813	0.000604981	0.116687747
}{\RVMSLSAA}

\pgfplotstableread{
x         y      y-min y-max
0.2	0.015487577	0.002218362	0.058901107
0.4 0.028875099	0.003801395	0.090988218
0.6 0.03583293	0.004911582	0.102559832
0.8 0.047249868	0.006051757	0.102019723
}{\RVMSSSAA}
\title{On the Value of Risk-Averse Multistage Stochastic Programming in Capacity Planning}
\author{Xian Yu\thanks{Corresponding author; Department of Integrated Systems Engineering, The Ohio State University, Columbus, OH, USA, Email: {\tt yu.3610@osu.edu};}~~~and Siqian Shen\thanks{Department of Industrial and Operations Engineering, University of Michigan at Ann Arbor, USA, Email: {\tt siqian@umich.edu}.}}
\date{}

\begin{document}
%%%%%%%%%%%%%%%%
\maketitle

\begin{abstract}
We consider a risk-averse stochastic capacity planning problem under uncertain demand in each period. Using a scenario tree representation of the uncertainty, we formulate a multistage stochastic integer program to adjust the capacity expansion plan dynamically as more information on the uncertainty is revealed. Specifically, in each stage, a decision maker optimizes capacity acquisition and resource allocation to minimize certain risk measures of maintenance and operational cost. We compare it with a two-stage approach that determines the capacity acquisition for all the periods up front.  Using expected conditional risk measures (ECRMs), we derive a \textit{tight} lower bound and an upper bound for the gaps between the optimal objective values of risk-averse multistage models and their two-stage counterparts. Based on these derived bounds, we present general guidelines on when to solve risk-averse two-stage or multistage models. Furthermore, we propose approximation algorithms to solve the two models more efficiently, which are asymptotically optimal under an expanding market assumption. We conduct numerical studies using randomly generated and real-world instances with diverse sizes, to demonstrate the tightness of the analytical bounds and efficacy of the approximation algorithms. We find that the gaps between risk-averse multistage and two-stage models increase as the variability of the uncertain parameters increases and decrease as the decision maker becomes more risk-averse. Moreover, stagewise-dependent scenario tree attains much higher gaps than stagewise-independent counterpart, while the latter produces tighter analytical bounds.
\end{abstract}

%Supplemental Material:
%Data Ethics & Reproducibility Note:

% Sample
%\KEYWORDS{Stochastic programming, Decision support,Uncertainty, Disaster response, Optimization}

% Fill in data. If unknown, outcomment the field
~\\
{\bf Keywords:} Multistage stochastic integer programming; risk-averse optimization; dynamic time-consistent risk measure; approximation algorithms; capacity planning
%\HISTORY{Received: Month DD, YYYY; Accepted: Month DD, YYYY; Published Online: Month DD, YYYY}

\maketitle
%%%%%%%%%%%%%%%%%%%%%%%%%%%%%%%%%%%%%%%%%%%%%%%%%%%%%%%%%%%%%%%%%%%%%%

% Text of your paper here

\section{Introduction}
\label{sec:intro}
The capacity planning problem, aiming to determine the optimal level of capacity acquisition and allocation, is a classical optimization problem solved in a broad spectrum of applications, including capacity expansion in supply chains \citep{aghezzaf2005capacity}, production capacity planning for semiconductor manufacturing \citep{swaminathan2000tool,hood2003capacity}, electric vehicles (EVs) charging station expansion \citep{bayram2015capacity}, and so on. In these applications, demand fluctuates spatially and temporally, and as the cost of resource allocation in each stage depends on the demand input, it is also random depending on the distribution of demand over time. To adapt to new demand, service providers need to expand the existing capacities in an efficient way. Therefore, estimating and utilizing uncertain demand in the decision processes of capacity planning is crucial for cost reduction and quality-of-service improvement.

In this paper, we focus on a finite time horizon, where the uncertain demand in each period is modeled by a random vector. A decision maker optimizes when and how to expand the capacity and how to allocate resources to meet the demand to minimize a certain risk measure of the maintenance and operational cost over multiple periods. We compare two modeling frameworks: two-stage stochastic programming and multistage stochastic programming. In the two-stage model, we determine capacity acquisition for all time periods up front as ``here-and-now'' decisions and decide optimal resource allocation as ``wait-and-see'' recourse decisions based on realized demand and acquired capacities. In the multistage model, the uncertain demand values are revealed gradually over time and the capacity-acquisition decisions are adapted to this process.

Both two-stage and multistage decision frameworks are commonly used in the stochastic capacity planning literature. We refer interested readers to \cite{luss1982operations, davis1987optimal,sabet2020strategic} for comprehensive reviews on capacity planning problems under uncertainty. 
To handle general parameter uncertainty, robust optimization and stochastic programming are two of the popular approaches. \cite{paraskevopoulos1991robust,aghezzaf2005capacity} discussed robust capacity planning problems, and  \cite{eppen1989or,fine1990optimal,berman1994stochastic} considered stochastic programming approaches using scenarios to model the uncertain parameters. Most of these approaches are based on the two-stage decision-making framework.
There are a limited number of papers considering multistage stochastic capacity planning problems. Both \cite{ahmed2003approximation} and  \cite{singh2009dantzig} formulated a multistage stochastic mixed-integer programming model for capacity-planning problems, where \cite{ahmed2003approximation} proposed an LP-relaxation-based approximation scheme with asymptotic optimality and \cite{singh2009dantzig} applied ``variable splitting'' and Dantzig-Wolfe decomposition to tackle the problem. 

Stochastic programming approaches mainly focus on minimizing the total cost on average. However, minimizing the expected cost does not necessarily avoid the rare occurrences of undesirably high cost, and in a situation where it is important to maintain reliable performance, we aim to evaluate and control the \textit{risk}. In particular, coherent risk measures \citep{artzner1999coherent} have been used in many risk-averse stochastic programs as they satisfy several natural and desirable properties.    \citet{schultz2006conditional,shapiro2009lectures,ahmed2006convexity,miller2011risk} extended two-stage stochastic programs with risk-neutral expectation-based objective functions to risk-averse ones. It becomes nontrivial to model multistage risk-averse stochastic programs because the risk can be measured based on 
the cumulative cost, a total sum of the individual risk from each stage, or the conditional risk in a nested way. When the risk is calculated in a nested way via dynamic risk measures, a desirable property is called \textit{time consistency} \citep{ruszczynski2010risk}, which ensures consistent risk preferences over time. 

\citet{pflug2005measuring} proposed a class of multiperiod risk measures, which can be calculated by solving a stochastic dynamic linear optimization problem, and they analyzed its convexity and duality structure. \citet{homem2016risk} extended the above risk measures to expected conditional risk measures (ECRMs) and proved some appealing properties. First, ECRMs, originally defined for each stage separately, can be rewritten in a nested form. Second, any risk-averse multistage stochastic programs with ECRMs using CVaR measure can be recast as a risk-neutral multistage stochastic program with additional variables and constraints, which can be efficiently solved using existing decomposition algorithms. Third, ECRMs are time consistent following the inherited optimality property proposed in  \citet{homem2016risk}. In this paper, we show that ECRMs also follow the definition of time consistency proposed in \cite{ruszczynski2010risk} in Appendix A in the Online Supplement. Due to these properties, we base our analysis on the dynamic time-consistent ECRMs using CVaR measure, but we note that all the analysis and results derived in this paper can be also applied to other CVaR-based risk measures, such as static CVaR or period-wise CVaR, which are not time consistent.

A natural question is then about the performance of the risk-averse two-stage and multistage stochastic programs. Specifically, we are interested in comparing the optimal objective values and computational effort between solving risk-averse two-stage and multistage models for stochastic capacity planning problems. Noting that the multistage models have larger feasible regions and thus will always have lower costs, we aim to bound the gap between optimal objective values of the two-stage and multistage models given specific risk measures of the cost and characteristics of the uncertainty. \cite{huang2009value} were the first to show analytical lower bounds for the value of multistage stochastic programming (VMS) compared to the two-stage approach for capacity planning problems with an \textit{expectation-based} objective function. They developed an asymptotically optimal approximation algorithm by exploiting a decomposable substructure in the problem. Our work differs from \citet{huang2009value} in the following aspects. First, instead of using expectation as the objective function, we consider a dynamic time-consistent risk measure (i.e., ECRM). Although the resulting risk-averse problem can be reformulated as a risk-neutral counterpart, the risk-neutral results derived from \citet{huang2009value} cannot be directly applied here due to a new structure induced by ECRM-related decision variables and constraints. Instead, we derive new bounds on the VMS under this risk-averse setting. Second, we propose both lower and upper bounds to the gaps between ECRM-based optimal objective values of two-stage and multistage models and provide an example to illustrate the tightness of the lower bound, where \citet{huang2009value} only provided lower bounds without tightness guarantee. Third, our approximation algorithms improve the one developed in \citet{huang2009value} by deriving closed-form solutions of the substructure problem in the risk-averse setting and strengthening the upper bound iteratively, where the one in \citet{huang2009value} relied on solving the substructure problem in each time. We further prove that the approximation schemes are asymptotically optimal when the demand increases over time.

There has been rich literature on proposing decomposition algorithms to solve risk-neutral and risk-averse multistage stochastic programs. \cite{pereira1991multi} were the first to develop Stochastic Dual Dynamic Programming (SDDP) algorithms for efficiently computing multistage stochastic linear programs under stagewise-independent scenario trees, which construct an under-approximation for the value function in each stage iteratively. We refer interested readers to \citet{philpott2008convergence, girardeau2014convergence, shapiro2011analysis, shapiro2013risk,guigues2016convergence} for studies on the convergence of the SDDP algorithm under different problem settings. Among them, \cite{shapiro2013risk,guigues2016convergence} extended the SDDP algorithm to solve a risk-averse multistage stochastic program. Stochastic Dual Dynamic integer Programming (SDDiP), first proposed by \citet{zou2019stochastic}, is an extension of SDDP to handle the nonconvexity arising in multistage stochastic integer programs based on a similar computational framework. 
The main goal of this paper is to bound the gap between two-stage and multistage stochastic programs and provide general guidelines on model choices. Using the substructure explored when deriving the bounds, we also propose approximation algorithms to solve risk-averse two-stage and multistage stochastic capacity planning problems more efficiently. Note that a common assumption made in the aforementioned SDDP-type algorithms is that the underlying stochastic process is stagewise independent. However, our analysis and proposed approximation algorithms are not restricted by this stagewise independence assumption. While our analytical results can be applied to a more general stochastic process, using a stage-wise independent scenario tree, we will compare the performance of our approximation algorithm with SDDiP algorithm in Section \ref{sec:time}.

The remainder of the paper is organized as follows. In Section \ref{sec:models}, we present the formulations of risk-averse two-stage and multistage capacity planning problems under ECRM. In Section \ref{sec:value-proof}, we derive two lower bounds and an upper bound for the gaps between the optimal ECRM-based objective values of two-stage and multistage stochastic capacity planning models. In Section \ref{sec:approx}, we propose approximation algorithms with performance guarantees. In Section \ref{sec:compu}, we conduct numerical studies using instances with diverse uncertainty patterns, to demonstrate the tightness of our bounds and the performance of the approximation algorithms. 
 Section \ref{sec:conclu} concludes the paper and states future research directions. Throughout the paper, we use bold symbols to denote vectors/matrices and use $[n]$ to denote the set $\{1,2,\ldots,n\}$.

\section{Problem Formulations}\label{sec:models}

Consider a general class of capacity planning problems, where we have $1,\ldots,T$ time periods, $1,\ldots,M$ facilities, and $1,\ldots,N$ customer sites. Let $c_{tij}$ be the cost of serving one unit of demand from customer site $j$ in facility $i$ at period $t$, $f_{ti}$ be the cost of maintaining one unit of resources in facility $i$ at period $t$, $h_{ti}$ be the number of customers that one unit of resources can serve in facility $i$ at period $t$, and $d_{tj}$ be the demand at customer site $j$ at period $t$ for all $i\in [M]$, $j \in [N]$, and $t \in [T]$. 

Define decision variables $x_{ti} \in \mathbb{Z}_+$ as the units of resources we invest in facility $i$ at period $t$, and decision variables $y_{tij}\in\mathbb{R}_+$ as the amount of demand from customer site $j$ served by facility $i$ at period $t$. 
The deterministic capacity planning problem is given by:
\begin{subequations}
\label{model:deter}
\begin{align}
\min\quad &\sum_{t=1}^T\sum_{i=1}^M f_{ti}\sum_{\tau=1}^tx_{\tau i}+\sum_{t=1}^T\sum_{i=1}^M\sum_{j=1}^N c_{tij}y_{tij} \label{eq:1}\\
\text{s.t.}\quad& \sum_{i=1}^M y_{tij}= d_{tj}, \ \forall j=1,\ldots,N,\ t=1,\ldots,T \label{eq:2}\\
&\sum_{j=1}^N y_{tij}\le h_{ti}\sum_{\tau=1}^t x_{\tau i},\ \forall i=1,\ldots,M,\ t=1,\ldots,T \label{eq:3}\\
&x_{ti}\in \mathbb{Z}_{+}, \ \forall i=1,\ldots,M,\ t=1,\ldots,T\label{eq:5}\\
&y_{tij}\in \mathbb{R}_{+}, \ \forall i=1,\ldots,M,\ j=1,\ldots,N,\ t=1,\ldots,T.\label{eq:6}
\end{align}
\end{subequations}
The objective function \eqref{eq:1} minimizes the total maintenance and operational cost over all time periods, where $\sum_{\tau=1}^tx_{\tau i}$ represents the number of units of resources available in facility $i$ at period $t$. Constraints \eqref{eq:2} require all the demands to be satisfied in each period. Constraints \eqref{eq:3} indicate that in each period, we can only serve customer demand within the current capacity of each facility.

In the remainder of this paper, we will work with a generic capacity planning model in vector forms below: 
\begin{subequations}\label{model:vector}
\begin{align}
\min_{\substack{\boldsymbol{x}_1,\ldots,\boldsymbol{x}_T\\\boldsymbol{y}_1,\ldots,\boldsymbol{y}_T}}\quad &\sum_{t=1}^T \boldsymbol{f}_{t}^{\mathsf T} \sum_{\tau=1}^t \boldsymbol{x}_{\tau}+\sum_{t=1}^T\boldsymbol{c}_t^{\mathsf T} \boldsymbol{y}_t\\
\text{s.t.}\quad& \boldsymbol{A}_t\boldsymbol{y}_t= \boldsymbol{d}_t, \ \forall t=1,\ldots,T \label{eq:2-2}\\
&\boldsymbol{B}_{t}\boldsymbol{y}_t\le \sum_{\tau=1}^t \boldsymbol{x}_{\tau},\ \forall t=1,\ldots,T \label{eq:2-3}\\
&\boldsymbol{x}_t\in \mathbb{Z}^{M}_{+}, \ \boldsymbol{y}_t\in \mathbb{R}^{MN\times 1}_{+},\ \forall t=1,\ldots,T,\nonumber
\end{align}
\end{subequations}
where $\boldsymbol{x}_t\in \mathbb{Z}^{M}_{+}$ represents the capacity-acquisition decisions (in number of units) for the resources and $\boldsymbol{y}_t\in \mathbb{R}^{MN\times 1}_{+}$ represents the operational level allocation of resources in period $t$. Matrices $\boldsymbol{A}_t\in\mathbb{R}^{N\times MN}$ and $\boldsymbol{B}_{t}\in\mathbb{R}^{M\times MN}$ correspond to the coefficients in Constraints \eqref{eq:2} and \eqref{eq:3}, respectively (i.e., $\boldsymbol{A}_{t}$ is a 0-1 matrix corresponding to Constraints (1b), and $\boldsymbol{B}_{t}$ is a matrix with entries 0 or $1/h_{ti}$ corresponding to Constraints (1c)).
Note that the difference between Model \eqref{model:vector} and Model (7)-(10) in \cite{huang2009value} is that instead of considering one-time acquisition costs of resources, we take into account the period-wise maintenance costs for all acquired resources, which is more realistic in applications such as EV charging station maintenance \citep{EVMaintenance}, power system equipment maintenance \citep{PowerMaintenance}, etc. 

\subsection{Scenario Tree Representation} 
In Model \eqref{model:vector}, the data we acquire at each period $t$ is the demand $\boldsymbol{d}_t\in\mathbb{R}_+^N$ for all $t=1,\ldots, T$. We consider that the data series $\{\boldsymbol{d}_2,\ldots,\boldsymbol{d}_T\}$ evolve according to a known probability distribution, and $\boldsymbol{d}_1\in \mathbb R^N_+$ is deterministic (see similar assumptions made in, e.g.,  \cite{huang2009value,shapiro2009lectures,zou2019stochastic}). In practice, $\boldsymbol{d}_1$ can be derived and forecasted as the average of historical demand. Consider a discrete distribution and assume that the number of realizations is finite, where such an approximation can be constructed by Monte Carlo sampling if the probability distribution is instead continuous and the resultant problem is called Sample Average Approximation (SAA) \citep{kleywegt2002sample}.  

In a multistage stochastic setting, the uncertain demand is revealed gradually, where we need to decide both capacity acquisition $\boldsymbol{x}_t$ and resource allocation $\boldsymbol{y}_t$ in each stage $t$ based on the currently realized demand $\boldsymbol{d}_t$. Correspondingly, the decision-making process can be described as follows:
\begin{align*}
\underbrace{\text{decision}\ (\boldsymbol{x}_1, \boldsymbol{y}_1)}_{\text{Stage 1}}&\to \underbrace{\text{observation}\ (\boldsymbol{d}_2)\to \text{decision}\ (\boldsymbol{x}_2, \boldsymbol{y}_2)}_{\text{Stage 2}}\to \text{observation}\ (\boldsymbol{d}_3)\to\\
\cdots&\to \text{decision}\ (\boldsymbol{x}_{T-1}, \boldsymbol{y}_{T-1})\to \underbrace{\text{observation}\ (\boldsymbol{d}_{T})\to \text{decision}\ (\boldsymbol{x}_T, \boldsymbol{y}_T)}_{\text{Stage}\ T}.
\end{align*}

\begin{figure}[ht!]
    \centering
			\includegraphics[width=0.3\textwidth]{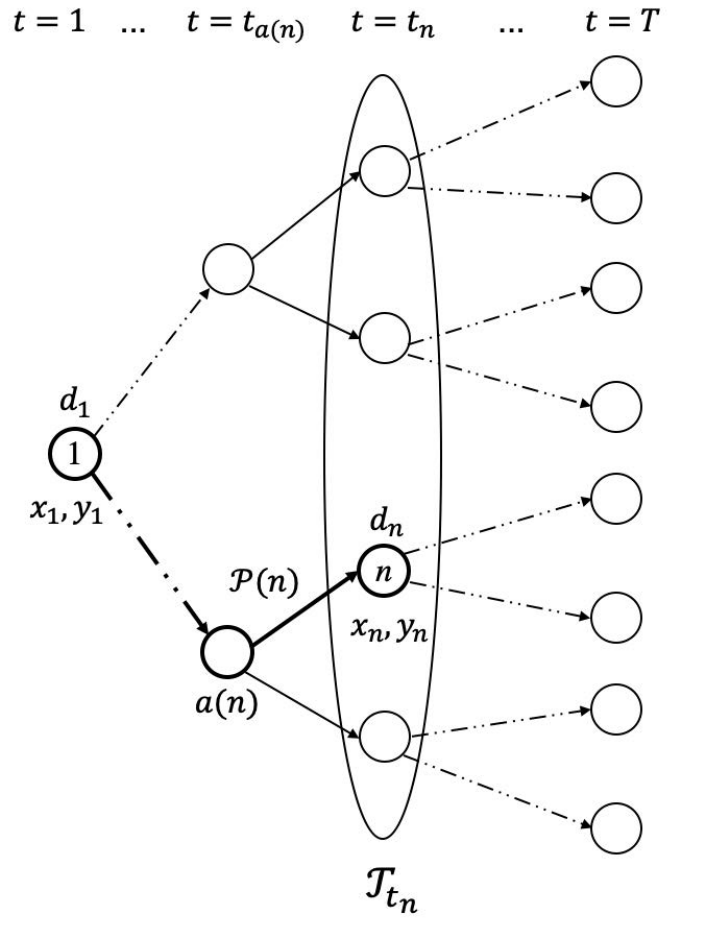}
    \caption{Illustration of a scenario tree and its related notation.}
    \label{fig:scenario-tree}
\end{figure}

To facilitate formulating the stochastic programs, we introduce a scenario tree representation of the uncertainty and decision variables \citep{shapiro2009lectures}. Let $\mathcal{T}$ be {the set of all nodes in} the scenario tree associated with the underlying stochastic process (see Figure \ref{fig:scenario-tree}). Each node $n$ in period $t>1$ has a unique parent node $a(n)$ in period $t-1$, and the set of children nodes of node $n$ is denoted by $\mathcal{C}(n)$. The set $\mathcal{T}_t$ denotes all the nodes corresponding to time period $t$, and $t_n$ is the time period corresponding to node $n$. Specifically, all nodes in the last period $\mathcal{T}_T$ are referred to as leaf nodes and denoted by $\mathcal{L}$.  The set of all nodes on the path from the root node to node $n$ is denoted by $\mathcal{P}(n)$. 
Each node $n$ is associated with a (unconditional) probability $p_n$, which is the probability of the realization of the $t_n$-period data sequence $\{\boldsymbol{d}_m\}_{m\in\mathcal{P}(n)}$. The probabilities of the nodes in each period sum up to one, i.e., $\sum_{n\in\mathcal{T}_t}p_n = 1,\ \forall t=1,\ldots, T$, and the probabilities of all children nodes sum up to the probability of the parent node, i.e., $\sum_{m\in\mathcal{C}(n)}p_m=p_n,\ \forall n\not\in\mathcal{L}$.
If $n$ is a leaf node, i.e., $n\in\mathcal{L}$, then $\mathcal{P}(n)$ corresponds to a scenario that represents a joint realization of the problem parameters over all periods. For each leaf node $n\in\mathcal{L}$, we denote the probability of the scenario $\mathcal{P}(n)$ as $p_n$.
We denote the capacity-acquisition decision variable at node $n$ by $\boldsymbol{x}_n$, and the resource-allocation decision variable at node $n$ by $\boldsymbol{y}_n$. {Figure \ref{fig:scenario-tree} depicts a scenario tree representation of $T$-period uncertainty and its related notation.}

\subsection{Expected Conditional Risk Measures (ECRMs)}
Next, we introduce the definition and key properties of coherent risk measures. Consider a probability space $(\Xi,\mathcal{F},P)$, and let $\mathcal{F}_1\subset\mathcal{F}_2\subset\ldots\subset\mathcal{F}_T$ be sub-sigma-algebras of $\mathcal{F}$ such that each $\mathcal{F}_t$ corresponds to the information available up to (and including) period $t$, with $\mathcal{F}_1=\{\emptyset,\Xi\}, \ \mathcal{F}_T=\mathcal{F}$. Let $\mathcal{Z}_t$ denote a space of $\mathcal{F}_t$-measurable functions from $\Xi$ to $\mathbb{R}$, and let $\mathcal{Z}:=\mathcal{Z}_1\times\cdots\times\mathcal{Z}_T$, e.g., $\boldsymbol{d}_{t}\in \mathcal{Z}_{t},\ \forall t=1,\ldots,T$. 

\begin{definition}\citep{artzner1999coherent}
A conditional risk measure $\rho_{t}^{\boldsymbol d_{[t-1]}}: \mathcal{Z}_{t}\to\mathcal{Z}_{t-1}$ is coherent if it satisfies the following four properties: (i) Monotonicity: If $Z_1,Z_2\in\mathcal{Z}_{t}$ and $Z_1\ge Z_2$, then $\rho_{t}^{\boldsymbol d_{[t-1]}}(Z_1)\ge\rho_{t}^{\boldsymbol d_{[t-1]}}(Z_2)$; (ii) Convexity: $\rho_{t}^{\boldsymbol d_{[t-1]}}(\gamma Z_1+(1-\gamma)Z_2)\le \gamma\rho_{t}^{\boldsymbol d_{[t-1]}}(Z_1)+(1-\gamma)\rho_{t}^{\boldsymbol d_{[t-1]}}(Z_2)$ for all $Z_1,\ Z_2\in\mathcal{Z}_{t}$ and all $\gamma\in[0, 1]$; (iii) Translation invariance: If $W\in\mathcal{Z}_{t-1}$ and $Z\in\mathcal{Z}_{t}$, then $\rho_{t}^{\boldsymbol d_{[t-1]}}(Z+W)=\rho_{t}^{\boldsymbol d_{[t-1]}}(Z)+W$; and (iv) Positive Homogeneity: If $\gamma\ge 0$ and $Z\in\mathcal{Z}_{t}$, then $\rho_{t}^{\boldsymbol d_{[t-1]}}(\gamma Z)=\gamma\rho_{t}^{\boldsymbol d_{[t-1]}}(Z)$. Here, $Z_1\ge Z_2$ if and only if $Z_1(\xi)\ge Z_2(\xi)$ for a.e. $\xi\in\Xi$.
\end{definition}
For notation simplicity, denote the period-wise cost $g_t(\boldsymbol{x}_{1:t},\boldsymbol{y}_t)= \boldsymbol{f}_{t}^{\mathsf T} \sum_{\tau=1}^t \boldsymbol{x}_{\tau}+\boldsymbol{c}_t^{\mathsf T} \boldsymbol{y}_t$ by $g_t$. We consider a multiperiod risk function $\mathbb{F}$ as a mapping from $\mathcal{Z}$ to $\mathbb{R}$ below: 
\small
\begin{align}
\mathbb{F}(g_1,\ldots,g_{T})=g_1+\rho_2(g_2)+\mathbb{E}_{\boldsymbol d_{[2]}}\left[{\rho_3^{\boldsymbol d_{[2]}}}(g_3)\right]+\mathbb{E}_{\boldsymbol d_{[3]}}\left[{\rho_4^{\boldsymbol d_{[3]}}}(g_4)\right]+\cdots+\mathbb{E}_{\boldsymbol d_{[T-1]}}\left[{\rho_{T}^{\boldsymbol d_{[T-1]}}}(g_{T})\right],\label{ECRMs}
\end{align}
\normalsize
where ${\rho_t^{\boldsymbol d_{[t-1]}}}$ is a {conditional} risk measure mapping from $\mathcal{Z}_t$ to $\mathcal{Z}_{t-1}$ to represent the risk given the information available up to (including) period $t-1$, {i.e., $\boldsymbol d_{[t-1]}=(\boldsymbol d_1,\ldots,\boldsymbol d_{t-1})$}. 
This class of multiperiod risk measures is called expected conditional risk measures (ECRMs), first introduced by \cite{homem2016risk}. We show that ECRMs \eqref{ECRMs} are time consistent following the definition in \cite{ruszczynski2010risk} in Appendix A in the Online Supplement. 

Using tower property of expectations \citep{wasserman2004all}, the multi-period risk function ECRM \eqref{ECRMs} can be written in a nested form below
\small
\begin{align}
\mathbb{F}(g_1,\ldots,g_{T})=g_1+\rho_2(g_2)+\mathbb{E}_{\boldsymbol d_{2}}\Big[{\rho_3^{\boldsymbol d_{[2]}}}(g_3)+\mathbb{E}_{\boldsymbol d_{3}|\boldsymbol d_{[2]}}\Big[{\rho_4^{\boldsymbol d_{[3]}}}(g_4)+\cdots
+\mathbb{E}_{\boldsymbol d_{T-1}|\boldsymbol d_{[T-2]}}\Big[{\rho_{T}^{\boldsymbol d_{[T-1]}}}(g _{T})\Big]\cdots\Big]\Big],\label{eq:obj}
\end{align}
\normalsize
where $\mathbb{E}_{\boldsymbol{d}_{t}|\boldsymbol d_{[t-1]}}$ represents the expectation with respect to the conditional probability distribution of $\boldsymbol{d}_{t}$ given realization $\boldsymbol d_{[t-1]}$.

For our problem, we consider a special class of single-period coherent risk measures -- a convex combination of conditional expectation and Conditional Value-at-Risk (CVaR) \citep[for the definition of CVaR, please refer to][]{rockafellar2000optimization}, i.e., for $t=2,\ldots,T$,
\begin{equation}
{\rho_{t}^{\boldsymbol d_{[t-1]}}}(g_t)=(1-\lambda_t)\mathbb{E}[g_t{|\boldsymbol d_{[t-1]}}]+\lambda_t {\text{CVaR}_{\alpha_t}^{\boldsymbol d_{[t-1]}}}[g_t],\label{eq:rho}
\end{equation}
where $\lambda_t\in[0,1]$ is a parameter that adjusts between optimizing on average and risk control, and $\alpha_t\in(0,1)$ is the confidence level. Notice that this risk measure is more general than CVaR and it includes CVaR as a special case when $\lambda_t=1$. 

Following the results by \cite{rockafellar2002conditional}, CVaR can be expressed as the following optimization problem:
\begin{equation}
\text{CVaR}_{\alpha_t}^{\boldsymbol d_{[t-1]}}[g_t]:=\inf_{\eta_t\in\mathbb{R}}\left\lbrace \eta_t+\frac{1}{1-\alpha_t}\mathbb{E}[[g_t-\eta_t]_{+}|\boldsymbol d_{[t-1]}]\right\rbrace,\label{eq:ccvar}
\end{equation}
where $[a]_{+}:=\max\{a,0\},$ and $\eta_t$ is an auxiliary variable. The minimum of the right-hand side of Eq. \eqref{eq:ccvar} is attained at $\eta_t^* = \text{VaR}_{\alpha_t}[g_t] := \inf\{v: \mathbb{P}(g_t\le v)\ge \alpha_t\}$, and thus CVaR is the mean of the upper $(1-\alpha_t)$-tail distribution of $g_t$, i.e., $\mathbb{E}[g_t|g_t>\eta_t^*]$. Selecting a large $\alpha_t$ value makes CVaR sensitive to rare but very high costs (we fix $\alpha_t=0.95$ and vary $\lambda_t$ to reflect different risk attitudes in our numerical experiments). To linearize $[g_t-\eta_t]_{+}$ in \eqref{eq:ccvar}, we replace it by a variable $u_t$ with two additional constraints: $u_t\ge 0,\ u_t\ge g_t-\eta_t$.

\subsection{Risk-Averse Two-Stage and Multistage Stochastic Programs}
Combining \eqref{eq:obj}, \eqref{eq:rho} and \eqref{eq:ccvar}, the objective function of our risk-averse multistage capacity planning model is specified as
\small
\begin{align}
\hspace{-0.2cm} z_R^{MS}&=\min_{\substack{{(\boldsymbol{x}_1,\boldsymbol{y}_1)\in {X}_1(\boldsymbol{d}_1)},\\\eta_2\in\mathbb{R}}}g_1(\boldsymbol{x}_1,\boldsymbol{y}_1)+\lambda_2\eta_2+\mathbb{E}_{\boldsymbol d_2}\Big[\min_{\substack{(\boldsymbol{x}_2,\boldsymbol{y}_2,u_2)\in{X}_2(\boldsymbol{x}_1,\eta_2,\boldsymbol{d}_2),\\\eta_3\in\mathbb{R}}}\Big\{\frac{\lambda_2}{1-\alpha_2}u_2+(1-\lambda_2)g_2(\boldsymbol{x}_{1:2},\boldsymbol{y}_2)+\lambda_3\eta_3+\cdots\nonumber\\
&+\mathbb{E}_{\boldsymbol d_{T-1}{|\boldsymbol{d}_{[1,T-2]}}}\Big[\min_{\substack{(\boldsymbol{x}_{T-1},\boldsymbol{y}_{T-1},u_{T-1})\in{X}_{T-1}(\boldsymbol{x}_{1:T-2},\eta_{T-1},\boldsymbol{d}_{T-1}),\\\eta_{T}\in\mathbb{R}}}\Big\{\frac{\lambda_{T-1}}{1-\alpha_{T-1}}u_{T-1}+(1-\lambda_{T-1})g_{T-1}(\boldsymbol{x}_{1:T-1},\boldsymbol{y}_t)+\lambda_T\eta_T\nonumber\\
&+\mathbb{E}_{\boldsymbol d_T{|{\boldsymbol{d}_{[1,T-1]}}}}\Big[\min_{(\boldsymbol{x}_T,\boldsymbol{y}_T,u_T)\in{X}_T(\boldsymbol{x}_{1:T-1},\eta_T,\boldsymbol{d}_T)}\Big\{\frac{\lambda_T}{1-\alpha_T}u_T+(1-\lambda_T)g_T(\boldsymbol{x}_{1:T},\boldsymbol{y}_T)\Big\}\Big]\cdots\Big\}\Big],\label{model:multirisk}
\end{align}
\normalsize
where the auxiliary variable $\eta_t\in\mathbb{R}$ is a function of $\boldsymbol{d}_1,\ldots,\boldsymbol{d}_{t-1}$, i.e., $\eta_t$ is a ``$(t-1)$-stage'' variable similar to $\boldsymbol{x}_{t-1}$ for all $t=2,\ldots,T$, and the auxiliary variable $u_t\in\mathbb{R}$ is a $t$-stage variable to represent the excess of $t$-stage cost of above $\eta_t$.
Sets ${X}_1(\boldsymbol{d}_1)=\{(\boldsymbol{x}_1, \boldsymbol{y}_1)\in  \mathbb{Z}^{M}_{+}\times \mathbb{R}^{MN\times 1}_{+}:  \boldsymbol{A}_1\boldsymbol{y}_1= \boldsymbol{d}_{1},\ \boldsymbol{B}_{1}\boldsymbol{y}_{1}\le \boldsymbol{x}_1\}$ and ${X}_t(\boldsymbol{x}_{1:t-1},\eta_t,\boldsymbol{d}_t)=\{(\boldsymbol{x}_t, \boldsymbol{y}_t,u_t)\in  \mathbb{Z}^{M}_{+}\times \mathbb{R}^{MN\times 1}_{+}\times\mathbb{R}_+:  \boldsymbol{A}_t\boldsymbol{y}_t= {\boldsymbol{d}_{t}},\ \boldsymbol{B}_t\boldsymbol{y}_t\le \sum_{\tau=1}^{t}\boldsymbol{x}_{\tau},\ u_t+\eta_t\ge g_t(\boldsymbol{x}_{1:t},\boldsymbol{y}_t)\},\ \forall t=2,\ldots,T$ are the feasible regions in each stage.

Plugging $g_t(\boldsymbol{x}_{1:t},\boldsymbol{y}_t)= \boldsymbol{f}_{t}^{\mathsf T} \sum_{\tau=1}^t \boldsymbol{x}_{\tau}+\boldsymbol{c}_t^{\mathsf T} \boldsymbol{y}_t$ into Eq. \eqref{model:multirisk} and using scenario-node-based notation, the risk-averse multistage model \eqref{model:multirisk} can be written in the following extensive form:
\begin{subequations}\label{model:multiriskaverse}
\begin{align}
z^{MS}_R=\min_{\substack{\boldsymbol{x}_n,\boldsymbol{y}_n, n\in \mathcal{T}\\\eta_n, n\not\in\mathcal{L}, u_n,n\not=1}} \quad& \sum_{n\in \mathcal{T}} p_n\left(\boldsymbol{\tilde{f}}_n^{\mathsf T}\sum_{m\in \mathcal{P}(n)}\boldsymbol{x}_m+\boldsymbol{\tilde{c}}_n^{\mathsf T}\boldsymbol{y}_n + \tilde{\lambda}_n \eta_n + \tilde{\alpha}_n u_n\right)\label{eq:multiriskaverse_obj}\\
\text{s.t.} \quad& \boldsymbol{A}_{t_n}\boldsymbol{y}_n= \boldsymbol{d}_n, \ \forall n\in \mathcal{T}\label{eq:multi-constraint_y}\\
& \boldsymbol{B}_{t_n}\boldsymbol{y}_n\le \sum_{m\in \mathcal{P}(n)}\boldsymbol{x}_m,\ \forall n\in \mathcal{T}\\
& u_n + \eta_{a(n)} \ge \boldsymbol{f}_{t_n}^{\mathsf T} \sum_{m\in \mathcal{P}(n)}\boldsymbol{x}_m+\boldsymbol{c}_{t_n}^{\mathsf T} \boldsymbol{y}_n,\ \forall n\not= 1 \label{eq:constraint_eta}\\
& \boldsymbol{x}_n\in \mathbb{Z}^{M}_{+},\ \boldsymbol{y}_n\in \mathbb{R}^{MN \times 1}_{+},\ \forall n\in \mathcal{T},\ u_n \ge 0,\ \forall n\not=1,\nonumber
\end{align}
\end{subequations}
where $\boldsymbol{\tilde{f}}_n = \boldsymbol{f}_{t_n}$ if $n = 1$ and $\boldsymbol{\tilde{f}}_n = (1-\lambda_{t_n})\boldsymbol{f}_{t_n}$ otherwise; $\boldsymbol{\tilde{c}}_n = \boldsymbol{c}_{t_n}$ if $n = 1$ and $\boldsymbol{\tilde{c}}_n = (1-\lambda_{t_n})\boldsymbol{c}_{t_n}$ otherwise; $\tilde{\lambda}_n = 0$ if $n\in \mathcal{L}$ and $\tilde{\lambda}_n = \lambda_{t_n+1}$ otherwise; $\tilde{\alpha}_n = 0$ if $n = 1$ and $\tilde{\alpha}_n = \frac{\lambda_{t_n}}{1-\alpha_{t_n}}$ otherwise.

Note that in the above multistage stochastic program, we have both capacity expansion $\boldsymbol{x}_n$ and allocation $\boldsymbol{y}_n$ decisions corresponding to each node $n$ in the scenario tree, $\eta_n$-variables are defined for all nodes $n$ except for the leaf nodes $\mathcal{L}$, and $u_n$-variables are defined for non-root nodes $n\not=1$. 
When $\lambda_t=0,\ \forall t=2,\ldots, T$, the risk measure $\rho_t$ becomes the expectation and \eqref{model:multiriskaverse} reduces to a risk-neutral multistage stochastic program.

In a two-stage stochastic program, we decide the capacity expansion plan for all time periods in the first stage, regardless of the parameter realizations, and in the second stage, we decide the capacity allocation plan based on the realized demand in each period. Thus, we do not allow any flexibility in the capacity-acquisition plan with respect to the demand realizations, i.e., $\boldsymbol{x}_n$ should take the same values for all the nodes in the same stage. This gives us the following extensive formulation of the two-stage capacity planning model:
\begin{subequations}
\label{model:tworiskaverse}
\begin{align}
z^{TS}_R=\min_{\substack{\boldsymbol{x}_n,\boldsymbol{y}_n, n\in \mathcal{T}\\\eta_n, n\not\in\mathcal{L}, u_n,n\not=1}}\quad & \sum_{n\in \mathcal{T}} p_n\left(\boldsymbol{\tilde{f}}_n^{\mathsf T}\sum_{m\in \mathcal{P}(n)}\boldsymbol{x}_m+\boldsymbol{\tilde{c}}_n^{\mathsf T}\boldsymbol{y}_n + \tilde{\lambda}_n \eta_n + \tilde{\alpha}_n u_n\right)\\
\text{s.t.}\quad&\text{\eqref{eq:multi-constraint_y}--\eqref{eq:constraint_eta}}\nonumber\\
& {\boldsymbol{x}_m}=\boldsymbol{x}_n,\ \forall m,\ n\in \mathcal{T}_t, \ t=1,...,T\label{eq:twostage_x}\\
& \boldsymbol{x}_n\in \mathbb{Z}^{M}_{+},\ \boldsymbol{y}_n\in \mathbb{R}^{MN \times 1}_{+},\ \forall n\in \mathcal{T},\ u_n \ge 0,\ \forall n\not=1.\nonumber
\end{align}
\end{subequations}

In the following proposition, we first show that the optimal objective values of Models \eqref{model:multiriskaverse} and \eqref{model:tworiskaverse} both increase as the decision maker becomes more risk-averse. All the proofs throughout this paper are presented in Appendix B in the Online Supplement.
\begin{proposition}\label{prop:risk}
    Both $z_R^{MS}$ and $z_R^{TS}$ increase as $\lambda$ increases.
\end{proposition}

We define the difference between the optimal objective values of the risk-averse two-stage and multistage formulations as the \textit{value of risk-averse multistage stochastic programming}: ${\rm VMS_{R}}=z^{TS}_R-z^{MS}_{R}$. From the extensive formulations, we observe that the risk-averse two-stage model \eqref{model:tworiskaverse} is the multistage model \eqref{model:multiriskaverse} with additional constraints \eqref{eq:twostage_x}. We thus conclude that ${\rm VMS_{R}}\ge 0$ as the multistage model provides more flexibility in the capacity acquisition decisions with respect to the uncertain parameter realizations. Because of Constraints \eqref{eq:twostage_x}, $\boldsymbol{x}_n$ variables can be indexed by the time stage $t$ instead and the two-stage model \eqref{model:tworiskaverse} involves $M\times T$ integer variables, while the multistage model \eqref{model:multiriskaverse} involves $M\times |\mathcal{T}|$ integer variables ($|\mathcal{T}|$ is the total number of nodes in the scenario tree $\mathcal{T}$). For any non-trivial scenario tree, $|\mathcal{T}|\gg T$, and thus solving the multistage model \eqref{model:multiriskaverse} requires much more computational effort than solving the two-stage counterpart \eqref{model:tworiskaverse}. If ${\rm VMS_{R}}$ is high, then it may be worth solving a more computationally expensive multistage model; otherwise, a two-stage model will be sufficient to provide a good enough solution.

\section{Value of Multistage {Risk-Averse} Stochastic Programming in Capacity Planning}
	\label{sec:value-proof}
Computing ${\rm VMS_{R}}$ exactly requires us to solve the multistage model, which is computationally expensive. In this section, we provide lower and upper bounds on ${\rm VMS_{R}}$, without the need to solve the multistage model. These bounds can serve as a priori estimates of ${\rm VMS_{R}}$ to analyze the trade-off between risk-averse two-stage and multistage models, and can be used to design guidelines on model choices. Next,
we first examine a substructure of our problem and provide its analytical optimal solutions in Section \ref{sec:substructure}. Based on that, we derive lower and upper bounds for ${\rm VMS_{R}}$ in Sections \ref{sec:VMS} and \ref{sec:VMS-upper}, respectively. Using the derived bounds, we design a flowchart to help decide which model to solve in Section \ref{sec:guidelines}.
\subsection{Analytical Solutions of the Substructure Problem}
\label{sec:substructure}
We examine an important substructure of Models \eqref{model:multiriskaverse} and \eqref{model:tworiskaverse} once we fix $(\boldsymbol{y}, \boldsymbol{u})$-variables. We denote the resultant problems with known $(\boldsymbol{y}_n^*, {u}_n^*)$ values as $\mbox{{\bf SP-RMS}}(\boldsymbol{y}_n^{*}, u_n^{*})$ and $\mbox{{\bf SP-RTS}}(\boldsymbol{y}_n^{*}, u_n^{*})$, which are defined as follows:
\begin{subequations}\label{eq:ms-substructure}
\begin{align}
\mbox{{\bf SP-RMS}}(\boldsymbol{y}_n^{*}, u_n^{*}): \ \min_{\substack{\boldsymbol{x}_n, n\in \mathcal{T}\\\eta_n,,n\not\in\mathcal{L}}} \quad& \sum_{n\in \mathcal{T}} p_n\left(\boldsymbol{\tilde{f}}_n^{\mathsf T}\sum_{m\in \mathcal{P}(n)}\boldsymbol{x}_m+ \tilde{\lambda}_n \eta_n\right)\\
\text{s.t.}\quad& \sum_{m\in \mathcal{P}(n)}\boldsymbol{x}_m\ge \boldsymbol{B}_{t_n}\boldsymbol{y}_n^*,\ \forall n\in \mathcal{T}\label{eq:constraint_xy_knowny}\\
& \eta_{a(n)} \ge \boldsymbol{f}_{t_n}^{\mathsf T} \sum_{m\in \mathcal{P}(n)}\boldsymbol{x}_m+\boldsymbol{c}_{t_n}^{\mathsf T} \boldsymbol{y}_n^* - u^{*}_n,\ \forall n\not= 1 \label{eq:constraint_eta_knowny}\\
&\boldsymbol{x}_n\in \mathbb{Z}_{+}^{M},\ \forall n\in\mathcal{T},\nonumber
\end{align}
\end{subequations}
and
\begin{align}
\mbox{{\bf SP-RTS}}(\boldsymbol{y}_n^{*}, u_n^{*}): \ \min_{\substack{\boldsymbol{x}_n, n\in \mathcal{T}\\\eta_n,,n\not\in\mathcal{L}}}\quad & \sum_{n\in \mathcal{T}} p_n\left(\boldsymbol{\tilde{f}}_n^{\mathsf T}\sum_{m\in \mathcal{P}(n)}\boldsymbol{x}_m+ \tilde{\lambda}_n \eta_n\right) \label{eq:ts-substructure}\\
\text{s.t.}\quad& \text{\eqref{eq:constraint_xy_knowny}--\eqref{eq:constraint_eta_knowny}}\nonumber\\
&\text{\eqref{eq:twostage_x}\ (Two-stage constraints for $\boldsymbol{x}$)}\nonumber\\
&\boldsymbol{x}_n\in \mathbb{Z}_{+}^{M},\ \forall n\in\mathcal{T}.\nonumber
\end{align}
Here, we denote the optimal objective values of Models \eqref{eq:ms-substructure} and \eqref{eq:ts-substructure} as $Q^M(\boldsymbol{y}_n^{*}, u_n^{*}),\ Q^{T}(\boldsymbol{y}_n^{*}, u_n^{*})$, respectively. The next proposition demonstrates the analytical forms of the optimal solutions to $\mbox{{\bf SP-RMS}}(\boldsymbol{y}_n^{*}, u_n^{*})$ and $\mbox{{\bf SP-RTS}}(\boldsymbol{y}_n^{*}, u_n^{*})$, and we present a detailed proof in Appendix B in the Online Supplement.
\begin{proposition}
\label{prop:substructure}
 Given any $(\boldsymbol{y}^*, \boldsymbol{u}^*)$ values, the optimal solutions of \eqref{eq:ms-substructure} and \eqref{eq:ts-substructure} have the following analytical forms:
	\begin{align*}
	&\boldsymbol{x}^{MS}_{1}=\lceil \boldsymbol{B}_{t_1} \boldsymbol{y}_1^{*}\rceil,\ \boldsymbol{x}^{MS}_{n}=\max_{m\in\mathcal{P}(n)}\lceil \boldsymbol{B}_{t_m} \boldsymbol{y}_m^{*}\rceil-\max_{m\in\mathcal{P}(a(n))}\lceil \boldsymbol{B}_{t_m} \boldsymbol{y}_m^{*}\rceil,\ \forall n\not=1,  \\
	&{\eta}^{MS}_n= \max_{m\in \mathcal{C}(n)}\left\lbrace \boldsymbol{f}_{t_m}^{\mathsf T}\sum_{l\in\mathcal{P}(m)}\boldsymbol{x}^{MS}_l+ \boldsymbol{c}_{t_m}^{\mathsf T} \boldsymbol{y}_m^{*}-u_m^{*}\right\rbrace,\ \forall n\not\in\mathcal{L},\\
	&\boldsymbol{x}^{TS}_1 = \lceil \boldsymbol{B}_{t_1} \boldsymbol{y}_1^{*}\rceil,\ \boldsymbol{x}^{TS}_{n}=\max_{m\in\mathcal{P}(n)}\lceil \max_{l\in \mathcal{T}_{t_m}}\boldsymbol{B}_{t_l} \boldsymbol{y}_l^{*}\rceil-\max_{m\in\mathcal{P}(a(n))}\lceil \max_{l\in \mathcal{T}_{t_m}}\boldsymbol{B}_{t_l} \boldsymbol{y}_l^{*}\rceil,\ \forall n\not=1,\\ 
	& \eta^{TS}_n = \max_{m\in \mathcal{C}(n)}\left\lbrace \boldsymbol{f}_{t_m}^{\mathsf T}\sum_{l\in\mathcal{P}(m)}\boldsymbol{x}^{TS}_l+ \boldsymbol{c}_{t_m}^{\mathsf T} \boldsymbol{y}^*_m-u_m^{*}\right\rbrace, \ \forall n\not\in\mathcal{L},
	\end{align*}
	and correspondingly we have $Q^M(\boldsymbol{y}_n^{*}, u_n^{*}) = \sum_{n\in \mathcal{T}} p_n\left(\boldsymbol{\tilde{f}}_n^{\mathsf T}\sum_{m\in \mathcal{P}(n)}\boldsymbol{x}^{MS}_m+ \tilde{\lambda}_n \eta^{MS}_n\right), \	Q^{T}(\boldsymbol{y}_n^{*},u_n^{*}) = \sum_{n\in \mathcal{T}} p_n\left(\boldsymbol{\tilde{f}}_n^{\mathsf T}\sum_{m\in \mathcal{P}(n)}\boldsymbol{x}^{TS}_m + \tilde{\lambda}_n \eta^{TS}_n\right)$.
\end{proposition}
\begin{remark}
From Proposition \ref{prop:substructure}, one can easily verify that $Q^T(\boldsymbol{y}_n^{*}, u_n^{*}) - Q^M(\boldsymbol{y}_n^{*}, u_n^{*})\ge 0$ as $\sum_{m\in \mathcal{P}(n)}\boldsymbol{x}^{TS}_m = \max_{m\in\mathcal{P}(n)}\lceil \max_{l\in \mathcal{T}_{t_m}}\boldsymbol{B}_{t_l} \boldsymbol{y}_l^{*}\rceil \ge \max_{m\in\mathcal{P}(n)}\lceil \boldsymbol{B}_{t_m} \boldsymbol{y}_m^{*}\rceil = \sum_{m\in \mathcal{P}(n)}\boldsymbol{x}^{MS}_m$. Next, we will use Proposition \ref{prop:substructure} to construct lower and upper bounds on the ${\rm VMS_R}$ in Theorem \ref{thm:risk-averse-bound}, Corollary \ref{cor:LP}, and Theorem \ref{thm:risk-averse-upper-bound}, as well as to design approximation algorithms to solve Model \eqref{model:multiriskaverse} in Algorithm \ref{alg:approx-multistage}.
\end{remark}

\subsection{Lower Bound on ${\rm VMS_R}$ for Risk-Averse Capacity Planning}\label{sec:VMS}
We now describe two lower bounds on the ${\rm VMS_R}$ for the risk-averse multistage and two-stage capacity planning models \eqref{model:multiriskaverse} and \eqref{model:tworiskaverse} based on the analysis of the substructure problems. Note that these lower bounds will be useful when they are sufficiently large, indicating that there is a need to solve the much harder multistage model.
\begin{theorem}
\label{thm:risk-averse-bound}
	Let $\{\boldsymbol{y}_n^{*}\}_{n\in \mathcal{T}}, \ \{u_n^{*}\}_{n\in \mathcal{T}\setminus\{1\}}$ be the decisions in an optimal solution to the two-stage model \eqref{model:tworiskaverse}, and let 
$\boldsymbol{x}^{MS},\boldsymbol{\eta}^{MS},\boldsymbol{x}^{TS},\boldsymbol{\eta}^{TS}$ follow the definitions in Proposition \ref{prop:substructure}, which are constructed by $\{\boldsymbol{y}_n^{*}\}_{n\in \mathcal{T}}, \ \{u_n^{*}\}_{n\in \mathcal{T}\setminus\{1\}}$.
	Then,
 \small
	\begin{align}
	{\rm VMS_R}&\ge {\rm VMS_R^{LB}}\nonumber\\
 &:=  \sum_{n\in \mathcal{T}\setminus\{1\}} p_n{(1-\lambda_{t_n})\boldsymbol{f}_{t_n}^{\mathsf T}}\left(\max_{m\in \mathcal{P}(n)} \lceil \max_{l\in \mathcal{T}_{t_m}}\boldsymbol{B}_{t_l} \boldsymbol{y}_l^{*}\rceil - \max_{m\in \mathcal{P}(n)}\lceil \boldsymbol{B}_{t_m} \boldsymbol{y}_m^{*}\rceil\right)+ {\sum_{n\in \mathcal{T}\setminus\mathcal{L}} p_n {\lambda}_{t_n+1}} \left(\eta^{TS}_n-\eta^{MS}_n\right). \label{eq:LB}
	\end{align}
 \normalsize
\end{theorem}

Later in Example \ref{ex:tight}, we show that ${\rm VMS_R^{LB}}$ (Eq. \eqref{eq:LB}) is tight. Next, we derive another lower bound ${\rm VMS_R^{LB1}}$ that is not necessarily tight but more computationally tractable, which utilizes the optimal solutions to the linear programming (LP) relaxation of the two-stage model \eqref{model:tworiskaverse}. Since we consider the LP relaxation of the two-stage model, we do not round up when constructing the two-stage solutions $\boldsymbol{x}_n^{TS}$.
\begin{corollary}\label{cor:LP}
Let $\{\boldsymbol{y}_n^{TSLP}\}_{n\in \mathcal{T}}, \ \{u_n^{TSLP}\}_{n\in \mathcal{T}\setminus\{1\}}$ be the decisions in an optimal solution to the LP relaxation of the two-stage model \eqref{model:tworiskaverse}, and let
	\begin{align*}
	&\boldsymbol{x}^{MS}_{1}=\lceil\boldsymbol{B}_{t_1} \boldsymbol{y}_1^{TSLP}\rceil,\ \boldsymbol{x}^{MS}_{n}=\max_{m\in\mathcal{P}(n)} \lceil\boldsymbol{B}_{t_m} \boldsymbol{y}_m^{TSLP}\rceil-\max_{m\in\mathcal{P}(a(n))} \lceil\boldsymbol{B}_{t_m} \boldsymbol{y}_m^{TSLP}\rceil,\ \forall n\not=1,  \\
	&{\eta}^{MS}_n= \max_{m\in \mathcal{C}(n)}\left\lbrace \boldsymbol{f}_{t_m}^{\mathsf T}\sum_{l\in\mathcal{P}(m)}\boldsymbol{x}^{MS}_l+ \boldsymbol{c}_{t_m}^{\mathsf T} \boldsymbol{y}_m^{TSLP}-u_m^{TSLP}\right\rbrace,\ \forall n\not\in\mathcal{L},\\
	&\boldsymbol{x}^{TS}_1 =  \boldsymbol{B}_{t_1} \boldsymbol{y}_1^{TSLP},\ \boldsymbol{x}^{TS}_{n}=\max_{m\in\mathcal{P}(n)} \max_{l\in \mathcal{T}_{t_m}}\boldsymbol{B}_{t_l} \boldsymbol{y}_l^{TSLP}-\max_{m\in\mathcal{P}(a(n))} \max_{l\in \mathcal{T}_{t_m}}\boldsymbol{B}_{t_l} \boldsymbol{y}_l^{TSLP},\ \forall n\not=1,\\ 
	& \eta^{TS}_n = \max_{m\in \mathcal{C}(n)}\left\lbrace \boldsymbol{f}_{t_m}^{\mathsf T}\sum_{l\in\mathcal{P}(m)}\boldsymbol{x}^{TS}_l+ \boldsymbol{c}_{t_m}^{\mathsf T} \boldsymbol{y}^{TSLP}_m-u_m^{TSLP}\right\rbrace, \ \forall n\not\in\mathcal{L}.
	\end{align*}
Then,
\small
	\begin{align*}
 \hspace{-0.3cm}
	{\rm VMS_R}\ge& {\rm VMS_R^{LB1}}\\
 :=&\boldsymbol{f}_{t_1}^{\mathsf T}\left(\boldsymbol{B}_{t_1} \boldsymbol{y}_1^{TSLP}-\lceil\boldsymbol{B}_{t_1} \boldsymbol{y}_1^{TSLP}\rceil\right)+\sum_{n\in \mathcal{T}\setminus\{1\}} p_n(1-\lambda_{t_n})\boldsymbol{{f}}_{t_n}^{\mathsf T}\left(\max_{m\in \mathcal{P}(n)}  \max_{l\in \mathcal{T}_{t_m}}\boldsymbol{B}_{t_l} \boldsymbol{y}_l^{TSLP} - \max_{m\in \mathcal{P}(n)} \lceil\boldsymbol{B}_{t_m} \boldsymbol{y}_m^{TSLP}\rceil\right) \\
	&+\sum_{n\in \mathcal{T}\setminus\mathcal{L}} p_n {\lambda}_{t_n+1} \left(\eta^{TS}_n-\eta^{MS}_n\right).
	\end{align*}
	\normalsize
\end{corollary}
Note that different than ${\rm VMS_R^{LB}}$, ${\rm VMS_R^{LB1}}$ is not guaranteed to be always non-negative. When ${\rm VMS_R^{LB1}}<0$, we can simply replace it with a trivial lower bound 0, i.e., ${\rm VMS_R}\ge \max\{{\rm VMS_R^{LB1}}, 0\}$. Computing ${\rm VMS_R^{LB}}$ and ${\rm VMS_R^{LB1}}$ requires solving the two-stage model \eqref{model:tworiskaverse} and its LP relaxation, respectively, both of which are more computationally tractable than solving the multistage counterpart.

\subsection{Upper Bound on ${\rm VMS_R}$ for Risk-Averse Capacity Planning}\label{sec:VMS-upper}
We now describe an upper bound on the ${\rm VMS_R}$ based on the optimal solution of the LP relaxation of the multistage model \eqref{model:multiriskaverse}. Since we consider the LP relaxation of the multistage model, we do not round up when constructing the multistage solutions $\boldsymbol{x}_n^{MS}$. Note that this upper bound will be useful when it is small enough, indicating that the two-stage model is sufficient to provide a good enough solution.

\begin{theorem}
\label{thm:risk-averse-upper-bound}
	Let $\{\boldsymbol{y}_n^{MSLP}\}_{n\in \mathcal{T}}, \ \{u_n^{MSLP}\}_{n\in \mathcal{T}\setminus\{1\}}$ be the decisions in an optimal solution to the LP relaxation of multistage model \eqref{model:multiriskaverse}, and let 
	\begin{align*}
	&\boldsymbol{x}^{MS}_{1}=\boldsymbol{B}_{t_1} \boldsymbol{y}_1^{MSLP},\ \boldsymbol{x}^{MS}_{n}=\max_{m\in\mathcal{P}(n)} \boldsymbol{B}_{t_m} \boldsymbol{y}_m^{MSLP}-\max_{m\in\mathcal{P}(a(n))} \boldsymbol{B}_{t_m} \boldsymbol{y}_m^{MSLP},\ \forall n\not=1,  \\
	&{\eta}^{MS}_n= \max_{m\in \mathcal{C}(n)}\left\lbrace \boldsymbol{f}_{t_m}^{\mathsf T}\sum_{l\in\mathcal{P}(m)}\boldsymbol{x}^{MS}_l+ \boldsymbol{c}_{t_m}^{\mathsf T} \boldsymbol{y}_m^{MSLP}-u_m^{MSLP}\right\rbrace,\ \forall n\not\in\mathcal{L},\\
	&\boldsymbol{x}^{TS}_1 =  \lceil\boldsymbol{B}_{t_1} \boldsymbol{y}_1^{MSLP}\rceil,\ \boldsymbol{x}^{TS}_{n}=\max_{m\in\mathcal{P}(n)} \lceil\max_{l\in \mathcal{T}_{t_m}}\boldsymbol{B}_{t_l} \boldsymbol{y}_l^{MSLP}\rceil-\max_{m\in\mathcal{P}(a(n))} \lceil\max_{l\in \mathcal{T}_{t_m}}\boldsymbol{B}_{t_l} \boldsymbol{y}_l^{MSLP}\rceil,\ \forall n\not=1,\\ 
	& \eta^{TS}_n = \max_{m\in \mathcal{C}(n)}\left\lbrace \boldsymbol{f}_{t_m}^{\mathsf T}\sum_{l\in\mathcal{P}(m)}\boldsymbol{x}^{TS}_l+ \boldsymbol{c}_{t_m}^{\mathsf T} \boldsymbol{y}^{MSLP}_m-u_m^{MSLP}\right\rbrace, \ \forall n\not\in\mathcal{L}.
	\end{align*}
	Then,
	\begin{align}
 \hspace{-0.3cm}
	{\rm VMS_R}\le& {\rm VMS_R^{UB}}
 :=\boldsymbol{f}_{t_1}^{\mathsf T}\left(\lceil\boldsymbol{B}_{t_1} \boldsymbol{y}_1^{MSLP}\rceil-\boldsymbol{B}_{t_1} \boldsymbol{y}_1^{MSLP}\right)\nonumber\\
 &+\sum_{n\in \mathcal{T}\setminus\{1\}} p_n{(1-\lambda_{t_n})\boldsymbol{f}_{t_n}^{\mathsf T}}\left(\max_{m\in \mathcal{P}(n)} \lceil \max_{l\in \mathcal{T}_{t_m}}\boldsymbol{B}_{t_l} \boldsymbol{y}_l^{MSLP}\rceil - \max_{m\in \mathcal{P}(n)} \boldsymbol{B}_{t_m} \boldsymbol{y}_m^{MSLP}\right)\nonumber\\
 &+ {\sum_{n\in \mathcal{T}\setminus\mathcal{L}} p_n {\lambda}_{t_n+1}} \left(\eta^{TS}_n-\eta^{MS}_n\right). \label{eq:UB}
	\end{align} 
\end{theorem}

The following example illustrates the usefulness of the derived lower and upper bounds.
\begin{example}\label{ex:tight}
We consider a special example with one facility ($h_1=50$) and one customer site (see Figure \ref{fig:eg}). The underlying scenario tree has two periods and the root node has two branches with equal probability. Let the customer demand be $d_1=0$ in node 1 (stage 1), $d_2=50$ in node 2 (stage 2), and $d_3=150$ in node 3 (stage 2). Denote $f$ as the unit maintenance cost for one stage and $c$ as the operational cost, respectively. Solving a two-stage model yields a total cost of $z_R^{TS}=3f+100c+50\lambda c$, and a multistage model gives a total cost of $z_R^{MS}=2f+100c+\lambda f+50\lambda c$. Thus ${\rm VMS_R}=z_R^{TS}-z_R^{MS}=(1-\lambda)f$. According to Eq. \eqref{eq:LB} and \eqref{eq:UB}, we have ${\rm VMS_R^{LB}}=(1-\lambda)f,\ {\rm VMS_R^{UB}}=(1+\lambda)f$, where the lower bound ${\rm VMS_R^{LB}}$ is actually tight. We refer interested readers to Appendix C in the Online Supplement for detailed computations of these costs and bounds.  
\begin{figure}[ht!]
    \centering
    \includegraphics[width=0.53\textwidth]{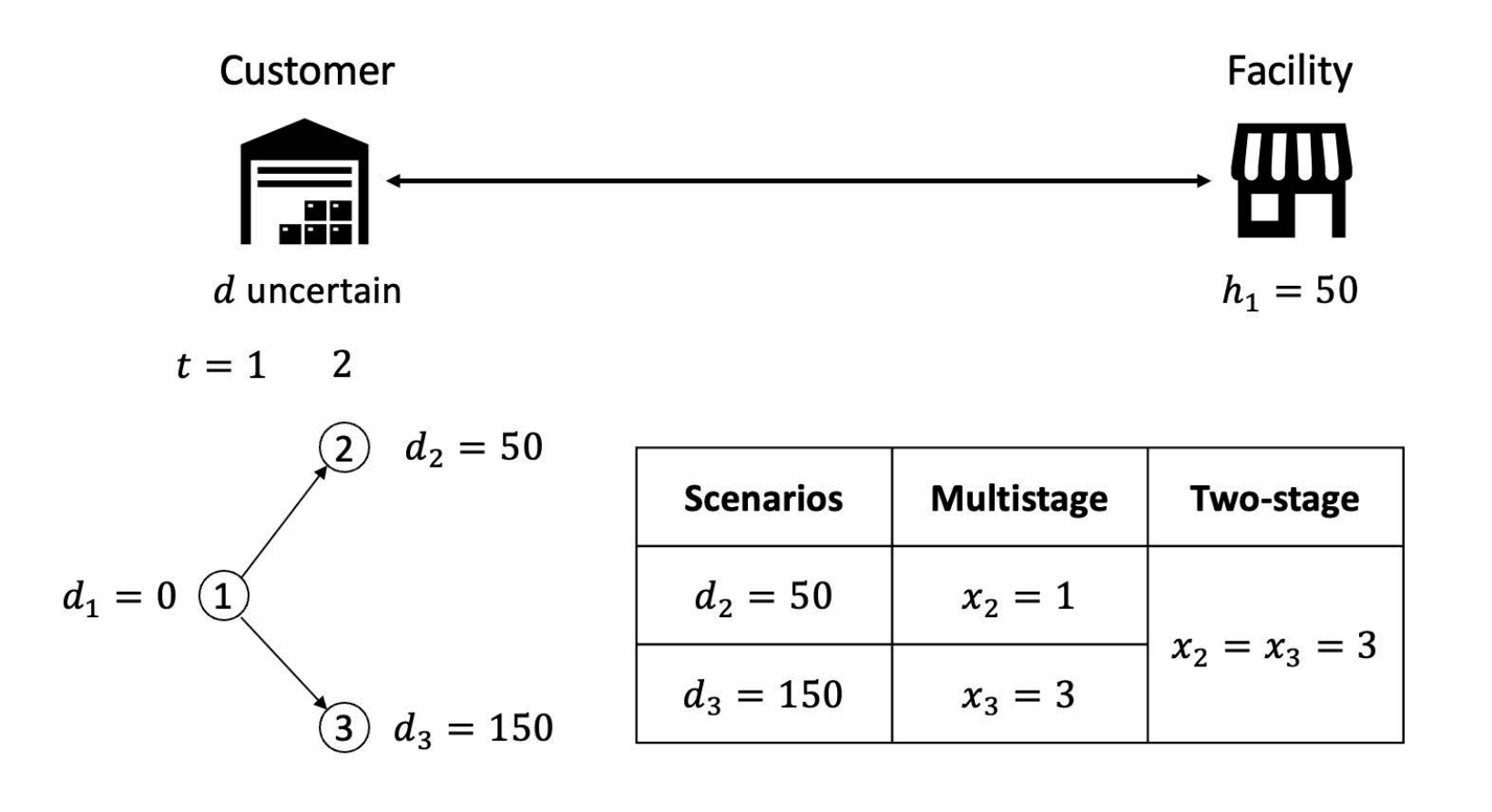}
    \caption{An instance to illustrate the gap between multistage and two-stage capacity planning models.}
    \label{fig:eg}
\end{figure}
\end{example}

\subsection{General Guidelines on Model Choices}\label{sec:guidelines}
\begin{figure}[ht!]
    \centering
    \includegraphics[width=0.9\textwidth]{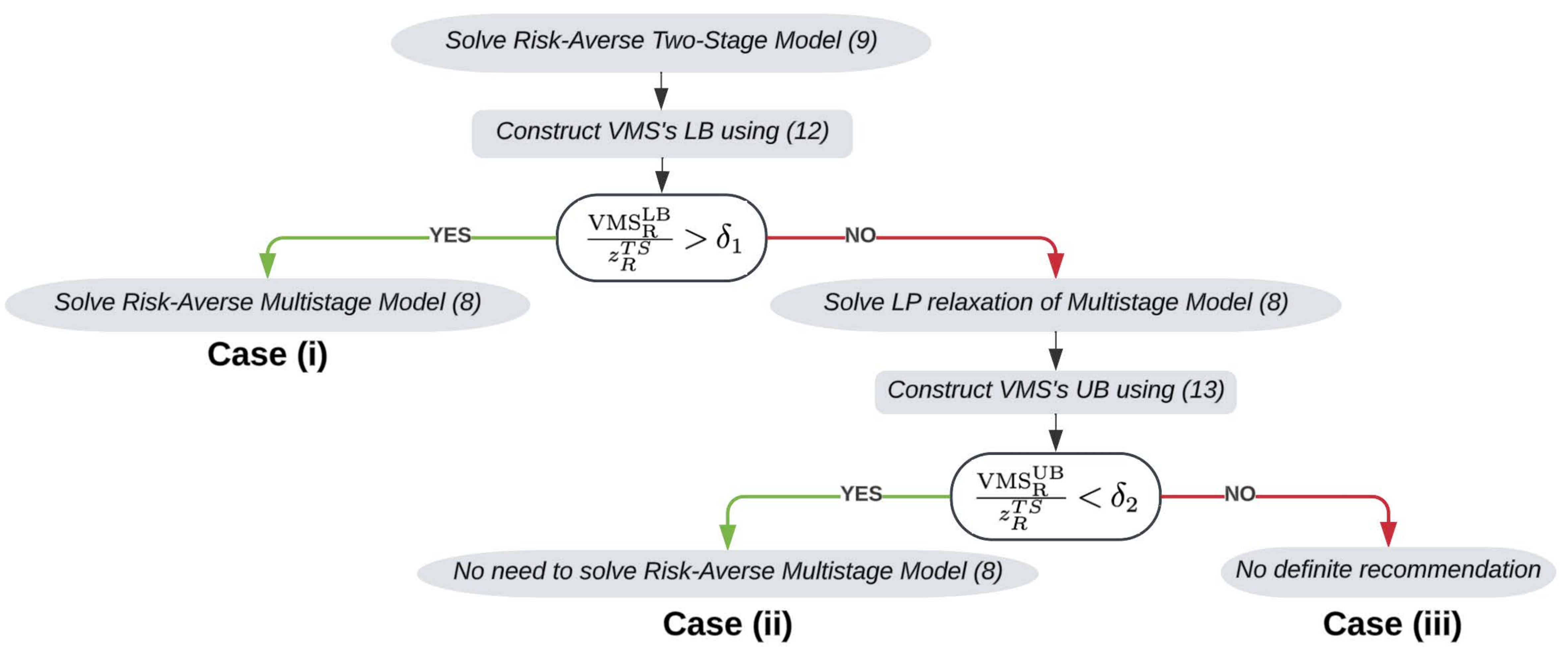}
    \caption{A flowchart to decide which model to solve.}
    \label{fig:flowchart}
\end{figure}

The bounds in Theorems \ref{thm:risk-averse-bound} and \ref{thm:risk-averse-upper-bound} are dependent on the input data and optimal solutions to the two-stage model and LP relaxation of the multistage model. Once the decision maker computes these bounds, they can compare the relative bounds with some pre-given thresholds $\delta_1$ and $\delta_2$ (e.g., $\delta_1=10\%,\ \delta_2=30\%$) and then decide if the two-stage optimal solution is good enough. We outline a general procedure in Figure \ref{fig:flowchart}. 
If the relative ${\rm VMS_R^{LB}}$ is sufficiently large (i.e.,  $\frac{\rm VMS_R^{LB}}{z^{TS}_R}>\delta_1$), then it indicates that the multistage model can reduce the optimal objective value by at least $\delta_1$ compared to the two-stage model and thus one should solve the multistage model. On the other hand, if the relative ${\rm VMS_R^{UB}}$ is sufficiently small (i.e., $\frac{\rm VMS_R^{UB}}{z^{TS}_R}<\delta_2$), then we suggest adopting the optimal solutions to the two-stage model without bearing additional computational effort of solving the multistage model. In the third case, if the lower bound ${\rm VMS_R^{LB}}$ is not large enough and the upper bound ${\rm VMS_R^{UB}}$ is not small enough, then there is no definite recommendation. Later in Section 5, we will show that most of our instances will fall into either Case (i) or Case (ii).  
Note that these two parameters $\delta_1$ and $\delta_2$ are user-defined and can represent the decision maker's trade-off between optimality and tractability.

\begin{example}
    Using the same setting in Example \ref{ex:tight} and letting $\lambda=0.5,\ \delta_1=10\%,\ \delta_2=30\%,\ f=1000,\ c=10$, we have $z_R^{TS}=4250$ and ${\rm VMS_R^{LB}}=500$. Since $\frac{\rm VMS_R^{LB}}{z^{TS}_R}=\frac{500}{4250}>\delta_1=10\%$, we will solve the risk-averse multistage model. In this case, the unit capacity expansion cost $f$ is relatively high and the true ${\rm VMS_R}=(1-\lambda)f$ is also high, indicating the value of the multistage model. On the other hand, if the unit capacity expansion cost is relatively low, e.g., $f=100$, then $z_R^{TS}=1550$ and ${\rm VMS_R^{LB}}=50,\ {\rm VMS_R^{UB}}=150$. Since $\frac{\rm VMS_R^{LB}}{z^{TS}_R}=\frac{50}{1550}<\delta_1=10\%$ and $\frac{\rm VMS_R^{UB}}{z^{TS}_R}=\frac{150}{1550}<\delta_2=30\%$, we will use the optimal solution of the two-stage model without solving the multistage model.
\end{example}

\section{Approximation Algorithms}
\label{sec:approx}
Up to this point, we have discussed how to bound the gap between risk-averse two-stage and multistage models based on their optimal objective values. However, when we enter Case (i), it is inherently hard to solve a risk-averse multistage model as we show in Theorem \ref{thm:hardness}. In this section, we propose approximation algorithms to solve the risk-averse multistage model \eqref{model:multiriskaverse} more efficiently, which can be also applied to solve risk-averse two-stage models \eqref{model:tworiskaverse}. 

\begin{theorem}
\label{thm:hardness}
    The deterministic capacity planning problem \eqref{model:deter} and its risk-averse multistage and two-stage counterparts \eqref{model:multiriskaverse} and \eqref{model:tworiskaverse} are NP-hard. 
\end{theorem}

Motivated by the computational intractability of Models \eqref{model:multiriskaverse} and \eqref{model:tworiskaverse}, we proceed to introduce approximation algorithms that can solve the risk-averse multistage and two-stage models efficiently by utilizing the decomposition structure we investigate in Section \ref{sec:substructure}. Next, we describe the main idea of the algorithm for the risk-averse multistage model \eqref{model:multiriskaverse} as follows: we first solve the LP relaxation of Model \eqref{model:multiriskaverse} to obtain a feasible solution $(\boldsymbol{y}_n^{MSLP},{u}_n^{MSLP})$, which is fed into the substructure problem $\mbox{{\bf SP-RMS}}(\boldsymbol{y}_n^{MSLP}, u_n^{MSLP})$ to obtain an optimal solution $(\boldsymbol{x}_n^1, \eta_n^1)$. Note that we can always find a feasible solution after plugging $(\boldsymbol{y}_n^{MSLP},{u}_n^{MSLP})$ into $\mbox{{\bf SP-RMS}}$ because we can make $\boldsymbol{x}$ and $\eta$ sufficiently large to satisfy Constraints \eqref{eq:constraint_xy_knowny}--\eqref{eq:constraint_eta_knowny}. We then solve Model \eqref{model:multiriskaverse} with fixed $(\boldsymbol{x}_n, \eta_n) = (\boldsymbol{x}_n^1, \eta_n^1)$ to derive an optimal solution $(\boldsymbol{y}_n^{1},{u}_n^{1})$, which together with $(\boldsymbol{x}_n^1, \eta_n^1)$ constitutes a feasible solution and thus an upper bound to Model \eqref{model:multiriskaverse}. This upper bound can be strengthened iteratively by repeating the process, {and we denote the feasible solution produced at the end of Algorithm \ref{alg:approx-multistage} by $(\boldsymbol{x}_n^{H}, \eta_n^{H}, \boldsymbol{y}_n^{H}, u_n^{H})$}. The detailed steps are described in Algorithm \ref{alg:approx-multistage}. We show the monotonicity of the upper bounds derived in Algorithm \ref{alg:approx-multistage} in Proposition \ref{prop:alg}, and show that the optimality gap of Algorithm \ref{alg:approx-multistage} can be upper bounded in Proposition \ref{prop:ratio}, which will eventually lead to an approximation ratio stated in Theorem \ref{thm:ratio}.

\begin{algorithm}[ht!]
% \OneAndAHalfSpacedXI
\caption{Approximation Algorithm for Risk-Averse Multistage Capacity Planning}
\begin{algorithmic}[1]
\label{alg:approx-multistage}
\STATE Solve the LP relaxation of the risk-averse multistage capacity planning problem \eqref{model:multiriskaverse} and let $(\boldsymbol{x}_n^{MSLP}, \eta_n^{MSLP}, \boldsymbol{y}_n^{MSLP}, u_n^{MSLP})_{n\in\mathcal{T}}$ be an optimal solution. If $\boldsymbol{x}_n^{MSLP}$ is integral for all $n\in\mathcal{T}$, stop and return $(\boldsymbol{x}_n^{MSLP}, \eta_n^{MSLP}, \boldsymbol{y}_n^{MSLP}, u_n^{MSLP})_{n\in\mathcal{T}}$.
\STATE Initialize $k=0$ and $(\boldsymbol{x}_n^0, \eta_n^0, \boldsymbol{y}_n^0, u_n^0)_{n\in\mathcal{T}} = (\boldsymbol{x}_n^{MSLP}, \eta_n^{MSLP}, \boldsymbol{y}_n^{MSLP}, u_n^{MSLP})_{n\in\mathcal{T}}$.
\WHILE{$||\boldsymbol{x}^k - \boldsymbol{x}^{k-1}|| \ge \epsilon,\ ||\eta^k - \eta^{k-1}|| \ge \epsilon,\ ||\boldsymbol{y}^k - \boldsymbol{y}^{k-1}|| \ge \epsilon,\ ||u^k - u^{k-1}|| \ge \epsilon$}
\STATE \label{alg:step4} Solve Problem $\mbox{{\bf SP-RMS}}(\boldsymbol{y}_n^{k}, u_n^{k})$
and let $\boldsymbol{x}_n^{k+1},\ {\eta}_n^{k+1}$ denote the corresponding optimal solutions. We have the analytical form of the optimal solutions as $\boldsymbol{x}^{k+1}_{1}=\lceil \boldsymbol{B}_{t_1} \boldsymbol{y}_1^{k}\rceil,\ \boldsymbol{x}_n^{k+1} = \max_{m\in\mathcal{P}(n)}\lceil \boldsymbol{B}_{t_m} \boldsymbol{y}_m^{k}\rceil-\max_{m\in\mathcal{P}(a(n))}\lceil \boldsymbol{B}_{t_m} \boldsymbol{y}_m^{k}\rceil,\ \forall n\not=1,\ {\eta}^{k+1}_n= \max_{m\in \mathcal{C}(n)}\left\lbrace \boldsymbol{f}_{t_m}^{\mathsf T}\sum_{l\in\mathcal{P}(m)}\boldsymbol{x}^{k+1}_l+ \boldsymbol{c}_{t_m}^{\mathsf T} \boldsymbol{y}_m^{k}-u_m^{k}\right\rbrace,\ \forall n\not\in\mathcal{L}$.

\STATE \label{alg:step5} Solve the following problem for each $n\in\mathcal{T},\ n\not=1$ independently
\begin{subequations}\label{eq:independent}
\begin{align}
		\min_{{\boldsymbol{y}_n}, u_n} \quad&   \tilde{\boldsymbol{c}}_{n}^{\mathsf T}\boldsymbol{y}_n + \tilde{\alpha}_nu_n \nonumber\\
\text{s.t.}	\quad	&    \boldsymbol{B}_{t_n}\boldsymbol{y}_n \le \sum_{m\in \mathcal{P}(n)}{\boldsymbol{x}^{k+1}_{m}},\\
		& \boldsymbol{A}\boldsymbol{y}_n =\boldsymbol{d}_n,\\
		& u_n - \boldsymbol{c}_{t_n}^{\mathsf T}\boldsymbol{y}_n \ge \boldsymbol{f}_{t_n}^{\mathsf T} \sum_{m\in \mathcal{P}(n)}\boldsymbol{x}^{k+1}_{m} - \eta^{k+1}_{a(n)},\label{indepent-u}\\
		& \boldsymbol{y}_n\in\mathbb{R}_{+},\ u_n\ge 0,\nonumber
\end{align}
\end{subequations}
and when $n=1$ we solve Problem \eqref{eq:independent} without the variables $u_n$ and constraints \eqref{indepent-u}.
Let $\boldsymbol{y}_n^{k+1},\ {u}_n^{k+1}$ be the optimal solutions.
\STATE Update $k = k+1$.
\ENDWHILE
\STATE Return $(\boldsymbol{x}_n^H, \eta_n^H, \boldsymbol{y}_n^H, u_n^H)_{n\in\mathcal{T}} := (\boldsymbol{x}_n^{k+1}, \eta_n^{k+1}, \boldsymbol{y}_n^{k+1}, u_n^{k+1})_{n\in\mathcal{T}}$.
\end{algorithmic}
\end{algorithm}
Throughout this section, we use $z^{MS}_{T,R}$ to denote the objective value of the multistage model \eqref{model:multiriskaverse} with the total number of stages being $T$. With some abuse of notation, we use $(\boldsymbol{x}_n^{*}, \eta_n^{*}, \boldsymbol{y}_n^{*}, u_n^{*})$ to denote an optimal solution to the multistage model \eqref{model:multiriskaverse} and use $z_{T,R}^{MS}(\boldsymbol{x}_n, \eta_n, \boldsymbol{y}_n, u_n)$ to denote the objective value of Model \eqref{model:multiriskaverse} with decision variable values $(\boldsymbol{x}_n, \eta_n, \boldsymbol{y}_n, u_n)$.
\begin{proposition}\label{prop:alg}
The objective value at the end of each iteration {in Algorithm \ref{alg:approx-multistage}} provides an upper bound to the optimal objective value of Model \eqref{model:multiriskaverse} and it satisfies
    $z_{T,R}^{MS}(\boldsymbol{x}_n^{*}, \eta_n^{*}, \boldsymbol{y}_n^{*}, u_n^{*})\le z_{T,R}^{MS}(\boldsymbol{x}_n^{k+1}, \eta_n^{k+1}, \boldsymbol{y}_n^{k+1}, u_n^{k+1}) \le z_{T,R}^{MS}(\boldsymbol{x}_n^{k}, \eta_n^{k}, \boldsymbol{y}_n^{k}, u_n^{k}),\ \forall k\ge 1$. 
\end{proposition}

\begin{proposition}\label{prop:ratio}
The optimality gap can be bounded above by a quantity only dependent on the maintenance cost, i.e., $z_{T,R}^{MS}(\boldsymbol{x}_n^{H}, \eta_n^{H}, \boldsymbol{y}_n^{H}, u_n^{H}) - z_{T,R}^{MS}(\boldsymbol{x}_n^{*}, \eta_n^{*}, \boldsymbol{y}_n^{*}, u_n^{*}) \le {\sum_{t=1}^T\sum_{i=1}^M f_{ti}}$.
\end{proposition}

\begin{theorem}\label{thm:ratio}
Algorithm \ref{alg:approx-multistage} has an approximation ratio of $$ {1+\frac{M\sum_{t=1}^Tf_{t,\rm max} }{M_{\rm min}\sum_{t=1}^Tf_{t,\rm min}+\sum_{t=1}^Tc_{t,\rm min}\min_{n\in\mathcal{T}_t}\{\sum_{j=1}^Nd_{n,j}\}}},$$ where $h_{\rm max} = \max_{i=1}^M\{h_{1i}\},\ f_{t, \rm max} = \max_{i=1}^M\{f_{ti}\},\ f_{t,\rm min} = \min_{i=1}^M\{f_{ti}\},\ {c_{t,\rm min}=\min_{i\in[M],j\in[N]}c_{tij}}$ and $M_{\rm min} = \lceil\frac{\sum_{j=1}^N d_{1j}}{h_{\rm max}}\rceil$ measures at least how many units of resources we need to cover the first-stage demand. 
\end{theorem}

\begin{corollary}\label{cor:ratio}
Assume that $f_{t,\rm max}=O(1),\ f_{t,\rm min}=O(1),\ c_{t,\rm min}=O(1),\ \min_{n\in\mathcal{T}_t}\sum_{j=1}^Nd_{n,j}=O(t)$ when $t\to\infty$. Then Algorithm \ref{alg:approx-multistage} is asymptotically optimal, i.e.,
    $\lim_{T\to\infty}\frac{z^{MS}_{T,R}(\boldsymbol{x}_n^{H}, \eta_n^{H}, \boldsymbol{y}_n^{H}, u_n^{H})}{z^{MS}_{T,R}(\boldsymbol{x}_n^{*}, \eta_n^{*}, \boldsymbol{y}_n^{*}, u_n^{*})}=1$.
\end{corollary}

Note that the assumptions in Corollary \ref{cor:ratio} are not particularly restrictive. They only require that the unit maintenance and operational costs can be bounded above by a value that does not depend on the time stage, and the demand grows at least linearly when time increases (the result still holds when the demand grows faster than linearly, e.g., when $\min_{n\in\mathcal{T}_t}\sum_{j=1}^Nd_{n,j}=O(t^2)$). The last condition can be achieved if we have an expanding market and the minimum of the demand is increased for each subsequent year (e.g., $\min_{n\in\mathcal{T}_t}\sum_{j=1}^Nd_{n,j}=\tilde{d}(1+2(t-1))$ with $\tilde{d}$ being the nominal demand in the first stage). We will test different demand patterns in Section \ref{sec:real}. Detailed proofs of Propositions \ref{prop:alg}, \ref{prop:ratio}, Theorem \ref{thm:ratio} {and Corollary \ref{cor:ratio}} are given in Appendix B in the Online Supplement.

One can also tailor Algorithm \ref{alg:approx-multistage} to solve the risk-averse two-stage model \eqref{model:tworiskaverse} by modifying Step 4 to use the analytical solutions of $\mbox{{\bf SP-RTS}}(\boldsymbol{y}_n^{k}, u_n^{k})$, as we discussed in Proposition \ref{prop:substructure}. 
We omit the details in the interest of brevity.

\section{Computational Results}
\label{sec:compu}
We test the risk-averse two-stage and multistage models on two types of networks -- a randomly generated grid network where we vary the parameter settings extensively and a real-world network based on the United States map with 49 candidate facilities and 88 customer sites \citep{daskin2011network}.
Specifically, we conduct sensitivity analysis and report results based on the synthetic data to illustrate the tightness of the analytical bounds and the efficacy and efficiency of the proposed approximation algorithms in Section \ref{sec:synthetic}. We also conduct a case study on EV charging station capacity planning based on the United States map to display the solution patterns under different settings of uncertainties in Section \ref{sec:real}.
We use Gurobi 10.0.0 coded in Python 3.11.0 for solving all mixed-integer programming models, where the computational time limit is set to one hour. Our numerical tests are conducted on a Macbook Pro with 8 GB RAM and an Apple M1 Pro chip. The source code and data files can be found in \cite{yu2024repo}.

\subsection{Result Analysis on Synthetic Data}
\label{sec:synthetic}
 We first introduce the experimental design and setup in Section \ref{sec:setup}, and report sensitivity analysis results in Section \ref{sec:sensi}. 
 Then we examine the tightness of the analytical bounds in Section \ref{sec:tight} and the performance of the approximation algorithms in Section \ref{sec:aa}, respectively. {Finally, in Section \ref{sec:time}, we report the computational time for solving the two risk-averse models.} 

\subsubsection{Experimental Design and Setup.}
\label{sec:setup}
We randomly sample $M$ potential facilities and $N$ customer sites on a $100\times 100$ grid and in the default setting, we have the number of stages ($T$) being 3, the number of facilities ($M$) being 5, the number of customer sites ($N$) being 10, the number of branches in each non-leaf node ($C$) being 2. The risk attitude parameters are set to $\lambda_t=0.5,\ \alpha_t=0.95,\ \forall t=2,\ldots, T$ at default. The operational costs between facilities and customer sites are calculated by their Manhattan distances times the unit travel cost.
We set the per stage maintenance costs $f_{ti}=6\times 10^4$ and all the facilities have the same unit capacity $h_{ti}=h=10^3$ for all $t=1,\ldots,T,\ i=1,\ldots,M$. For each customer site $j=1,\ldots, N$ and stage $t=1,\ldots,T$, we uniformly sample the demand mean from $U(1000(2t-1),5000(2t-1))$ and then multiply each mean by a fixed number ($\sigma=0.8$ at default) to generate its demand standard deviation. Lastly, we sample demand data following a truncated Normal distribution with the generated mean and standard deviation, while negative demand values are deleted.  We consider two types of scenario trees: 
\begin{itemize}
    \item Stagewise-dependent (SD): at every stage $t=1,\ldots,T-1$, every node $n\in\mathcal{T}_t$ is associated with a different set of children nodes $\mathcal{C}(n)$, i.e., $\mathcal{C}(n)\not = \mathcal{C}(m),\ \forall n,m\in\mathcal{T}_t$;
    \item Stagewise-independent (SI): at every stage $t=1,\ldots,T-1$, every node $n\in\mathcal{T}_t$ is associated with an identical set of children nodes $\mathcal{C}(n)$, i.e., $\mathcal{C}(n) = \mathcal{C}(m),\ \forall n,m\in\mathcal{T}_t$;
\end{itemize}
Here, SD represents the most general case where in each stage $t$, we have at most $C^{t}$ different realizations of the uncertainty $\boldsymbol{d}_t$. On the other hand, SI assumes that the stochastic process $(\boldsymbol{d}_1,\boldsymbol{d}_2,\ldots,\boldsymbol{d}_T)$ is stagewise independent and thus we have at most $C$ different realizations of the uncertainty in each stage $t$.

\subsubsection{Sensitivity Analysis on ${\rm VMS}$.}
\label{sec:sensi}
\begin{figure}[ht!]
    \centering
    \begin{subfigure}{0.4\textwidth}
 \resizebox{\textwidth}{!}{%
\begin{tikzpicture}
  \begin{axis}
  [
    xlabel={Number of branches $C$},
    ylabel={${\rm RVMS}$},
    xtick={2,3,4,5},
    yticklabel=
{\pgfmathparse{\tick*100}\pgfmathprintnumber{\pgfmathresult}\%},
cycle list name=black white,
   legend pos= outer north east
    % legend cell align={left}
]
\addplot [only marks] 
  plot [error bars/.cd, y dir=both, y explicit]
  table [y error minus=y-min] {\RVMSC};
    \addplot coordinates {
(2,	0.081269906)
(3, 0.120356703)
(4, 0.138398513)
(5, 0.141935673)
    };\pgfplotsset{cycle list shift=5}
    \addplot [only marks] 
  plot [error bars/.cd, y dir=both, y explicit]
  table [y error minus=y-min] {\RVMSCSAA};
    \addplot coordinates {
(2,	0.047249868)
(3, 0.059624061)
(4, 0.079003529)
(5, 0.082601612)
    };
    \legend{SI-${\rm RVMS}$,SD-${\rm RVMS}$}
  \end{axis}
\end{tikzpicture}%
}
 \caption{different numbers of branches $C$}
\end{subfigure}
\begin{subfigure}{0.4\textwidth}
 \resizebox{\textwidth}{!}{%
 \pgfplotsset{scaled y ticks=false}
\begin{tikzpicture}
  \begin{axis}
  [
    xlabel={Number of stages $T$},
    ylabel={${\rm RVMS}$},
    xtick={3,4,5,6,7,8,9,10},
    yticklabel= {\pgfmathparse{\tick*100}\pgfmathprintnumber{\pgfmathresult}\%},
    cycle list name=black white,
     legend pos= outer north east
    ]
    \addplot [only marks] 
  plot [error bars/.cd, y dir=both, y explicit]
  table [y error minus=y-min] {\RVMST};
    \addplot coordinates {
(3,	0.081269906)
(4, 0.121580874)
(5, 0.161268633)
(6, 0.181477098)
(7, 0.206564113)
(8, 0.228041379)
(9, 0.232445938)
(10, 0.249707189)
    };\pgfplotsset{cycle list shift=5}
    \addplot [only marks] 
  plot [error bars/.cd, y dir=both, y explicit]
  table [y error minus=y-min] {\RVMSTSAA};
    \addplot coordinates {
(3,	0.047249868)
(4, 0.043323422)
(5, 0.042672586)
(6, 0.042325944)
(7, 0.044753596)
(8, 0.045210809)
(9, 0.04005351)
(10, 0.062513911)
    };
    \legend{SI-${\rm RVMS}$, SD-${\rm RVMS}$}
  \end{axis}
\end{tikzpicture}%
}
\caption{different numbers of stages $T$}
\end{subfigure}
\begin{subfigure}{0.4\textwidth}
 \resizebox{\textwidth}{!}{%
\begin{tikzpicture}
    \begin{axis}
  [
    xlabel={Risk attitude $\lambda$},
    ylabel={${\rm RVMS}$},
    xtick={0,0.2,0.4,0.6,0.8,1},
    yticklabel=
{\pgfmathparse{\tick*100}\pgfmathprintnumber{\pgfmathresult}\%},
    cycle list name=black white,
    legend pos= outer north east
    ]
    \addplot [only marks] 
  plot [error bars/.cd, y dir=both, y explicit]
  table [y error minus=y-min] {\RVMSL};
    \addplot coordinates {
(0,	0.121421711)
(0.2, 0.103277735)
(0.4, 0.090164587)
(0.6, 0.077791966)
(0.8, 0.058526046)
(1, 0.042120239)
    };\pgfplotsset{cycle list shift=5}
    \addplot [only marks] 
  plot [error bars/.cd, y dir=both, y explicit]
  table [y error minus=y-min] {\RVMSLSAA};
    \addplot coordinates {
(0, 0.081540354)
(0.2, 0.072638059)
(0.4, 0.051665765)
(0.6, 0.036500582)
(0.8, 0.017707293)
(1, 0.00196813)
    };
    \legend{SI-${\rm RVMS}$, SD-${\rm RVMS}$}
  \end{axis}
\end{tikzpicture}%
}
\caption{different risk attitudes $\lambda$}
\end{subfigure}
\begin{subfigure}{0.4\textwidth}
 \resizebox{\textwidth}{!}{%
 \pgfplotsset{scaled y ticks=false}
\begin{tikzpicture}
    \begin{axis}
  [
    xlabel={Demand standard deviation $\sigma$},
    ylabel={${\rm RVMS}$},
    xtick={0.2,0.4,0.6,0.8},
    yticklabel={\pgfmathparse{\tick*100}\pgfmathprintnumber{\pgfmathresult}\%},
    cycle list name=black white,
     legend pos= outer north east
    ]
    \addplot [only marks] 
  plot [error bars/.cd, y dir=both, y explicit]
  table [y error minus=y-min] {\RVMSS};
    \addplot coordinates {
(0.2,	0.033305674)
(0.4, 0.058298967)
(0.6, 0.079716841)
(0.8, 0.081269906)
    };\pgfplotsset{cycle list shift=5}
    \addplot [only marks] 
  plot [error bars/.cd, y dir=both, y explicit]
  table [y error minus=y-min] {\RVMSSSAA};
    \addplot coordinates {
(0.2,	0.015487577)
(0.4,0.028875099)
(0.6, 0.03583293)
(0.8, 0.047249868)
    };
    \legend{SI-${\rm RVMS}$, SD-${\rm RVMS}$}
  \end{axis}
\end{tikzpicture}%
}
\caption{different standard deviations $\sigma$}
\end{subfigure}
    \caption{Statistics of ${\rm RVMS}$ over 100 instances with different numbers of branches $C$, stages $T$, risk attitudes $\lambda$ and demand standard deviations $\sigma$.}
    \label{fig:sensitivity_RVMS}
\end{figure}
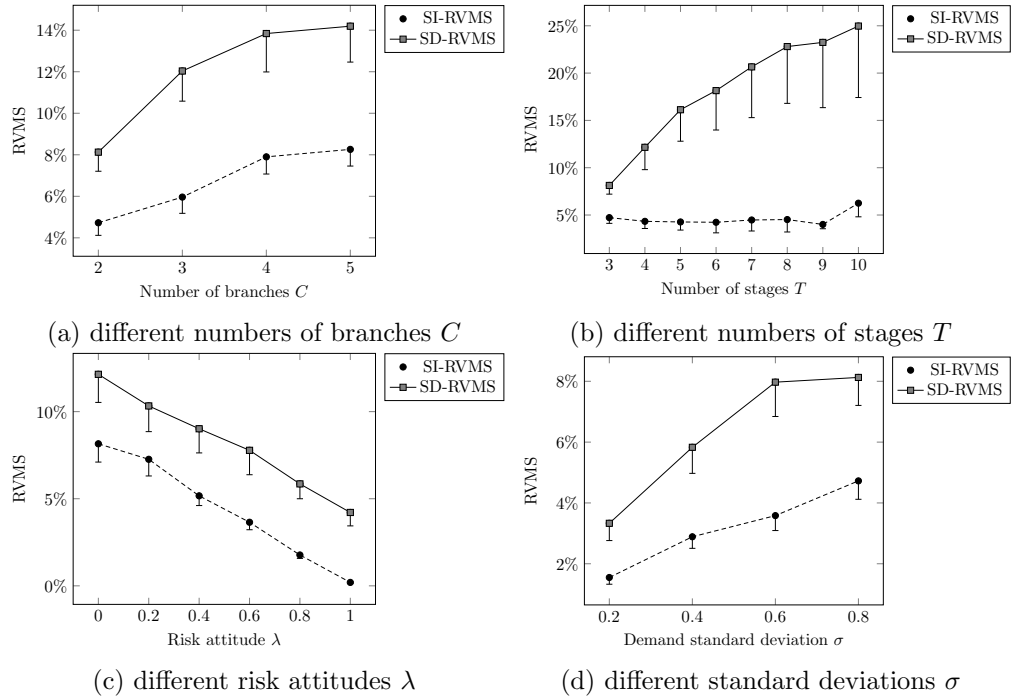

To compare the two-stage and multistage models, we define the relative value of risk-averse multistage stochastic programming as ${\rm RVMS} = \frac{z_R^{TS} - z_R^{MS}}{z_R^{TS}}$ and a lower bound on the RVMS as ${\rm RVMS_{LB}}=\frac{\rm VMS_R^{LB}}{z_R^{TS}}$ using ${\rm VMS_R^{LB}}$ computed by Eq. \eqref{eq:LB}. To illustrate what types of demand scenarios make the multistage model more valuable, we compare the above two types of scenario trees to evaluate ${\rm RVMS}$, namely SD-${\rm RVMS}$, and SI-${\rm RVMS}$, respectively. We first vary the number of branches $C$ from 2 to 5, the number of stages $T$ from 3 to 6, the risk attitude $\lambda$ from 0 to 1, and the standard deviation $\sigma$ from 0.2 to 0.8 to see how ${\rm RVMS}$ changes with respect to different parameter settings. The corresponding results are presented in Figure \ref{fig:sensitivity_RVMS}, where we plot the mean of ${\rm RVMS}$ over 100 independently generated instances, and ${\rm RVMS_{LB}}$ is represented by the lower end of the error bars.

From the figure, when we increase the demand variability such as the number of branches, {the number of stages} and the standard deviation, ${\rm RVMS}$ values with stagewise-dependent scenario trees increase drastically (see Figures \ref{fig:sensitivity_RVMS}(a), (b) and (d)). 
A higher ${\rm RVMS}$ indicates that the multistage model is much more valuable than the two-stage counterpart. 
Comparing different types of scenario trees, SD scenario trees always gain much higher ${\rm RVMS}$ than the stagewise-independent ones, and adding more stages makes no significant changes in SI-RVMS, as can be seen in Figure \ref{fig:sensitivity_RVMS}(b). This is because, in SI scenario trees, the number of realizations in each stage does not depend on the number of stages, and having a deeper scenario tree would not necessarily increase the demand variability. Following these observations, if the stochastic capacity planning problem has a large number of branches, stages, or high demand standard deviation with a stagewise-dependent scenario tree, then it is worth solving a multistage model because the gap between a two-stage formulation and a multi-stage formulation is very high.
Notably, from Figure \ref{fig:sensitivity_RVMS}(c), ${\rm RVMS}$ decreases approximately linearly with respect to the risk attitude parameter $\lambda$, suggesting that as we become more risk-averse, the performances of the two models tend to be closer. This can be understood from the following facts: (i) a two-stage model restricts all capacity-expansion decisions to be identical in the same stage and (ii) a risk-averse multistage model aims to reduce the upper tail of the cost distribution in each stage. As a result, the risk-averse multistage model tends to reduce the dispersion of capacity-expansion decisions in each stage as we become more risk-averse, which has similar effects to two-stage models. Because of this, if a decision maker is extremely risk-averse (e.g., $\lambda$ $\approx$ 1), given the small gaps between the two models, we suggest solving a risk-averse two-stage model without bearing the additional computational effort.  Moreover, the ${\rm RVMS_{LB}}$ obtained by all instances are pretty close to the true RVMS, and stagewise-independent scenario trees obtain tighter lower bounds than stagewise-dependent ones. 
We will take a closer look at these lower bounds in Section \ref{sec:tight}.

Following the flowchart in Figure \ref{fig:flowchart} and setting $\delta_1 = 10\%,\ \delta_2=30\%$, we have the following three cases: Case (i) ${\rm RVMS_{LB}} > 10\%$; Case (ii) ${\rm RVMS_{LB}} \le 10\%$ and ${\rm RVMS_{UB}} < 30\%$; and Case (iii) ${\rm RVMS_{LB}} \le 10\%$ and ${\rm RVMS_{UB}} \ge 30\%$. We record the number of instances in each case with SD scenario trees and varying numbers of branches and stages in Tables \ref{tab:percentage-C}--\ref{tab:percentage-T}. From these tables, in most of the instances, we have either Case (i) or Case (ii). As we increase the number of branches or stages, most cases shift to Case (i) from Case (ii), which indicates a higher need for solving multistage models. We also present the results when we set $\delta_1=5\%,\ \delta_2=20\%$ and change other parameters in Appendix D in the Online Supplement.
\begin{table}[ht!]
% \OneAndAHalfSpacedXI
  \centering
  \caption{Percentage of instances in different cases with SD scenario trees and varying number of branches $C$ when $\delta_1 = 10\%,\ \delta_2=30\%$}
    \begin{tabular}{rrrr}
    \hline
    Branches C & Case (i) & Case (ii) & Case (iii) \\
    \hline
    2     & 20/100 & 73/100 & 7/100 \\
    3     & 56/100 & 25/100 & 19/100 \\
    4     & 69/100 & 8/100 & 23/100 \\
    5     & 79/100 & 2/100 & 19/100 \\
    \hline
    \end{tabular}%
  \label{tab:percentage-C}%
\end{table}%
\begin{table}[ht!]
% \OneAndAHalfSpacedXI
  \centering
  \caption{Percentage of instances in different cases with SD scenario trees varying number of stages $T$ when $\delta_1 = 10\%,\ \delta_2=30\%$}
    \begin{tabular}{rrrr}
    \hline
    Stages T & Case (i) & Case (ii) & Case (iii) \\
    \hline
    3     & 20/100 & 73/100 & 7/100 \\
    4     & 40/100 & 35/100 & 25/100 \\
    5     & 85/100 & 5/100 & 10/100 \\
    6     & 95/100 & 0/100 & 5/100 \\
    7     & 98/100	& 0/100	& 2/100 \\
    8     & 100/100	& 0/100	 & 0/100\\
    \hline
    \end{tabular}%
  \label{tab:percentage-T}%
\end{table}%

\subsubsection{Tightness of the Analytical Bounds on the Synthetic Data Set}
\label{sec:tight}
\begin{figure}
    \centering
    \begin{subfigure}[b]{0.3\textwidth}
        \includegraphics[width=\textwidth]{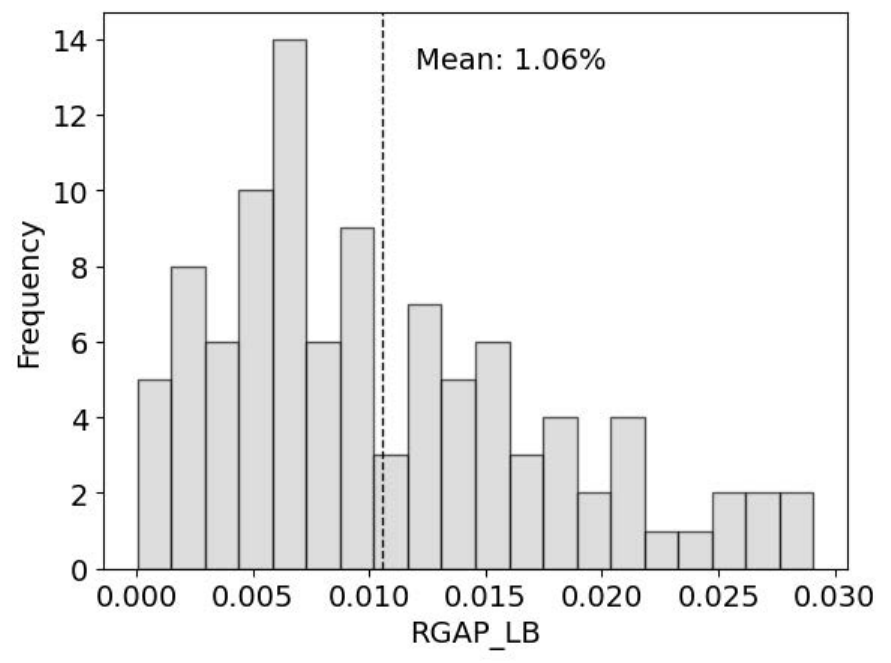}
        \caption{$\rm{RGAP_{LB}}$ with SD trees}
    \end{subfigure}
                \begin{subfigure}[b]{0.3\textwidth}
        \includegraphics[width=\textwidth]{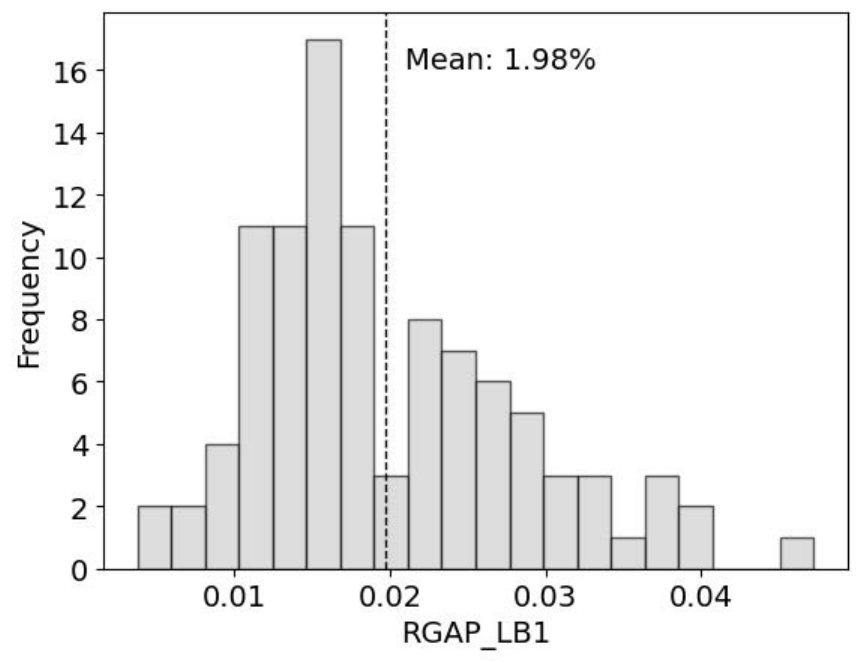}
        \caption{$\rm{RGAP_{LB1}}$ with SD trees}
    \end{subfigure}
                    \begin{subfigure}[b]{0.3\textwidth}
        \includegraphics[width=\textwidth]{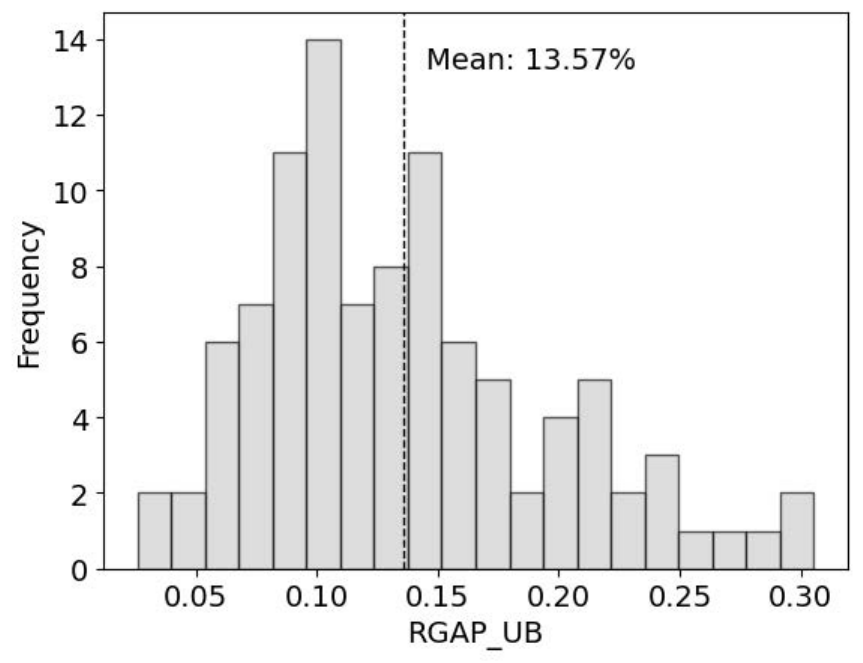}
        \caption{$\rm{RGAP_{UB}}$ with SD trees}
    \end{subfigure}
        \begin{subfigure}[b]{0.3\textwidth}
        \includegraphics[width=\textwidth]{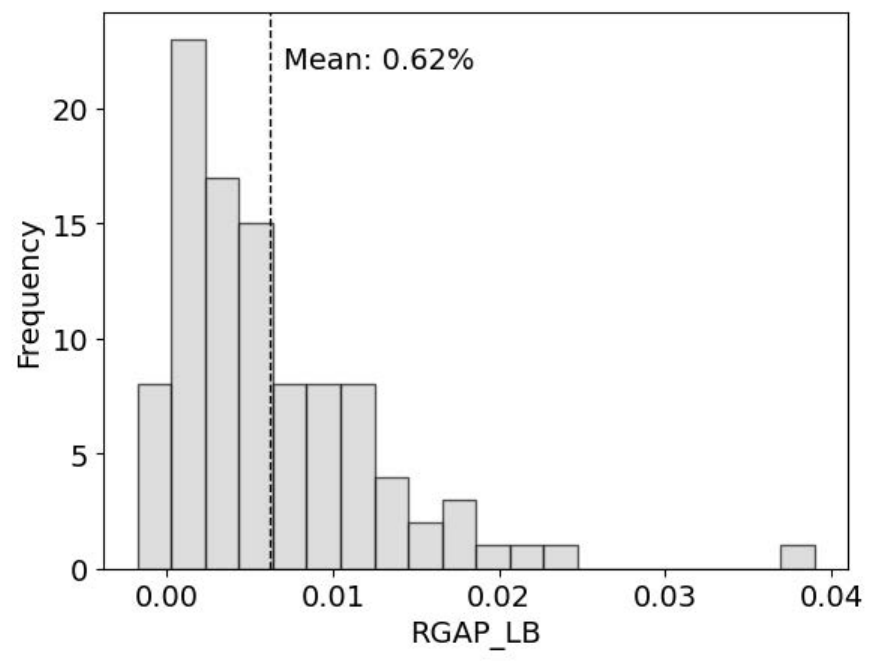}
        \caption{$\rm{RGAP_{LB}}$ with SI trees}
    \end{subfigure}
                \begin{subfigure}[b]{0.3\textwidth}
        \includegraphics[width=\textwidth]{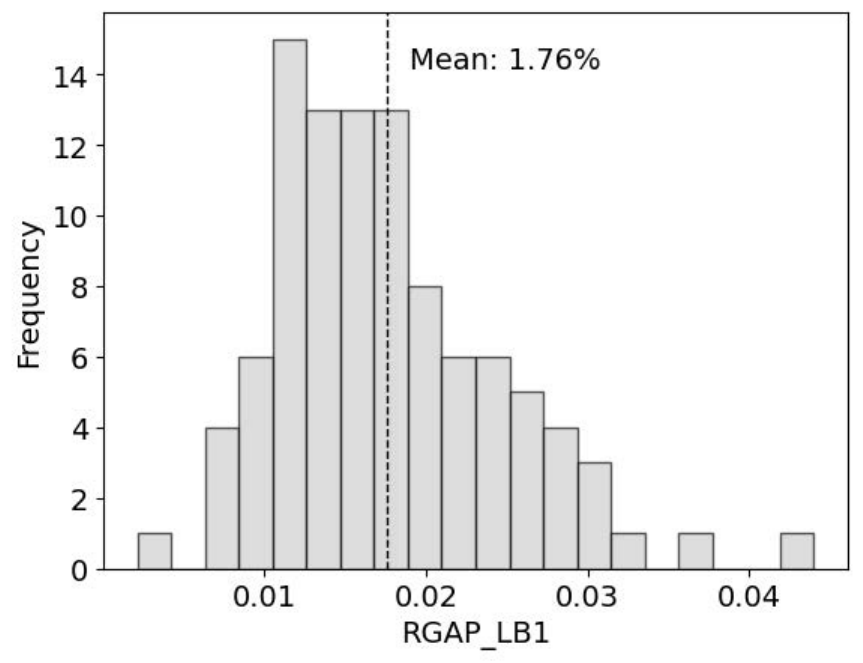}
        \caption{$\rm{RGAP_{LB1}}$ with SI trees}
    \end{subfigure}
                    \begin{subfigure}[b]{0.3\textwidth}
        \includegraphics[width=\textwidth]{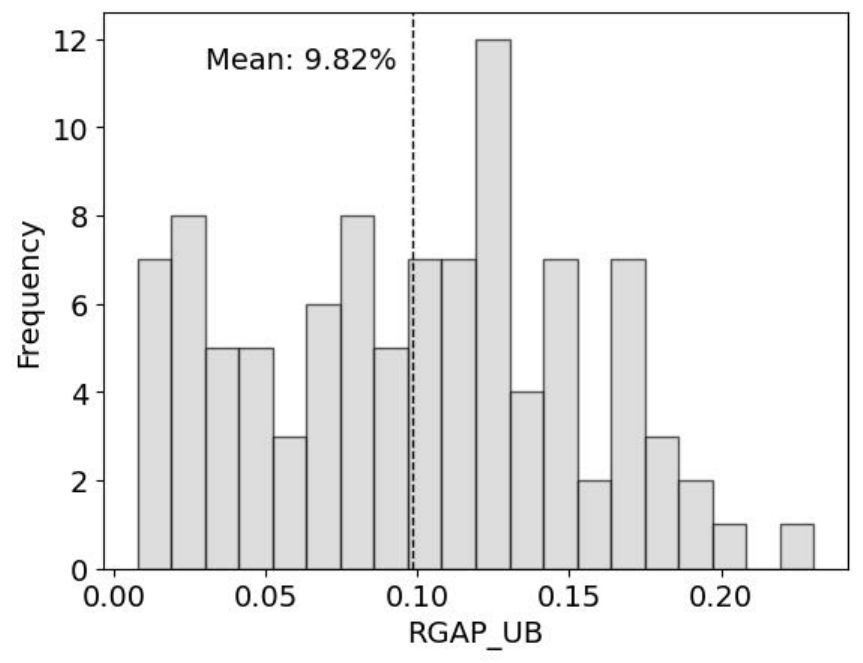}
        \caption{$\rm{RGAP_{UB}}$ with SI trees}
    \end{subfigure}
    \caption{Tightness of analytical bounds with SD and SI scenario trees.}
    \label{fig:tightness}
\end{figure}

To examine the tightness of the analytical bounds derived in Theorem \ref{thm:risk-averse-bound}, Corollary \ref{cor:LP} and Theorem \ref{thm:risk-averse-upper-bound}, we define the relative gaps of the analytical bounds as ${\rm RGAP_{LB}} = {\rm RVMS} - {\rm RVMS_{LB}},\ {\rm RGAP_{LB1}} = {\rm RVMS} - {\rm RVMS_{LB1}},\ {\rm RGAP_{UB}} = {\rm RVMS_{UB}} - {\rm RVMS}$, respectively. A lower ${\rm RGAP}$ means that the analytical bound recovers the true $\rm{VMS}$ better. Using the default setting and the two scenario trees defined in Section \ref{sec:setup}, we plot the histograms of ${\rm RGAP_{LB}},\ {\rm RGAP_{LB1}},\ {\rm RGAP_{UB}}$ over 100 independently generated instances in Figure \ref{fig:tightness}. From the figure, ${\rm RGAP_{LB}}$ with SI scenario trees obtains the lowest ${\rm RGAP}$ mean (0.62\%), where in most instances ${\rm RGAP_{LB}}$ is within 1\%. Stagewise-dependent scenario trees obtain a slightly higher mean of ${\rm RGAP_{LB}}$. This is because RVMS in stagewise-dependent scenario trees is much higher than the stagewise-independent counterparts (see Figure \ref{fig:sensitivity_RVMS}), and as a result, the optimal solutions of the two-stage models are farther away from the optimal ones of the multistage models. As can be seen from the proof of Theorem \ref{thm:risk-averse-bound}, the tightness of the lower bounds depends on the gap between these two optimal solutions. On the other hand, ${\rm VMS_R^{LB1}}$ produces a slightly higher ${\rm RGAP}$, with a mean of 1.98\% in the SD scenario trees and 1.76\% in the SI scenario trees. Compared to the tightness of lower bounds, ${\rm VMS_R^{UB}}$ is much looser, with an RGAP mean of 13.57\% in SD scenario trees and 9.82\% in SI scenario trees.

\subsubsection{Performance of the Approximation Algorithm on the Synthetic Data Set}
\label{sec:aa}
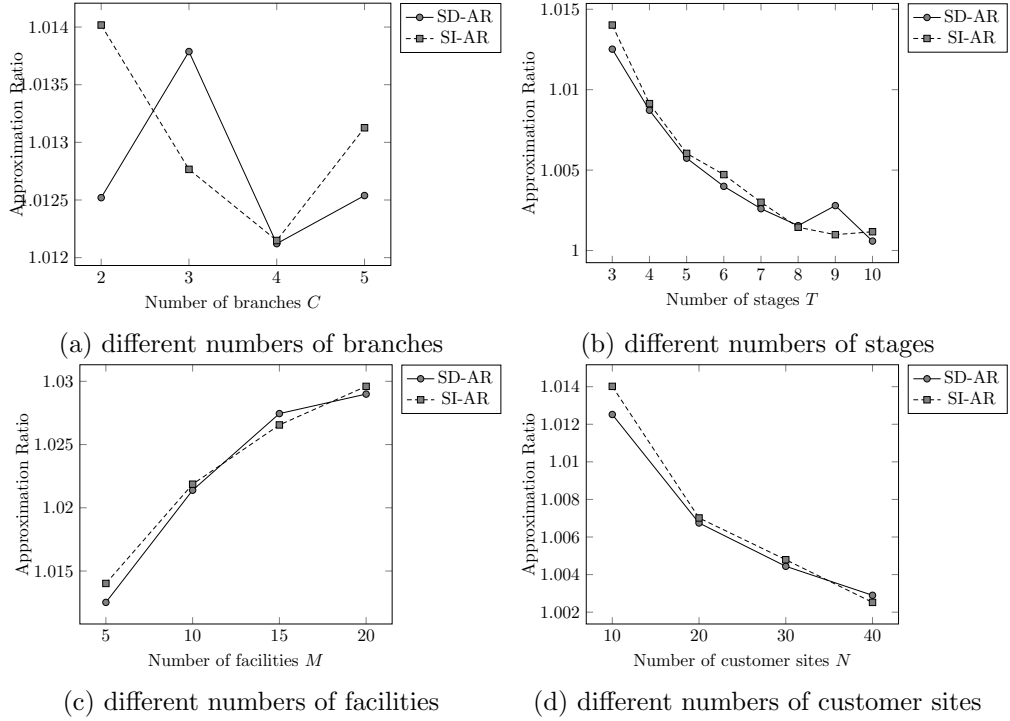
\begin{figure}[ht!]
    \centering
    \begin{subfigure}{0.4\textwidth}
 \resizebox{\textwidth}{!}{%
\begin{tikzpicture}
  \begin{axis}
  [yticklabel style={/pgf/number format/fixed,
                  /pgf/number format/precision=4},
    xlabel={Number of branches $C$},
    ylabel={Approximation Ratio},
    xtick={2,3,4,5},
cycle list name=black white,
legend pos= outer north east
]
    \addplot coordinates {
(2,	1.012519803)
(3, 1.013787096)
(4, 1.012121144)
(5, 1.012538309)
    };\pgfplotsset{cycle list shift=5}
    \addplot coordinates {
(2,	1.014017216)
(3, 1.012765429)
(4, 1.012148279)
(5, 1.013126254)
    };
    \legend{SD-AR,SI-AR}
  \end{axis}
\end{tikzpicture}%
}
\caption{different numbers of branches}
\end{subfigure}
\begin{subfigure}{0.4\textwidth}
 \resizebox{\textwidth}{!}{%
 \pgfplotsset{scaled y ticks=false}
\begin{tikzpicture}
  \begin{axis}
  [yticklabel style={/pgf/number format/fixed,
                  /pgf/number format/precision=3},
    xlabel={Number of stages $T$},
    ylabel={Approximation Ratio},
    xtick={3,4,5,6,7,8,9,10},
    cycle list name=black white,
    legend pos= outer north east
    ]
    \addplot coordinates {
(3,	1.012519803)
(4, 1.008722353)
(5, 1.005741382)
(6, 1.003990685)
(7, 1.002601263)
(8, 1.001546789)
(9, 1.00279318)
(10, 1.000588126)
    };\pgfplotsset{cycle list shift=5}
    \addplot coordinates {
(3,	1.014017216)
(4, 1.009133412)
(5, 1.00603948)
(6, 1.004721164)
(7, 1.003000734)
(8, 1.001448225)
(9, 1.000986006)
(10, 1.001167967)
    };
   \legend{SD-AR,SI-AR}
  \end{axis}
\end{tikzpicture}%
}
\caption{different numbers of stages}
\end{subfigure}
\begin{subfigure}{0.4\textwidth}
 \resizebox{\textwidth}{!}{%
 \pgfplotsset{scaled y ticks=false}
\begin{tikzpicture}
    \begin{axis}
  [yticklabel style={/pgf/number format/fixed,
                  /pgf/number format/precision=3},
    xlabel={Number of facilities $M$},
    ylabel={Approximation Ratio},
    xtick={5,10,15,20},
    cycle list name=black white,
    legend pos= outer north east
    ]
    \addplot coordinates {
(5,	1.012519803)
(10, 1.021390524)
(15, 1.027447951)
(20, 1.028990251)
    };\pgfplotsset{cycle list shift=5}
    \addplot coordinates {
(5,	1.014017216)
(10, 1.021865974)
(15, 1.026550333)
(20, 1.029601179)
    };
    \legend{SD-AR,SI-AR}
  \end{axis}
\end{tikzpicture}%
}
\caption{different numbers of facilities}
\end{subfigure}
\begin{subfigure}{0.4\textwidth}
 \resizebox{\textwidth}{!}{%
 \pgfplotsset{scaled y ticks=false}
\begin{tikzpicture}
    \begin{axis}
  [yticklabel style={/pgf/number format/fixed,
                  /pgf/number format/precision=3},
    xlabel={Number of customer sites $N$},
    ylabel={Approximation Ratio},
    xtick={10,20,30,40,50},
    cycle list name=black white,
    legend pos= outer north east
    ]
    \addplot coordinates {
(10,1.012519803)
(20, 1.00674601)
(30, 1.004441524)
(40, 1.002898082)
    };\pgfplotsset{cycle list shift=5}
    \addplot coordinates {
(10,	1.014017216)
(20, 1.007014511)
(30, 1.004791857)
(40, 1.002520881)
    };
   \legend{SD-AR,SI-AR}
  \end{axis}
\end{tikzpicture}%
}
\caption{different numbers of customer sites}
\end{subfigure}
    \caption{Statistics of the approximation ratios over 100 instances with different numbers of branches $C$, stages $T$, risk attitudes $\lambda$, demand standard deviations $\sigma$, facilities $M$ and customer sites $N$.}
    \label{fig:sensitivity_ratio}
\end{figure}

Next, we vary the number of branches $C$ from 2 to 5, number of stages $T$ from 3 to 6, number of facilities $M$ from 5 to 20 and number of customer sites $N$ from 10 to 40 to see how the empirical approximation ratio changes with respect to different parameter settings. The results are presented in Figure \ref{fig:sensitivity_ratio}, where we plot the mean of empirical approximation ratios (i.e., $\frac{z_R^{MS}(\boldsymbol{x}_n^{H}, \eta_n^{H}, \boldsymbol{y}_n^{H}, u_n^{H})}{z_R^{MS}(\boldsymbol{x}_n^{*}, \eta_n^{*}, \boldsymbol{y}_n^{*}, u_n^{*})}$) over 100 independently generated instances.
From the figure, the approximation ratios decrease gradually when we increase the number of branches $C$ and the number of stages $T$ (see Figures \ref{fig:sensitivity_ratio}(a) and (b)). 
Moreover, from Figures \ref{fig:sensitivity_ratio}(c) and (d), the approximation ratios are clearly positively related to the number of facilities $M$ and negatively impacted by the number of customer sites $N$, which is because increasing $N$ would increase the {total demand and} the number of units of resources needed in the first stage ($M_{\rm min}$) and our approximation ratio has an upper bound $1+ {\frac{M\sum_{t=1}^Tf_{t,\rm max}}{M_{\rm min}\sum_{t=1}^Tf_{t,\rm min}+\sum_{t=1}^Tc_{t,\rm min}\min_{n\in\mathcal{T}_t}\{\sum_{j=1}^Nd_{n,j}\}}}$. In all of these instances, the approximation algorithm achieves an approximation ratio of at most 1.03.

\subsubsection{Computational Time Comparison}
\label{sec:time}
We end this section by comparing the computational time of two-stage and multistage models solved via Gurobi, and multistage models solved by our approximation algorithms (AA). We also use the Python package \cite{ding2019python} to solve the multistage model via SDDiP algorithm \cite{zou2019stochastic} and terminate the algorithm after 10 iterations when the objective values become stable. Since SDDiP algorithm assumes that the underlying stochastic process is stagewise independent, we only consider SI scenario trees when evaluating their performance in this section. We note that our approximation algorithms can be also applied to solve SD scenario trees and we present the computational time comparison between our approximation algorithm and directly solving the problem using Gurobi 10.0.0 under SD scenario trees in Appendix D in the Online Supplement. We fix the risk parameters $\lambda=0.5,\ \alpha=0.95$ and the standard deviation $\sigma=0.8$ while varying the number of branches $C$ from 2 to 5 in Table \ref{tab:time-C}, number of stages $T$ from 3 to 6 in Table \ref{tab:time-T}, 
respectively. The percentages shown in parenthesis are the optimality gaps produced by Gurobi within the 1-hour time limit. We also display the gaps of the optimal objective value between our approximation algorithm/SDDiP and Gurobi in Column ``Gap''.
Note that we solve Step \ref{alg:step5} for all node $n\in\mathcal{T}$ as a whole in Algorithm \ref{alg:approx-multistage}, and one may further speed up the approximation algorithms by utilizing parallel computing techniques. From the tables, we observe that the approximation algorithms scale very well with respect to the problem size, while Gurobi cannot solve larger-scale risk-averse multistage models to optimality within the 1-hour time limit.
Compared to multistage models, two-stage models are less computationally expensive and therefore can be preferred over multistage models when their objective gaps are relatively small. While our approximation algorithms provide an upper bound on the multistage model by producing a feasible solution, SDDiP algorithm returns a lower bound on the multistage model by constructing under-approximations of the value functions. As the problem size increases, SDDiP algorithm fails to produce a high-quality lower bound within several hours (with the highest gap being $17.02\%$), while our approximation algorithms can always find a high-quality upper bound within several seconds (with all the gaps below 2\%).

\begin{table}[ht!]
% \OneAndAHalfSpacedXII
  \centering
  \caption{Computational time comparison with SI scenario trees and different number of branches $C$}
  \resizebox{\textwidth}{!}{
    \begin{tabular}{r|rr|rr|rrr|rrr}
         \hline
          & \multicolumn{2}{c|}{$z_R^{TS}$ via Gurobi} & \multicolumn{2}{c|}{$z_R^{MS}$ via Gurobi} & \multicolumn{3}{c|}{$z_R^{MS}$ via AA} & \multicolumn{3}{c}{$z_R^{MS}$ via SDDiP}\\
   $C$     & Time (sec.)  & Obj. (\$) & Time (sec.)  & Obj. (\$) & Time (sec.)  & Obj. (\$) & Gap & Time (sec.)  & Obj. (\$) & Gap\\
    \hline
    2     & 1.46  & 13,482K & 3.97  & 13,330K & \bf{0.13}  & 13,592K & 1.97\% & 90.06 & 13,270K & -0.45\% \\
    3     & 6.40  & 21,714K & 3600 (0.02\%) & 20,137K & \bf{0.31}  & 20,474K & 1.67\% & 679.13 & 20,098K & -0.19\% \\
    4     & 3.55  & 24,363K & 3600 (0.11\%) & 22,322K & \bf{0.43}  & 22,712K & 1.75\% & 2102.27 & 22,272K & -0.22\%\\
    5     & 10.01 & 24,142K & 3600 (0.19\%) & 22,520K & \bf{0.66}  & 22,907K & 1.72\% & 3430.90 & 22,467K & -0.24\% \\
    \hline
    \end{tabular}%
    }
  \label{tab:time-C}%
\end{table}%

\begin{table}[ht!]
% \OneAndAHalfSpacedXII
  \centering
  \caption{Computational time comparison with SI scenario trees and different number of stages $T$}
  \resizebox{\textwidth}{!}{
    \begin{tabular}{r|rr|rr|rrr|rrr}
          \hline
          & \multicolumn{2}{c|}{$z_R^{TS}$ via Gurobi} & \multicolumn{2}{c|}{$z_R^{MS}$ via Gurobi} & \multicolumn{3}{c|}{$z_R^{MS}$ via AA} & \multicolumn{3}{c}{$z_R^{MS}$ via SDDiP}\\
    $T$     & Time (sec.)  & Obj. (\$) & Time (sec.)  & Obj. (\$) & Time (sec.)  & Obj. (\$) & Gap & Time (sec.)  & Obj. (\$) & Gap \\
    \hline
    3     & 1.46  & 13,482K & 3.97 & 13,330K & \bf{0.13}  & 13,592K & 1.97\% & 90.06 & 13,270K & -0.45\% \\
    4     & 8.49  & 26,484K & 3600 (0.18\%) & 26,212K & \bf{0.67}  & 26,568K &1.36\% & 847.48 & 26,101K & -0.42\% \\
    5     & 3600 (0.05\%) & 44,003K & 3600 (0.14\%) & 43,205K & \bf{1.07}  & 43,623K & 0.97\% & 1028.06 & 43,063K & -0.33\% \\
    6     & 83.80 & 63,440K & 3600 (0.17\%) & 62,624K & \bf{6.97}  & 63,177K & 0.88\% & 1108.69 & 62,467K & -0.25\% \\
    7     & 640.56 & 138,372K & 3600 (0.13\%) & 112,508K & \bf{5.80} & 113,207K &0.62\% & 5556.12 & 93,360K & -17.02\% \\
    8     & 3600 (0.01\%) & 201,522K & 3600 (0.14\%) & 150,002K & \bf{14.64} & 150,718K &0.48\% & 13859.29 & 125,095K &-16.60\%\\
    \hline
    \end{tabular}%
    }
  \label{tab:time-T}%
\end{table}%

\subsection{Case Study on Real-World EV Charging Station Capacity Planning}
\label{sec:real}
\subsubsection{Experimental Design and Setup}
We consider the 49-node and 88-node data sets described in \citet{daskin2011network} with $M=49,\ N=88,\ T=5,\ C=2$ where each stage is one year. The maintenance cost for each charger is set to \$100 per year \cite{EVmainten}.
The capacity of each charger $h_{ti}$ is set to $6\times 360$ assuming that each charger can charge 6 EVs from empty per day \cite{EVchargingspeed}.
The resource-allocation costs are set equal to the great-circle distance times the travel cost per mile per unit of demand, i.e., $c_{ij}=\textrm{dist}(i,j) * 0.00001$.
We then collect the population data in each city and multiply them by 6\% times 120 days, assuming that 6\% of the population currently own EVs \cite{EVpercentage} and will charge their EVs every 3 days \cite{EVchargingfrequency}, which gives nominal demand $\tilde{d}$ in each customer site in the beginning year of the planning horizon. 
We consider four demand patterns (described in Column ``Pattern'' in Table \ref{table:pattern}), all of which follow truncated Normal distributions. In Patterns III and IV, the nominal demand $\tilde{d}$ is increased with a rate of $200\%$ for each subsequent year, which will reach $6\%\times (1+2\times4)=54\%$ of EVs market share by the end of the fifth year, matching President Biden’s goal of having 50 percent of all new vehicle sales be electric \cite{BidenEVGoal}; in Patterns II and IV, the standard deviation is increased with the same rate. 
\begin{table}[ht!]
% \OneAndAHalfSpacedXI
	\centering
	\caption{Demand patterns}
	\begin{tabular}{lr}
		\hline
Pattern & Distribution at stage $t$ \\
\hline
I.\ Constant mean, constant standard deviation & $\mathcal N(\tilde{d},\tilde{d}\cdot \sigma)$\\
 II.\ Constant mean, increasing standard deviation & $\mathcal N(\tilde{d},\tilde{d}\cdot (\sigma+2(t-1)))$\\
 III.\ Increasing mean, constant standard deviation & $\mathcal N(\tilde{d}\cdot (1+2(t-1)),\tilde{d}\cdot \sigma)$\\
 IV.\ Increasing mean, increasing standard deviation & $\mathcal N(\tilde{d}\cdot (1+2(t-1)),\tilde{d}\cdot (\sigma+2(t-1)))$\\
           \hline
\end{tabular}%
\label{table:pattern}%
\end{table}%

\subsubsection{Result Comparison with Different Demand Patterns}
With the baseline setting $C=2,\ \sigma=0.8,\ \lambda_t=0.5,\ \alpha_t=0.95$ and SD scenario trees, we present the optimal solutions and cost breakdown of two-stage and multistage models under different demand patterns in Table \ref{tab:real-sol}.  We use Gurobi to solve the risk-averse two-stage (TS) and multistage (MS) models directly, and we compare these two models with the multistage models solved via approximation algorithms (AA). Columns ``$|\boldsymbol{x}_1|$'', ``$|\boldsymbol{x}_2|$'', ``$|\boldsymbol{x}_3|$'', ``$|\boldsymbol{x}_4|$'' and ``$|\boldsymbol{x}_5|$'' display the average number of chargers installed in each stage across all scenarios. Columns ``Maintenance (\$)'' and  ``Operational (\$)'' show the maintenance and operational cost without considering the risk parameters, and Column ``Obj. (\$)'' presents the overall risk-averse optimal objective values, where we mark the lowest ones among the three models in bold. The last column ``RVMS/AR'' shows the relative VMS or the approximation ratios.  Comparing the three models, multistage models solved by Gurobi always achieve the least maintenance cost, because they have more flexibility in deciding the capacity-expansion plans. In terms of operational cost, the multistage models solved by approximation algorithms always achieve the minimum among the three. 
In all of the cases, our approximation algorithm aligns with the optimal solutions and objective values of the multistage models very well and achieves an approximation ratio of at most 1.00004.

We also find that the demand patterns with increasing standard deviation always obtain higher ${\rm RVMS_R}$ compared to the ones with constant standard deviation, which agrees with our findings in Section \ref{sec:sensi}. In Pattern I, most chargers are invested in the first stage because the demand mean is constant and we do not need to build much more chargers in later stages. Under this setting, we achieve an RVMS of 5.19\%. In Patterns III and IV where the demand mean is increased at a rate of 200\%, more chargers are installed in later stages, and the ${\rm RVMS}$ increases to $12.02\%$ and 18.61\%, respectively. Note that this gap can be further increased with more branches $C$ and higher standard deviations $\sigma$.  

\begin{table}[ht!]
 % \OneAndAHalfSpacedXII
  \centering
  \caption{Optimal solutions and costs of two-stage and multistage models under different demand patterns}
  \resizebox{\textwidth}{!}{
    \begin{tabular}{clrrrrrrrrrr}
     \hline
    \multicolumn{1}{l}{Pattern} & Model & $|\boldsymbol{x}_1|$  & $|\boldsymbol{x}_2|$  & $|\boldsymbol{x}_3|$  & $|\boldsymbol{x}_4|$  & $|\boldsymbol{x}_5|$  & Maintenance (\$) & Operational (\$) & Obj. (\$)  & Time & RVMS/AR \\
    \hline
    \multirow{3}[0]{*}{I} & TS    & 149,502 & 33,681 & 3,845 & 16,686 & 20,869 & \$94,801K & \$2,368K & \$97,336K & 8.90  & 5.19\% \\
          & MS    & 149,489 & 29,303 & 961   & 7,018 & 6,509 & \textbf{\$88,808K} & \$2,492K & \textbf{\$92,284K} & 8.94  &  \\
          & AA    & 149,489 & 29,315 & 962   & 7,018 & 6,509 & \$88,813K & \textbf{\$2,491K} & \$92,288K & 58.72 & 1.00004 \\
          \hline
    \multirow{3}[0]{*}{II} & TS    & 149,523 & 261,128 & 265,695 & 389,223 & 488,955 & \$385,661K & \$9,891K & \$396,742K & 8.68  & 14.82\% \\
          & MS    & 149,469 & 253,004 & 199,195 & 279,369 & 266,388 & \textbf{\$318,207K} & \$9,062K & \textbf{\$337,928K} & 57.36 &  \\
          & AA    & 149,494 & 252,997 & 199,198 & 279,371 & 266,388 & \$318,218K & \textbf{\$9,062K} & \$337,939K & 56.65 & 1.00003 \\
          \hline
    \multirow{3}[0]{*}{III} & TS    & 149,531 & 399,975 & 385,524 & 490,802 & 595,216 & \$508,095K & \$12,256K & \$521,440K & 8.81  & 12.02\% \\
          & MS    & 149,469 & 386,866 & 305,341 & 384,100 & 374,120 & \textbf{\$435,315K} & \$11,303K & \textbf{\$458,753K} & 750.35 &  \\
          & AA    & 149,495 & 386,858 & 305,345 & 384,101 & 374,120 & \$435,326K & \textbf{\$11,303K} & \$458,763K & 59.20 & 1.00002 \\
          \hline
    \multirow{3}[0]{*}{IV} & TS    & 149,528 & 1,082,360 & 2,149,745 & 4,075,429 & 6,453,154 & \$2,613,033K & \$62,620K & \$2,683,077K & 9.00  & 18.61\% \\
          & MS    & 149,470 & 1,057,952 & 1,800,916 & 3,158,538 & 3,797,729 & \textbf{\$2,049,671K} & \$55,468K & \textbf{\$2,183,774K} & 67.52 &  \\
          & AA    & 149,495 & 1,057,946 & 1,800,916 & 3,158,537 & 3,797,730 & \$2,049,681K & \textbf{\$55,467K} & \$2,183,784K & 56.94 & 1.000005 \\
          \hline
    \end{tabular}%
    }
  \label{tab:real-sol}%
\end{table}%

\subsubsection{Effect of the Risk Parameter}
\begin{figure}[ht!]
\centering
\begin{subfigure}{0.42\textwidth}
 \resizebox{0.95\textwidth}{!}{%
\begin{tikzpicture}
\centering
  \begin{axis}
  [
    xlabel={Risk parameter $\lambda$},
    ylabel={Optimal Objective Values},
    xtick={0, 0.2, 0.4, 0.6, 0.8, 1},
cycle list name=black white,
legend pos= outer north east
]
    \addplot coordinates {
(0, 90520249.016607)
(0.2, 90549055.520645)
(0.4, 90574818.785963)
(0.6, 90600477.883901)
(0.8, 90623275.1128686)
(1, 90649977.3990519)
    };\pgfplotsset{cycle list shift=5}
    \addplot coordinates {
(0, 85459770.0825577)
(0.2, 85811286.2541412)
(0.4, 86163337.676473)
(0.6, 86508817.4987871)
(0.8, 86846701.4455332)
(1, 87178874.4049586)
    };
    \legend{Two-stage,Multistage}
  \end{axis}
\end{tikzpicture}%
}
\caption{Optimal objective values in two models.}
\end{subfigure}
\begin{subfigure}{0.42\textwidth}
 \resizebox{\textwidth}{!}{%
\begin{tikzpicture}
  \begin{axis}
  [yticklabel style={/pgf/number format/fixed,
                  /pgf/number format/precision=3},
    xlabel={Risk parameter $\lambda$},
    ylabel={95\% Percentile},
    axis y line* = left,
    xtick={0, 0.2, 0.4, 0.6, 0.8, 1},
cycle list name=black white
]
    \addplot coordinates {
(0, 18082978)
(0.2, 18075310)
(0.4, 18068612)
(0.6, 18064585)
(0.8, 18061029)
(1, 18058181)
    };
     \label{percentile}
    \end{axis}
    \pgfplotsset{cycle list shift=5}
    \begin{axis}
      [yticklabel style={/pgf/number format/fixed,
                  /pgf/number format/precision=4},
    xlabel={Risk parameter $\lambda$},
    ylabel={Mean},
    axis y line* = right,
    xtick={0, 0.2, 0.4, 0.6, 0.8, 1},
cycle list name=black white,
legend pos= outer north east
]
\addlegendimage{/pgfplots/refstyle=percentile}\addlegendentry{95\% Percentile}
    \addplot coordinates {
(0, 17091954)
(0.2, 17091261)
(0.4, 17093566)
(0.6, 17096689)
(0.8, 17100078)
(1, 17109673)
    };\addlegendentry{Mean}
  \end{axis}
\end{tikzpicture}%
}
\caption{Stagewise cost in multistage models.}
\end{subfigure}
\caption{Optimal objective value and stagewise cost comparison between two-stage and multistage models with varying risk attitude $\lambda$.}
\label{fig:lambda-obj}
\end{figure}
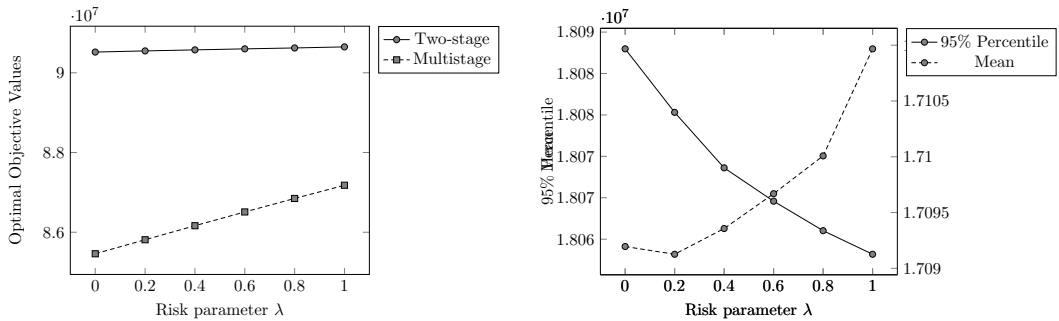

Finally, we discuss the effect of the risk attitude on the optimal solutions. We display the optimal objective values of two-stage and multistage models with varying risk attitude $\lambda$ in Figure \ref{fig:lambda-obj}(a) and the 95\% percentile and mean of the stagewise cost in multistage models with varying risk attitude $\lambda$ in Figure \ref{fig:lambda-obj}(b). From Figure \ref{fig:lambda-obj}(a), the optimal objective values of both two-stage and multistage models increase as we become more risk-averse, agreeing with our results in Proposition \ref{prop:risk}. Moreover, the gap between the two models reduces as we increase $\lambda$, as shown in Section \ref{sec:sensi}, Figure \ref{fig:sensitivity_RVMS}(c). From Figure \ref{fig:lambda-obj}(b), when we increase $\lambda$, the mean of the stagewise cost increases, while the 95\% percentile decreases in most cases, representing a switch from risk-neutral to risk-averse attitude.

\section{Conclusion}
\label{sec:conclu}
We considered a general class of multiperiod capacity planning problems under uncertain demand in each period. We compared a multistage stochastic programming model where the capacities of facilities can be determined dynamically throughout the uncertainty realization process, with a two-stage model where decision makers have to fix capacity acquisition at the beginning of the planning horizon. Using expected conditional risk measures (ECRMs), we bounded the gaps between the optimal objective values of risk-averse two-stage and multistage models from below and above and provided an example to show that one of the lower bounds is tight. We also proposed an approximation algorithm to solve the risk-averse multistage models more efficiently, which are asymptotically optimal in an expanding market. Our numerical tests indicate that the $\rm{RVMS}$ increases as the uncertainty variability increases and decreases as the decision maker becomes more risk-averse. Moreover, stagewise-dependent scenario trees attain higher $\rm{RVMS}$ than the stagewise-independent counterparts. On the other hand,
the analytical lower bounds can recover the true gaps very well. Moreover, the approximation ratios are pretty close to 1 in the case study based on a real-world network.

There are several interesting directions to investigate for future research. The risk measure we used in this paper is the ECRM, of which the risk is measured separately for each stage. There are other risk-measure choices that can be applied here, such as the nested risk measures. More research can be done to explore the relationship between these two models under different choices of risk measures. Moreover, this paper assumes that the uncertainty has a known distribution, while it is more realistic to assume that the probability distribution of the uncertainty belongs to an ambiguity set. Therefore, a distributionally robust optimization framework can be considered for two-stage and multistage capacity planning problems. As a middle ground between two-stage and multistage decision frameworks, \cite{basciftci2019adaptive} studied an adaptive two-stage stochastic programming framework, where each component of the decision policy has its own revision point, before which the decisions are determined until this revision point, and after which they are revised for adjusting to the uncertainty. They analyzed this approach on a capacity expansion planning problem and derived bounds on the value of the proposed adaptive two-stage stochastic programming. Similarly, we can obtain a relaxation of the multi-stage stochastic program by ignoring the non-anticipativity constraints at some non-root stages, and such a relaxation provides a lower bound to the multistage model. In this paper, we focus on the comparison between risk-averse two-stage and multistage stochastic programs and leave the exploration of general adaptive stochastic programs as future research.

\paragraph{Acknowledgment}
The authors are grateful for the constructive feedback from the editorial team and two anonymous reviewers. The work of Dr. Xian Yu was partially supported by the United States National Science Foundation Grant \#2331782. 
% Leave this (end of acknowledgment)

% \bibliographystyle{informs2014}

%\THEEndNotes
% \begingroup \parindent 0pt \parskip 0.0ex \def\enotesize{\normalsize} \theendnotes \endgroup

% Appendix here
% Options are (1) APPENDIX (with or without general title) or
%             (2) APPENDICES (if it has more than one unrelated sections)
% Outcomment the appropriate case if necessary
%
% \begin{APPENDIX}{<Title of the Appendix>}
% \end{APPENDIX}
%
%   or
%
% \begin{APPENDICES}
% \section{<Title of Section A>}
% \section{<Title of Section B>}
% etc
% \end{APPENDICES}

% Acknowledgments here
% \ACKNOWLEDGMENT{We would like to express our sincere gratitude to [acknowledge individuals, organizations, or institutions] for their invaluable contributions to this research. We are also grateful to [mention any additional acknowledgements, such as technical assistance, data providers, or colleagues] for their support and assistance throughout the course of this work.}

% References here (outcomment the appropriate case)

% CASE 1: BiBTeX used to constantly update the references
%   (while the paper is being written).
%\bibliographystyle{informs2014} % outcomment this and next line in Case 1
%\bibliography{<your bib file(s)>} % if more than one, comma separated

%\bibliographystyle{informs2014} % outcomment this and next line in Case 1
%\bibliography{sample} % if more than one, comma separated

% CASE 2: BiBTeX used to generate mypaper.bbl (to be further fine tuned)
%\input{mypaper.bbl} % outcomment this line in Case 2

%If you don't use BiBTex, you can manually itemize references as shown below.

%\bibliographystyle{nonumber}

\bibliographystyle{apalike}
\bibliography{Xian_bib}

\newpage
\appendix
\begin{center}
    \Large \textbf{Online Supplement to ``On the Value of Risk-Averse Multistage Stochastic Programming in Capacity Planning''}\\
    \normalsize \textbf{Xian Yu and Siqian Shen}
\end{center}
This online supplement is organized as follows.
\begin{itemize}
    \item Appendix \ref{e-companion:ECRMs}: Time Consistency of ECRMs.
    \item Appendix \ref{e-companion:proofs}: Details of All Needed Proofs.
    \item Appendix \ref{e-companion:example}: An Example to Illustrate the Usefulness of the Derived Bounds.
        \item Appendix \ref{e-companion:results}: Additional Computational Results.
\end{itemize}

\section{Time Consistency of ECRMs}\label{e-companion:ECRMs}
\normalsize
We follow the definition in \cite{ruszczynski2010risk} and prove the time consistency of ECRMs \eqref{ECRMs}.
Consider the probability space $(\Xi,\mathcal{F},P)$, and let $\mathcal{F}_1\subset\mathcal{F}_2\subset\ldots\subset\mathcal{F}_T$ be sub-sigma-algebras of $\mathcal{F}$ such that each $\mathcal{F}_t$ corresponds to the information available up to (and including) stage $t$, with $\mathcal{F}_1=\{\emptyset,\Xi\}, \ \mathcal{F}_T=\mathcal{F}$. Let $\mathcal{Z}_t$ denote a space of $\mathcal{F}_t$-measurable functions from $\Xi$ to $\mathbb{R}$, and let $\mathcal{Z}_{1,T}:=\mathcal{Z}_1\times\cdots\times\mathcal{Z}_T$.
\begin{definition}
A mapping $\rho_{t,T}:\mathcal{Z}_{t,T}\to\mathcal{Z}_t$ where $1\le t\le T$ is called a \textit{conditional risk measure}, if it has the following monotonicity property:
$\rho_{t,T}(Z)\le\rho_{t,T}(W)$ for all $Z,W\in\mathcal{Z}_{t,T}$ such that $Z\le W$.
\end{definition}
\begin{definition}
A dynamic risk measure is a sequence of monotone one-step conditional risk measures $\rho_{t,T}: \mathcal{Z}_{t,T}\to\mathcal{Z}_t,\ 1\le t\le T$.
\end{definition}
\begin{definition}\label{def:time-consistent}
A dynamic risk measure $\{\rho_{t,T}\}_{t=1}^{T}$ is called time consistent if, for all $1\le l<k\le T$ and all sequences $Z, \ W\in \mathcal{Z}_{l, T}$, the conditions $Z_i = W_i,\ \forall i=l,\ldots, k-1, \text{and}\ 
    \rho_{k,T}(Z_k,\ldots,Z_T)\le \rho_{k,T}(W_k,\ldots,W_T)$
imply that $\rho_{l,T}(Z_l,\ldots,Z_T)\le 
     \rho_{l,T}(W_l,\ldots,W_T).$
\end{definition}
\begin{theorem}\label{thm:ECRMs}
The ECRMs defined in \eqref{ECRMs} are time consistent, if each $\rho_t^{d_{[1,t-1]}}$ is a coherent one-step conditional risk measure.
\end{theorem}
\proof{Proof of Theorem \ref{thm:ECRMs}}
According to Eq. \eqref{eq:obj} and the translation-invariant property of $\rho_t^{d_{[t-1]}}$, the risk function \eqref{ECRMs} can be recast as
\begin{align*}
\mathbb{F}(g_1,\ldots,g_{T})=g_1+\rho_2\Big(g_2+\mathbb{E}_{\boldsymbol d_{2}}\circ{\rho_3^{d_{[2]}}}\Big(g_3+\mathbb{E}_{\boldsymbol d_{3}|d_{[2]}}\circ{\rho_4^{d_{[3]}}}\Big(g_4+\cdots\nonumber+\mathbb{E}_{\boldsymbol d_{T-1}|d_{[T-2]}}\circ{\rho_{T}^{d_{[T-1]}}}\Big(g _{T}\Big)\Big)\cdots\Big)\Big),
\end{align*}
To simplify the notation, we define $\tilde{\rho}_t^{d_{[t-2]}}: = \mathbb{E}_{\boldsymbol{d}_{t-1}|d_{[t-2]}}\circ \rho_t^{d_{[t-1]}}$, which maps from $\mathcal{Z}_t$ to $\mathcal{Z}_{t-2}$. Then, the multiperiod risk function $\mathbb{F}$ can be recast as
\begin{align*}
\mathbb{F}(g_1,\ldots,g_{T})=g_1+\rho_2\Big(g_2+\tilde{\rho}_3^{d_{[1]}}\Big(g_3+\tilde{\rho}_4^{d_{[2]}}\Big(g_4+\cdots
+\tilde{\rho}_T^{d_{[T-2]}}\Big(g _{T}\Big)\Big)\cdots\Big)\Big),
\end{align*}
Define a dynamic risk measure $\{\mathbb{F}_{t,T}\}_{t=1}^T$ as
$\mathbb{F}_{1,T}=\mathbb{F}$, and for $2\le t\le T$, we have
\begin{align*}
    \mathbb{F}_{t,T}(g_t,\ldots,g_T) = g_t+\tilde{\rho}_{t+1}^{d_{[t-1]}}\Big(g_{t+1}+\tilde{\rho}_{t+2}^{d_{[t]}}\Big(g_{t+2}+\cdots
+\tilde{\rho}_T^{d_{[T-2]}}\Big(g _{T}\Big)\Big)\cdots\Big)
\end{align*}
Then for $1\le l < k\le T$,
\begin{align*}
    \mathbb{F}_{l,T}(g_l,\ldots,g_T) =g_l + \tilde{\rho}_{l+1}^{d_{[l-1]}}\Big(g_{l+1} +\ldots+\tilde{\rho}_k^{d_{[k-2]}}\Big(\mathbb{F}_{k,T}(g_k,\ldots,g_T)\Big)\Big).
\end{align*}
If $\mathbb{F}_{k,T}(Z_k,\ldots,Z_T)\le \mathbb{F}_{k,T}(W_k,\ldots,W_T)$ and $Z_i=W_i,\ \forall i=l,\ldots, k-1$, we have
$\mathbb{F}_{l,T}(Z_l,\ldots,Z_T)\le \mathbb{F}_{l,T}(W_l,\ldots,W_T)$ because of the monotonicity of $\tilde{\rho}_t^{d_{[t-2]}},\ \forall t\ge {l+1}$.
This completes the proof.
 
\endproof

\section{Details of All Needed Proofs}\label{e-companion:proofs}
\normalsize

{\color{black}
\proof{Proof of Proposition \ref{prop:risk}}
According to Eq. \eqref{eq:multi-constraint_y}--\eqref{eq:constraint_eta}, the feasible region of the multistage model \eqref{model:multiriskaverse} is independent of parameter $\lambda$. Given two sets of parameters $0\le (\lambda_2,\lambda_3,\ldots,\lambda_T)<(\hat{\lambda}_2,\hat{\lambda}_3,\ldots,\hat{\lambda}_T)\le 1$ (component-wise inequality), define $\mathbb{F}$ and $\hat{\mathbb{F}}$ as the ECRM value related to $\boldsymbol{\lambda}$ and $\hat{\boldsymbol{\lambda}}$, and $z_{R}^{MS}$ and $\hat{z}_{R}^{MS}$ as the corresponding optimal objective values of model \eqref{model:multiriskaverse}, respectively. For any feasible solution $(\boldsymbol{x},\boldsymbol{y})$, according to definitions \eqref{ECRMs} and \eqref{eq:rho}, we have for all $t=2,\ldots,T$, $(1-\lambda_t)\mathbb{E}[g_t(\boldsymbol{x}_t,\boldsymbol{y}_t)|\boldsymbol{d}_{t-1}]+\lambda_t{\rm CVaR}_{\alpha_t}^{\boldsymbol{d}_{[t-1]}}[g_t(\boldsymbol{x}_t,\boldsymbol{y}_t)] < (1-\hat{\lambda}_t)\mathbb{E}[g_t(\boldsymbol{x}_t,\boldsymbol{y}_t)|\boldsymbol{d}_{t-1}]+\hat{\lambda}_t{\rm CVaR}_{\alpha_t}^{\boldsymbol{d}_{[t-1]}}[g_t(\boldsymbol{x}_t,\boldsymbol{y}_t)]$ because $\mathbb{E}[g_t(\boldsymbol{x}_t,\boldsymbol{y}_t)|\boldsymbol{d}_{t-1}]< {\rm CVaR}_{\alpha_t}^{\boldsymbol{d}_{[t-1]}}[g_t(\boldsymbol{x}_t,\boldsymbol{y}_t)]$ when $0<\alpha_t<1$. Due to the monotonicity of expectation, we conclude that $\mathbb{F}(g_1(\boldsymbol{x}_1,\boldsymbol{y}_1),\ldots,g_T(\boldsymbol{x}_T,\boldsymbol{y}_T)) < \hat{\mathbb{F}}(g_1(\boldsymbol{x}_1,\boldsymbol{y}_1),\ldots,g_T(\boldsymbol{x}_T,\boldsymbol{y}_T))$. Thus, $z_{R}^{MS}<\hat{z}_{R}^{MS}$. The same reasoning also applies to the optimal objective value of the two-stage model $z_{R}^{TS}$.
  
}
\proof{Proof of Proposition \ref{prop:substructure}}
We first check the feasibility of solutions $\{\boldsymbol{x}^{MS}_{n}\}_{n\in \mathcal{T}},\ \{\eta_n^{MS}\}_{n\in \mathcal{T}\setminus\mathcal{L}}$ to $\mbox{{\bf SP-RMS}}(\boldsymbol{y}_n^{*}, u_n^{*})$. Note that for all $n\in\mathcal{T}$, we have $\boldsymbol{x}^{MS}_{n} \in \mathbb{Z}_{+}^{M}$ and $\sum_{m\in \mathcal{P}(n)}\boldsymbol{x}^{MS}_m = \max_{m\in\mathcal{P}(n)}\lceil \boldsymbol{B}_{t_m} \boldsymbol{y}_m^{*}\rceil \ge \boldsymbol{B}_{t_n} \boldsymbol{y}_n^{*}$. By the definition of $\eta_n^{MS}$, we have ${\eta}^{MS}_n \ge  \boldsymbol{f}^{\mathsf T}_{t_m}\sum_{l\in\mathcal{P}(m)}\boldsymbol{x}^{MS}_l+ \boldsymbol{c}_{t_m}^{\mathsf T} \boldsymbol{y}_m^{*}-u_m^{*},\ \forall m\in \mathcal{C}(n),\ n\not\in\mathcal{L}$, which satisfy constraints \eqref{eq:constraint_eta_knowny} automatically. As a result, we obtain 
\begin{align*}
Q^M(\boldsymbol{y}_n^{*}, u_n^{*}) &\le \sum_{n\in \mathcal{T}} p_n\left(\boldsymbol{\tilde{f}}_n^{\mathsf T}\sum_{m\in \mathcal{P}(n)}\boldsymbol{x}^{MS}_m+ \tilde{\lambda}_n \eta^{MS}_n\right)\\
&= \sum_{n\in \mathcal{T}} p_n\left(\boldsymbol{\tilde{f}}_n^{\mathsf T}\max_{m\in\mathcal{P}(n)}\lceil \boldsymbol{B}_{t_m} \boldsymbol{y}_m^{*}\rceil+ \tilde{\lambda}_n \max_{m\in \mathcal{C}(n)}\left\lbrace \boldsymbol{f}^{\mathsf T}_{t_m}\max_{l\in\mathcal{P}(m)}\lceil \boldsymbol{B}_{t_l} \boldsymbol{y}_l^{*}\rceil+ \boldsymbol{c}_{t_m}^{\mathsf T} \boldsymbol{y}_m^{*}-u_m^{*}\right\rbrace\right).
\end{align*}

Next, we show the optimality of solutions $\{\boldsymbol{x}^{MS}_{n}\}_{n\in \mathcal{T}},\ \{\eta_n^{MS}\}_{n\in \mathcal{T}\setminus\mathcal{L}}$ to $\mbox{{\bf SP-RMS}}(\boldsymbol{y}_n^{*}, u_n^{*})$. Note that for any feasible solution $\boldsymbol{x}_{n}\in \mathbb{Z}^M_+, \ \forall n\in\mathcal{T}$, from constraints \eqref{eq:constraint_xy_knowny}, we get $\sum_{l\in \mathcal{P}(n)}\boldsymbol{x}_l\ge \sum_{l\in \mathcal{P}(m)}\boldsymbol{x}_l\ge \lceil\boldsymbol{B}_{t_m}\boldsymbol{y}_m^*\rceil,\ \forall m\in \mathcal{P}(n)$ (we can raise $\boldsymbol{B}_{t_m}\boldsymbol{y}_m^*$ to $\lceil\boldsymbol{B}_{t_m}\boldsymbol{y}_m^*\rceil$ because of the integrality of $\boldsymbol{x}_n$-variables) and thus $\sum_{l\in \mathcal{P}(n)}\boldsymbol{x}_l\ge \max_{m\in \mathcal{P}(n)}\lceil\boldsymbol{B}_{t_m}\boldsymbol{y}_m^*\rceil$. From constraints \eqref{eq:constraint_eta_knowny}, we obtain ${\eta}_n \ge  \max_{m\in \mathcal{C}(n)}\{\boldsymbol{f}^{\mathsf T}_{t_m}\sum_{l\in\mathcal{P}(m)}\boldsymbol{x}_l+ \boldsymbol{c}_{t_m}^{\mathsf T} \boldsymbol{y}_m^{*}-u_m^{*}\}\ge \max_{m\in \mathcal{C}(n)}\{\boldsymbol{f}^{\mathsf T}_{t_m}\max_{l\in \mathcal{P}(m)}\lceil\boldsymbol{B}_{t_l}\boldsymbol{y}_l^*\rceil+ \boldsymbol{c}_{t_m}^{\mathsf T} \boldsymbol{y}_m^{*}-u_m^{*}\},\ \forall n\not\in\mathcal{L}$. With that, we conclude 
\begin{align*}
Q^M(\boldsymbol{y}_n^{*}, u_n^{*}) \ge \sum_{n\in \mathcal{T}} p_n\left(\boldsymbol{\tilde{f}}_n^{\mathsf T}\max_{m\in\mathcal{P}(n)}\lceil \boldsymbol{B}_{t_m} \boldsymbol{y}_m^{*}\rceil+ \tilde{\lambda}_n \max_{m\in \mathcal{C}(n)}\left\lbrace \boldsymbol{f}^{\mathsf T}_{t_m}\max_{l\in\mathcal{P}(m)}\lceil \boldsymbol{B}_{t_l} \boldsymbol{y}_l^{*}\rceil+ \boldsymbol{c}_{t_m}^{\mathsf T} \boldsymbol{y}_m^{*}-u_m^{*}\right\rbrace\right).
\end{align*}

Therefore, $\{\boldsymbol{x}^{MS}_{n}\}_{n\in \mathcal{T}},\ \{\eta_n^{MS}\}_{n\in \mathcal{T}\setminus\mathcal{L}}$ are optimal solutions to the risk-averse multistage problem \eqref{eq:ms-substructure}. In the risk-averse two-stage problem \eqref{eq:ts-substructure}, because the investment decisions $\boldsymbol{x}$ are identical for all nodes in the same stage, one only needs to replace the right-hand side of constraints \eqref{eq:constraint_xy_knowny} with $\max_{l\in \mathcal{T}_{t_n}}\boldsymbol{B}_{t_l} \boldsymbol{y}_l^{*}$. 
% Moreover, to satisfy the two-stage constraints for $\boldsymbol{\eta}$, we need to ensure that ${\eta}_n \ge  \max_{m\in \mathcal{T}_{t_n+1}}\{\boldsymbol{f}^{\mathsf T}_{t_m}\sum_{l\in\mathcal{P}(m)}\boldsymbol{x}_l+ \boldsymbol{c}_{t_m}^{\mathsf T} \boldsymbol{y}_m^{*}-u_m^{*}\},\ \forall n\not\in\mathcal{L}$. 
Following the same analysis used for risk-averse multistage problem \eqref{eq:ms-substructure} and replacing $\boldsymbol{B}_{t_n} \boldsymbol{y}_n^{*}$ with $\max_{l\in \mathcal{T}_{t_n}}\boldsymbol{B}_{t_l} \boldsymbol{y}_l^{*}$, it can be shown that $\{\boldsymbol{x}^{TS}_{n}\}_{n\in \mathcal{T}},\ \{\eta_n^{TS}\}_{n\in \mathcal{T}\setminus\mathcal{L}}$ are optimal solutions to the risk-averse two-stage problem \eqref{eq:ts-substructure}. This completes the proof.   
\endproof

\proof{Proof of Theorem \ref{thm:risk-averse-bound}}
First of all, $\rm VMS_{R}^{LB}\ge 0$ and thus provides a nontrivial lower bound. Let $\left\lbrace {\boldsymbol{y}_n^{*}}\right\rbrace_{{n\in\mathcal{T}}},\ \{u_n^{*}\}_{n\in \mathcal{T}\setminus\{1\}}$ be the values of decisions made in an optimal solution to the two-stage model \eqref{model:tworiskaverse}. Because of the optimality of decisions $\boldsymbol{y}_n^{*},\ u_n^{*}$, we have $
	z^{TS}_R
    =\sum_{n\in \mathcal{T}} p_n(\boldsymbol{\tilde{c}}_n^{\mathsf T}\boldsymbol{y^{*}_n} + \tilde{\alpha}_n u^{*}_n)+Q^{T}(\boldsymbol{y}_n^{*},u_n^{*})$,
	where $Q^{T}(\boldsymbol{y}_n^{*},u_n^{*})$ is the optimal objective value to the substructure problem defined in \eqref{eq:ts-substructure}.
Since $\left\lbrace {\boldsymbol{y}_n^{*}}\right\rbrace_{{n\in\mathcal{T}}},\ \{u_n^{*}\}_{n\in \mathcal{T}\setminus\{1\}}$ are feasible solutions for the multistage problem \eqref{model:multiriskaverse}, we have
\begin{align}
z^{MS}_R\le \sum_{n\in \mathcal{T}} p_n(\boldsymbol{\tilde{c}}_n^{\mathsf T}\boldsymbol{y^{*}_n} + \tilde{\alpha}_n u^{*}_n)+Q^{M}(\boldsymbol{y}_n^{*},u_n^{*}),\label{eq:ineq}
\end{align}
where $Q^{M}(\boldsymbol{y}_n^{*},u_n^{*})$ is the optimal objective value to the substructure problem defined in \eqref{eq:ms-substructure}.

Combining all steps above, we have
\begin{align*}
{\rm VMS_R}=&z^{TS}_R-z^{MS}_R\\
\ge& Q^{T}(\boldsymbol{y}_n^{*},u_n^{*})-Q^M(\boldsymbol{y}_n^{*},u_n^{*})\\
\overset{(a)}{=}& \sum_{n\in \mathcal{T}} p_n\left(\boldsymbol{\tilde{f}}_n^{\mathsf T}\left(\max_{m\in \mathcal{P}(n)} \lceil \max_{l\in \mathcal{T}_{t_m}}\boldsymbol{B}_{t_l} \boldsymbol{y}_l^{*}\rceil - \max_{m\in \mathcal{P}(n)}\lceil \boldsymbol{B}_{t_m} \boldsymbol{y}_m^{*}\rceil\right)+ \tilde{\lambda}_n \left(\eta^{TS}_n-\eta^{MS}_n\right)\right)\\
=& {\sum_{n\in \mathcal{T}\setminus\{1\}} p_n(1-\lambda_{t_n})\boldsymbol{f}_{t_n}^{\mathsf T}\left(\max_{m\in \mathcal{P}(n)} \lceil \max_{l\in \mathcal{T}_{t_m}}\boldsymbol{B}_{t_l} \boldsymbol{y}_l^{*}\rceil - \max_{m\in \mathcal{P}(n)}\lceil \boldsymbol{B}_{t_m} \boldsymbol{y}_m^{*}\rceil\right)+ {\sum_{n\in \mathcal{T}\setminus\mathcal{L}} p_n {\lambda}_{t_n+1}} \left(\eta^{TS}_n-\eta^{MS}_n\right)}
\end{align*}
where the equality $(a)$ follows from Proposition \ref{prop:substructure}.
This completes the proof.   
\endproof

{
\proof{Proof of Corollary \ref{cor:LP}}
Because $\{\boldsymbol{y}_n^{TSLP}\}_{n\in \mathcal{T}}, \ \{u_n^{TSLP}\}_{n\in \mathcal{T}\setminus\{1\}}$ is an optimal solution to the LP relaxation of the two-stage model \eqref{model:tworiskaverse}, we have $
	z^{TS}_R
	%&=&\sum_{n\in \mathcal{T}} p_n(\boldsymbol{\tilde{c}}_n^{\mathsf T}\boldsymbol{y^{TS}_n} + \tilde{\alpha}_n u^{TS}_n)+\sum_{n\in \mathcal{T}} p_n(\boldsymbol{\tilde{f}}_n^{\mathsf T}\sum_{m\in \mathcal{P}(n)}\boldsymbol{x}^{TS}_m + \tilde{\lambda} \eta^{TS}_n)\\
    \ge\sum_{n\in \mathcal{T}} p_n(\boldsymbol{\tilde{c}}_n^{\mathsf T}\boldsymbol{y}^{TSLP}_n + \tilde{\alpha}_n u^{TSLP}_n)+Q^{TLP}(\boldsymbol{y}_n^{TSLP},u_n^{TSLP})$,
where $Q^{TLP}(\boldsymbol{y}_n^{TSLP},u_n^{TSLP})$ stands for the optimal objective value of the LP relaxation of the subproblem $\mbox{{\bf SP-RTS}}(\boldsymbol{y}_n^{TSLP}, u_n^{TSLP})$. It can be easily verified that the constructed $\boldsymbol{x}^{TS},\boldsymbol{\eta}^{TS}$ are optimal to the LP relaxation of the subproblem $\mbox{{\bf SP-RTS}}(\boldsymbol{y}_n^{TSLP}, u_n^{TSLP})$, and thus $Q^{TLP}(\boldsymbol{y}_n^{TSLP},u_n^{TSLP})=\sum_{n\in \mathcal{T}} p_n\left(\boldsymbol{\tilde{f}}_n^{\mathsf T}\sum_{m\in \mathcal{P}(n)}\boldsymbol{x}^{TS}_m+ \tilde{\lambda}_n \eta^{TS}_n\right)$.

Since $\left\lbrace {\boldsymbol{y}_n^{TSLP}}\right\rbrace_{{n\in\mathcal{T}}},\ \{u_n^{TSLP}\}_{n\in \mathcal{T}\setminus\{1\}}$ are feasible solutions for the multistage problem \eqref{model:multiriskaverse}, we have $
z^{MS}_R\le \sum_{n\in \mathcal{T}} p_n(\boldsymbol{\tilde{c}}_n^{\mathsf T}\boldsymbol{y}^{TSLP}_n + \tilde{\alpha}_n u^{TSLP}_n)+Q^{M}(\boldsymbol{y}_n^{TSLP},u_n^{TSLP})$,
where $Q^{M}(\boldsymbol{y}_n^{TSLP},u_n^{TSLP})$ is the optimal objective value to the substructure problem \eqref{eq:ms-substructure}. The rest of the proof follows the proof of Theorem \ref{thm:risk-averse-bound}.
 
\endproof

}

\proof{Proof of Theorem \ref{thm:risk-averse-upper-bound}}
Because $\{\boldsymbol{y}_n^{MSLP}\}_{n\in \mathcal{T}}, \ \{u_n^{MSLP}\}_{n\in \mathcal{T}\setminus\{1\}}$ constitute an optimal solution to the LP relaxation of the multistage model \eqref{model:multiriskaverse}, we have $
	z^{MS}_R
	%&=&\sum_{n\in \mathcal{T}} p_n(\boldsymbol{\tilde{c}}_n^{\mathsf T}\boldsymbol{y^{TS}_n} + \tilde{\alpha}_n u^{TS}_n)+\sum_{n\in \mathcal{T}} p_n(\boldsymbol{\tilde{f}}_n^{\mathsf T}\sum_{m\in \mathcal{P}(n)}\boldsymbol{x}^{TS}_m + \tilde{\lambda} \eta^{TS}_n)\\
    \ge\sum_{n\in \mathcal{T}} p_n(\boldsymbol{\tilde{c}}_n^{\mathsf T}\boldsymbol{y}^{MSLP}_n + \tilde{\alpha}_n u^{MSLP}_n)+Q^{MLP}(\boldsymbol{y}_n^{MSLP},u_n^{MSLP})$,
where $Q^{MLP}(\boldsymbol{y}_n^{MSLP},u_n^{MSLP})$ stands for the optimal objective value of the LP relaxation of the subproblem $\mbox{{\bf SP-RMS}}(\boldsymbol{y}_n^{MSLP}, u_n^{MSLP})$ \eqref{eq:ms-substructure}. It can be easily verified that the constructed $\boldsymbol{x}^{MS},\boldsymbol{\eta}^{MS}$ are optimal to the LP relaxation of the subproblem, and thus $Q^{MLP}(\boldsymbol{y}_n^{MSLP},u_n^{MSLP})=\sum_{n\in \mathcal{T}} p_n\left(\boldsymbol{\tilde{f}}_n^{\mathsf T}\sum_{m\in \mathcal{P}(n)}\boldsymbol{x}^{MS}_m+ \tilde{\lambda}_n \eta^{MS}_n\right)$.

Since $\left\lbrace {\boldsymbol{y}_n^{MSLP}}\right\rbrace_{{n\in\mathcal{T}}},\ \{u_n^{MSLP}\}_{n\in \mathcal{T}\setminus\{1\}}$ are feasible solutions to the two-stage problem \eqref{model:tworiskaverse}, we have $
z^{TS}_R\le \sum_{n\in \mathcal{T}} p_n(\boldsymbol{\tilde{c}}_n^{\mathsf T}\boldsymbol{y}^{MSLP}_n + \tilde{\alpha}_n u^{MSLP}_n)+Q^{T}(\boldsymbol{y}_n^{MSLP},u_n^{MSLP})$,
where $Q^{T}(\boldsymbol{y}_n^{MSLP},u_n^{MSLP})$ is the optimal objective value to the substructure problem \eqref{eq:ts-substructure}.
Combining all steps above, we have
\begin{align*}
% &{\rm VMS_R}\\
&z^{TS}_R-z^{MS}_R\\
\le& Q^{T}(\boldsymbol{y}_n^{MSLP},u_n^{MSLP})-Q^{MLP}(\boldsymbol{y}_n^{MSLP},u_n^{MSLP})\\
{=}& \sum_{n\in \mathcal{T}} p_n\left(\boldsymbol{\tilde{f}}_n^{\mathsf T}\left(\max_{m\in \mathcal{P}(n)} \lceil \max_{l\in \mathcal{T}_{t_m}}\boldsymbol{B}_{t_l} \boldsymbol{y}_l^{MSLP}\rceil - \max_{m\in \mathcal{P}(n)} \boldsymbol{B}_{t_m} \boldsymbol{y}_m^{MSLP}\right)+ \tilde{\lambda}_n \left(\eta^{TS}_n-\eta^{MS}_n\right)\right)\\
=& \boldsymbol{f}_{t_1}^{\mathsf T}\left(\lceil\boldsymbol{B}_{t_1} \boldsymbol{y}_1^{MSLP}\rceil-\boldsymbol{B}_{t_1} \boldsymbol{y}_1^{MSLP}\right)+\sum_{n\in \mathcal{T}\setminus\{1\}} p_n{(1-\lambda_{t_n})\boldsymbol{f}_{t_n}^{\mathsf T}}\left(\max_{m\in \mathcal{P}(n)} \lceil \max_{l\in \mathcal{T}_{t_m}}\boldsymbol{B}_{t_l} \boldsymbol{y}_l^{MSLP}\rceil - \max_{m\in \mathcal{P}(n)} \boldsymbol{B}_{t_m} \boldsymbol{y}_m^{MSLP}\right)\\
&+ {\sum_{n\in \mathcal{T}\setminus\mathcal{L}} p_n {\lambda}_{t_n+1}} \left(\eta^{TS}_n-\eta^{MS}_n\right).
\end{align*}
This completes the proof.   
\endproof

\proof{Proof of Theorem \ref{thm:hardness}}
We first note that any instance of the NP-hard knapsack problem with $M$ items can be polynomially transformed to the deterministic capacity planning problem \eqref{model:deter} with $T=1$. Moreover, Model \eqref{model:deter} is a single-scenario version of the risk-averse two-stage \eqref{model:tworiskaverse} and multistage \eqref{model:multiriskaverse} counterparts with $\lambda_t=0,\ \forall t=2,\ldots, T$. With that, we prove the NP-hardness of the aforementioned models.
  
\endproof

\proof{Proof of Proposition \ref{prop:alg}}
First of all, at the end of each iteration $k\ge 1$, $(\boldsymbol{x}_n^{k}, \eta_n^{k}, \boldsymbol{y}_n^{k}, u_n^{k})$ constitutes a feasible solution to the risk-averse multistage problem \eqref{model:multiriskaverse} and thus provides an upper bound. We have
\begin{align*}
    & z_{T,R}^{MS}(\boldsymbol{x}_n^{k+1}, \eta_n^{k+1}, \boldsymbol{y}_n^{k+1}, u_n^{k+1}) - z_{T,R}^{MS}(\boldsymbol{x}_n^{k}, \eta_n^{k}, \boldsymbol{y}_n^{k}, u_n^{k})\\
    = & \left(z_{T,R}^{MS}(\boldsymbol{x}_n^{k+1}, \eta_n^{k+1}, \boldsymbol{y}_n^{k+1}, u_n^{k+1}) - 
    z_{T,R}^{MS}(\boldsymbol{x}_n^{k+1}, \eta_n^{k+1}, \boldsymbol{y}_n^{k}, u_n^{k})\right) +
    \left(z_{T,R}^{MS}(\boldsymbol{x}_n^{k+1}, \eta_n^{k+1}, \boldsymbol{y}_n^{k}, u_n^{k}) -
    z_{T,R}^{MS}(\boldsymbol{x}_n^{k}, \eta_n^{k}, \boldsymbol{y}_n^{k}, u_n^{k})\right)\\
    \le & 0.
\end{align*}
The first term $z_{T,R}^{MS}(\boldsymbol{x}_n^{k+1}, \eta_n^{k+1}, \boldsymbol{y}_n^{k+1}, u_n^{k+1}) - 
    z_{T,R}^{MS}(\boldsymbol{x}_n^{k+1}, \eta_n^{k+1}, \boldsymbol{y}_n^{k}, u_n^{k})$ is non-positive because $(\boldsymbol{y}_n^{k+1}, u_n^{k+1})$ is an optimal solution to Step \ref{alg:step5}, while $(\boldsymbol{y}_n^{k}, u_n^{k})$ is a feasible solution to Step \ref{alg:step5} due to Step \ref{alg:step4}; the second term $z_{T,R}^{MS}(\boldsymbol{x}_n^{k+1}, \eta_n^{k+1}, \boldsymbol{y}_n^{k}, u_n^{k}) -
    z_{T,R}^{MS}(\boldsymbol{x}_n^{k}, \eta_n^{k}, \boldsymbol{y}_n^{k}, u_n^{k})$ is non-positive because $(\boldsymbol{x}_n^{k+1}, \eta_n^{k+1})$ is an optimal solution to Step \ref{alg:step4}, while $(\boldsymbol{x}_n^{k}, \eta_n^{k})$ is feasible to Step \ref{alg:step4} due to Step \ref{alg:step5} at the previous iteration. This completes the proof.   
\endproof

\proof{Proof of Proposition \ref{prop:ratio}}
First of all, we have
\begin{align*}
0\le &z_{T,R}^{MS}(\boldsymbol{x}_n^{H}, \eta_n^{H}, \boldsymbol{y}_n^{H}, u_n^{H}) - z_{T,R}^{MS}(\boldsymbol{x}_n^{*}, \eta_n^{*}, \boldsymbol{y}_n^{*}, u_n^{*})\\
\overset{(a)}{\le} & z_{T,R}^{MS}(\boldsymbol{x}_n^{H}, \eta_n^{H}, \boldsymbol{y}_n^{H}, u_n^{H}) - z_{T,R}^{MS}(\boldsymbol{x}_n^{MSLP}, \eta_n^{MSLP}, \boldsymbol{y}_n^{MSLP}, u_n^{MSLP})\\
\overset{(b)}{\le} & z_{T,R}^{MS}(\boldsymbol{x}_n^{1}, \eta_n^{1}, \boldsymbol{y}_n^{1}, u_n^{1}) - z_{T,R}^{MS}(\boldsymbol{x}_n^{MSLP}, \eta_n^{MSLP}, \boldsymbol{y}_n^{MSLP}, u_n^{MSLP})\\
= & z_{T,R}^{MS}(\boldsymbol{x}_n^{1}, \eta_n^{1}, \boldsymbol{y}_n^{1}, u_n^{1}) - 
z_{T,R}^{MS}(\boldsymbol{x}_n^{1}, \eta_n^{1}, \boldsymbol{y}_n^{MSLP}, u_n^{MSLP}) \nonumber\\
& +z_{T,R}^{MS}(\boldsymbol{x}_n^{1}, \eta_n^{1}, \boldsymbol{y}_n^{MSLP}, u_n^{MSLP}) -
z_{T,R}^{MS}(\boldsymbol{x}_n^{MSLP}, \eta_n^{MSLP}, \boldsymbol{y}_n^{MSLP}, u_n^{MSLP})\\
\overset{(c)}{\le} & z_{T,R}^{MS}(\boldsymbol{x}_n^{1}, \eta_n^{1}, \boldsymbol{y}_n^{MSLP}, u_n^{MSLP}) -
z_{T,R}^{MS}(\boldsymbol{x}_n^{MSLP}, \eta_n^{MSLP}, \boldsymbol{y}_n^{MSLP}, u_n^{MSLP})\\
= & Q^M(\boldsymbol{y}_n^{MSLP}, u_n^{MSLP}) - Q^{MLP}(\boldsymbol{y}_n^{MSLP}, u_n^{MSLP})
\end{align*}
where $(a)$ is true because $(\boldsymbol{x}_n^{MSLP}, \eta_n^{MSLP}, \boldsymbol{y}_n^{MSLP}, u_n^{MSLP})$ provides an optimal solution to the LP relaxation problem, $(b)$ is true because of Proposition \ref{prop:alg}, and $(c)$ is true because $(\boldsymbol{y}_n^{1}, u_n^{1})$ is an optimal solution to Step \ref{alg:step5} in iteration 0 while $(\boldsymbol{y}_n^{MSLP}, u_n^{MSLP})$ is feasible. Here, $Q^M(\boldsymbol{y}_n^{MSLP}, u_n^{MSLP})$ is defined as the optimal objective value to the subproblem $\mbox{{\bf SP-RMS}}(\boldsymbol{y}_n^{MSLP}, u_n^{MSLP})$ \eqref{eq:ms-substructure}, and $Q^{MLP}(\boldsymbol{y}_n^{MSLP}, u_n^{MSLP})$ stands for the optimal objective value of its LP relaxation.

Based on Proposition \ref{prop:substructure}, we have a closed form for $Q^M(\boldsymbol{y}_n^{MSLP}, u_n^{MSLP})$ and we can apply the same analysis to obtain a closed form for $Q^{MLP}(\boldsymbol{y}_n^{MSLP}, u_n^{MSLP})$, where we replace all $\lceil \boldsymbol{B}_{t_n} \boldsymbol{y}_n^{MSLP}\rceil$ with $\boldsymbol{B}_{t_n} \boldsymbol{y}_n^{MSLP}$. As a result, 
\begin{align*}
    & Q^M(\boldsymbol{y}_n^{MSLP}, u_n^{MSLP}) - Q^{MLP}(\boldsymbol{y}_n^{MSLP}, u_n^{MSLP})\\
= & \sum_{n\in \mathcal{T}} p_n\left(\boldsymbol{\tilde{f}}_n^{\mathsf T}\max_{m\in\mathcal{P}(n)}\lceil \boldsymbol{B}_{t_m} \boldsymbol{y}_m^{MSLP}\rceil+ \tilde{\lambda}_n \max_{m\in \mathcal{C}(n)}\left\lbrace \boldsymbol{f}^{\mathsf T}_{t_m}\max_{l\in\mathcal{P}(m)}\lceil \boldsymbol{B}_{t_l} \boldsymbol{y}_l^{MSLP}\rceil+ \boldsymbol{c}_{t_m}^{\mathsf T} \boldsymbol{y}_m^{MSLP}-u_m^{MSLP}\right\rbrace\right) - \\
&\sum_{n\in \mathcal{T}} p_n\left(\boldsymbol{\tilde{f}}_n^{\mathsf T}\max_{m\in\mathcal{P}(n)} \boldsymbol{B}_{t_m} \boldsymbol{y}_m^{MSLP}+ \tilde{\lambda}_n \max_{m\in \mathcal{C}(n)}\left\lbrace \boldsymbol{f}^{\mathsf T}_{t_m}\max_{l\in\mathcal{P}(m)} \boldsymbol{B}_{t_l} \boldsymbol{y}_l^{MSLP}+ \boldsymbol{c}_{t_m}^{\mathsf T} \boldsymbol{y}_m^{MSLP}-u_m^{MSLP}\right\rbrace\right)\\
= & \sum_{n\in \mathcal{T}} p_n\boldsymbol{\tilde{f}}_n^{\mathsf T}\left(\max_{m\in\mathcal{P}(n)} \lceil \boldsymbol{B}_{t_m} \boldsymbol{y}_m^{MSLP}\rceil - \max_{m\in\mathcal{P}(n)} \boldsymbol{B}_{t_m} \boldsymbol{y}_m^{MSLP}\right) \\
&+ \sum_{n\in \mathcal{T}} p_n \tilde{\lambda}_n\Big(\max_{m\in \mathcal{C}(n)}\left\lbrace \boldsymbol{f}^{\mathsf T}_{t_m}\max_{l\in\mathcal{P}(m)}\lceil \boldsymbol{B}_{t_l} \boldsymbol{y}_l^{MSLP}\rceil+ \boldsymbol{c}_{t_m}^{\mathsf T} \boldsymbol{y}_m^{MSLP}-u_m^{MSLP}\right\rbrace\\
&\hspace{6em}- \max_{m\in \mathcal{C}(n)}\left\lbrace \boldsymbol{f}^{\mathsf T}_{t_m}\max_{l\in\mathcal{P}(m)} \boldsymbol{B}_{t_l} \boldsymbol{y}_l^{MSLP}+ \boldsymbol{c}_{t_m}^{\mathsf T} \boldsymbol{y}_m^{MSLP}-u_m^{MSLP}\right\rbrace \Big)\\
\overset{(a)}{\le} & \sum_{n\in\mathcal{T}}p_n\boldsymbol{\tilde{f}}_n^{\mathsf T}\boldsymbol{1} + \sum_{n\in\mathcal{T}}p_n\tilde{\lambda}_n\boldsymbol{f}^{\mathsf T}_{t_n+1} \boldsymbol{1}\\
= & {\sum_{i=1}^Mf_{1i} + \sum_{t=2}^T(1-\lambda_t)\sum_{i=1}^Mf_{ti} + \sum_{t=2}^T\lambda_t(\sum_{i=1}^Mf_{ti})}\\
= & {\sum_{t=1}^T\sum_{i=1}^M f_{ti}}
\end{align*}
where $(a)$ is true because $\max_{m\in\mathcal{P}(n)}\lceil \boldsymbol{B}_{t_m} \boldsymbol{y}_m^{MSLP}\rceil - \max_{m\in\mathcal{P}(n)} \boldsymbol{B}_{t_m} \boldsymbol{y}_m^{MSLP}\le \boldsymbol{1}$ and by letting $m^*(n) = \arg\max_{m\in \mathcal{C}(n)}\left\lbrace \boldsymbol{f}^{\mathsf T}_{t_m}\max_{l\in\mathcal{P}(m)}\lceil \boldsymbol{B}_{t_m} \boldsymbol{y}_l^{MSLP}\rceil+ \boldsymbol{c}_{t_m}^{\mathsf T} \boldsymbol{y}_m^{MSLP}-u_m^{MSLP}\right\rbrace$, we have 
\small
\begin{align*}
&\max_{m\in \mathcal{C}(n)}\left\lbrace \boldsymbol{f}^{\mathsf T}_{t_m}\max_{l\in\mathcal{P}(m)}\lceil \boldsymbol{B}_{t_l} \boldsymbol{y}_l^{MSLP}\rceil+ \boldsymbol{c}_{t_m}^{\mathsf T} \boldsymbol{y}_m^{MSLP}-u_m^{MSLP}\right\rbrace - \max_{m\in \mathcal{C}(n)}\left\lbrace \boldsymbol{f}^{\mathsf T}_{t_m}\max_{l\in\mathcal{P}(m)} \boldsymbol{B}_{t_l} \boldsymbol{y}_l^{MSLP}+ \boldsymbol{c}_{t_m}^{\mathsf T} \boldsymbol{y}_m^{MSLP}-u_m^{MSLP}\right\rbrace \\
\le 
&\left( \boldsymbol{f}^{\mathsf T}_{t_{n}+1}\max_{l\in\mathcal{P}(m^*(n))}\lceil \boldsymbol{B}_{t_l} \boldsymbol{y}_l^{MSLP}\rceil+ \boldsymbol{c}_{m^*(n)}^{\mathsf T} \boldsymbol{y}_{m^*(n)}^{MSLP}-u_{m^*(n)}^{MSLP}\right) - \left( \boldsymbol{f}^{\mathsf T}_{t_{n}+1}\max_{l\in\mathcal{P}(m^*(n))} \boldsymbol{B}_{t_l} \boldsymbol{y}_l^{MSLP}+ \boldsymbol{c}_{m^*(n)}^{\mathsf T} \boldsymbol{y}_{m^*(n)}^{MSLP}-u_{m^*(n)}^{MSLP}\right)\\
= & \boldsymbol{f}^{\mathsf T}_{t_n+1} \left(\max_{l\in\mathcal{P}(m^*(n))}\lceil \boldsymbol{B}_{t_l} \boldsymbol{y}_l^{MSLP}\rceil - \max_{l\in\mathcal{P}(m^*(n))} \boldsymbol{B}_{t_l} \boldsymbol{y}_l^{MSLP}\right)\\
\le & \boldsymbol{f}^{\mathsf T}_{t_n+1} \boldsymbol{1}.
\end{align*} 
\normalsize
This completes the proof.

% \begin{repeattheorem}[Theorem \ref{thm:ratio}]
% Algorithm \ref{alg:approx-multistage} has an approximation ratio of $$1+ {\frac{M\sum_{t=1}^Tf_{t,\rm max} }{M_{\rm min}\sum_{t=1}^Tf_{t,\rm min}+\sum_{t=1}^Tc_{t,\rm min}\min_{n\in\mathcal{T}_t}\{\sum_{j=1}^Nd_{n,j}\}}},$$ where $h_{\rm max} = \max_{i=1}^M\{h_{1i}\},\ f_{t, \rm max} = \max_{i=1}^M\{f_{ti}\},\ f_{t,\rm min} = \min_{i=1}^M\{f_{ti}\},\ {c_{t,\rm min}=\min_{i\in[M],j\in[N]}c_{tij}}$ and $M_{\rm min} = \lceil\frac{\sum_{j=1}^N d_{1j}}{h_{\rm max}}\rceil$ measures at least how many facilities we need to cover the first-stage demand. 
% \end{repeattheorem}

\proof{Proof of Theorem \ref{thm:ratio}}
Based on Proposition \ref{prop:ratio}, we have
\begin{align*}
    \frac{z_{T,R}^{MS}(\boldsymbol{x}_n^{H}, \eta_n^{H}, \boldsymbol{y}_n^{H}, u_n^{H})}{z_{T,R}^{MS}(\boldsymbol{x}_n^{*}, \eta_n^{*}, \boldsymbol{y}_n^{*}, u_n^{*})} 
    \le& \frac{z_{T,R}^{MS}(\boldsymbol{x}_n^{*}, \eta_n^{*}, \boldsymbol{y}_n^{*}, u_n^{*}) + {\sum_{t=1}^T\sum_{i=1}^M f_{ti}}}{z_{T,R}^{MS}(\boldsymbol{x}_n^{*}, \eta_n^{*}, \boldsymbol{y}_n^{*}, u_n^{*})} \\
    = &1+ \frac{{\sum_{t=1}^T\sum_{i=1}^M f_{ti}}}{z_{T,R}^{MS}(\boldsymbol{x}_n^{*}, \eta_n^{*}, \boldsymbol{y}_n^{*}, u_n^{*})}\\
    \overset{(a)}{\le} &1+ {\frac{M\sum_{t=1}^Tf_{t,\rm max}}{M_{\rm min}\sum_{t=1}^Tf_{t,\rm min}+\sum_{t=1}^Tc_{t,\rm min}\min_{n\in\mathcal{T}_t}\{\sum_{j=1}^Nd_{n,j}\}}}.
\end{align*} 
Next, we show that $(a)$ is true (i.e., $z_{T,R}^{MS}(\boldsymbol{x}_n^{*}, \eta_n^{*}, \boldsymbol{y}_n^{*}, u_n^{*})\ge M_{\rm min}\sum_{t=1}^Tf_{t,\rm min}{+\sum_{t=1}^Tc_{t,\rm min}\min_{n\in\mathcal{T}_t}\{\sum_{j=1}^Nd_{n,j}\}}$). 
As mentioned before, the optimal solution of $\eta_n$ is attained at $\eta^*_n = {\rm VaR}_{\alpha_{t_n+1}}[\boldsymbol{{f}}^{\mathsf T}_{t_m}\sum_{l\in \mathcal{P}(m)}\boldsymbol{x}_l^*+\boldsymbol{{c}}_m^{\mathsf T}\boldsymbol{y}_m^*: m\in\mathcal{C}(n)]$ by the definition of $\rm CVaR$ \eqref{eq:ccvar}.
Our analysis then follows 
\begin{align*}
    z^{MS}_{T,R}(\boldsymbol{x}_n^{*}, \eta_n^{*}, \boldsymbol{y}_n^{*}, u_n^{*})=&\sum_{n\in \mathcal{T}} p_n(\boldsymbol{\tilde{f}}_n^{\mathsf T}\sum_{m\in \mathcal{P}(n)}\boldsymbol{x}_m^*+\boldsymbol{\tilde{c}}_n^{\mathsf T}\boldsymbol{y}_n^* + \tilde{\lambda}_n \eta_n^* + \tilde{\alpha}_n u_n^*)\\
    \ge& \sum_{n\in \mathcal{T}} p_n(\boldsymbol{\tilde{f}}_n^{\mathsf T}\sum_{m\in \mathcal{P}(n)}\boldsymbol{x}_m^*+\boldsymbol{\tilde{c}}_{n}^{\mathsf T}\boldsymbol{y}_n^* + \tilde{\lambda}_n \eta_n^* )\\
    =& \boldsymbol{{f}}^{\mathsf T}_{t_1}\boldsymbol{x}_1^* + \boldsymbol{{c}}_{t_1}^{\mathsf T}\boldsymbol{y}_1^* + \sum_{n\not=1}(1-\lambda_{t_n}) p_n(\boldsymbol{{f}}^{\mathsf T}_{t_n}\sum_{m\in \mathcal{P}(n)}\boldsymbol{x}_m^* + \boldsymbol{{c}}_{t_n}^{\mathsf T}\boldsymbol{y}_n^*) \\
    &+ \sum_{n\not\in\mathcal{L}} \lambda_{t_n+1}p_n {\rm VaR}_{\alpha_{t_n+1}}[\boldsymbol{{f}}^{\mathsf T}_{t_m}\sum_{l\in \mathcal{P}(m)}\boldsymbol{x}_l^*+\boldsymbol{{c}}_{t_m}^{\mathsf T}\boldsymbol{y}_m^*: m\in\mathcal{C}(n)]\\
    \overset{(a)}{\ge} & \boldsymbol{{f}}^{\mathsf T}_{t_1}\boldsymbol{x}_1^* + \boldsymbol{{c}}_{t_1}^{\mathsf T}\boldsymbol{y}_1^* + \sum_{t=2}^T\min_{n\in\mathcal{T}_t}\{\boldsymbol{{f}}^{\mathsf T}_{t_n}\sum_{m\in \mathcal{P}(n)}\boldsymbol{x}_m^* + \boldsymbol{{c}}_{t_n}^{\mathsf T}\boldsymbol{y}_n^*\}\sum_{n\in\mathcal{T}_t}(1-\lambda_{t_n})p_n \\
    & + \sum_{t=1}^{T-1}\min_{n\in\mathcal{T}_{t+1}}\{\boldsymbol{{f}}^{\mathsf T}_{t_n}\sum_{m\in \mathcal{P}(n)}\boldsymbol{x}_m^* + \boldsymbol{{c}}_{t_n}^{\mathsf T}\boldsymbol{y}_n^*\}\sum_{n\in\mathcal{T}_t}\lambda_{t_n+1}p_n \\
    \overset{(b)}{=} & \boldsymbol{{f}}^{\mathsf T}_{t_1}\boldsymbol{x}_1^* + \boldsymbol{{c}}_{t_1}^{\mathsf T}\boldsymbol{y}_1^* + \sum_{t=2}^T \min_{n\in\mathcal{T}_t}\{\boldsymbol{{f}}^{\mathsf T}_{t_n}\sum_{m\in \mathcal{P}(n)}\boldsymbol{x}_m^* + \boldsymbol{{c}}_{t_n}^{\mathsf T}\boldsymbol{y}_n^*\}\\
    \ge & \sum_{t=1}^T \min_{n\in\mathcal{T}_t}\{\boldsymbol{{f}}^{\mathsf T}_{t_n}\sum_{m\in \mathcal{P}(n)}\boldsymbol{x}_m^*\} {+\sum_{t=1}^T\min_{n\in\mathcal{T}_t}\{\boldsymbol{{c}}_{t_n}^{\mathsf T}\boldsymbol{y}_n^*}\}\\
    %\overset{(c)}{\ge} & \sum_{t=1}^T \boldsymbol{{f}}^{\mathsf T}_{t}\boldsymbol{x}_1^*{+\sum_{t=1}^Tc_{t,\rm min}\min_{n\in\mathcal{T}_t}\{\sum_{j=1}^Nd_{n,j}\}}\\
    \overset{(c)}{\ge} & \sum_{t=1}^Tf_{t, \rm min}\sum_{i=1}^M x_{1i}^*{+\sum_{t=1}^Tc_{t,\rm min}\min_{n\in\mathcal{T}_t}\{\sum_{j=1}^Nd_{n,j}\}}\\
    \overset{(d)}{\ge} & \sum_{t=1}^Tf_{t,\rm min}M_{\rm min}{+\sum_{t=1}^Tc_{t,\rm min}\min_{n\in\mathcal{T}_t}\{\sum_{j=1}^Nd_{n,j}\}},
\end{align*}
where $(a)$ is true because $\boldsymbol{{f}}^{\mathsf T}_{t_n}\sum_{m\in \mathcal{P}(n)}\boldsymbol{x}_m^* + \boldsymbol{{c}}_{t_n}^{\mathsf T}\boldsymbol{y}_n^*\ge \min_{m\in\mathcal{T}_{t_n}}\{\boldsymbol{{f}}^{\mathsf T}_{t_m}\sum_{l\in \mathcal{P}(m)}\boldsymbol{x}_l^* + \boldsymbol{{c}}_{t_m}^{\mathsf T}\boldsymbol{y}_m^*\},\ \forall n\in\mathcal{T}$ and ${\rm VaR}_{\alpha_{t_n+1}}[\boldsymbol{{f}}^{\mathsf T}_{t_m}\sum_{l\in \mathcal{P}(m)}\boldsymbol{x}_l^*+\boldsymbol{{c}}_{t_m}^{\mathsf T}\boldsymbol{y}_m^*: m\in\mathcal{C}(n)]\ge \min_{n\in\mathcal{T}_{t_n+1}}\{\boldsymbol{{f}}^{\mathsf T}_{t_n}\sum_{m\in \mathcal{P}(n)}\boldsymbol{x}_m^* + \boldsymbol{{c}}_{t_n}^{\mathsf T}\boldsymbol{y}_n^*\},\ \forall n\not\in\mathcal{L}$; $(b)$ is true because $\sum_{n\in\mathcal{T}_t}(1-\lambda_{t_n})p_n = (1-\lambda_t)\sum_{n\in\mathcal{T}_t}p_n = 1-\lambda_t$ and $\sum_{n\in\mathcal{T}_t}\lambda_{t_n+1}p_n = \lambda_{t+1}\sum_{n\in\mathcal{T}_t}p_n = \lambda_{t+1}$; $(c)$ is true because $\sum_{m\in \mathcal{P}(n)}\boldsymbol{x}_m^*\ge \boldsymbol{x}_1^*$ and {$\boldsymbol{{c}}_{t_n}^{\mathsf T}\boldsymbol{y}_n^*\ge c_{t,\rm min}\sum_{i=1}^M\sum_{j=1}^N{y}_{n,ij}^*=c_{t,\rm min}\sum_{j=1}^Nd_{n,j}$}; $(d)$ is true because $h_{\rm max}\sum_{i=1}^Mx_{1i}^* \ge \sum_{i=1}^Mh_{1i}x_{1i}^* \ge \sum_{i=1}^M\sum_{j=1}^Ny_{1ij}^* = \sum_{j=1}^Nd_{1j}$ from constraints \eqref{eq:2} and \eqref{eq:3}. This completes the proof.    
\endproof

{
\proof{Proof of Corollary \ref{cor:ratio}}
First of all, $\frac{z_{T,R}^{MS}(\boldsymbol{x}_n^{H}, \eta_n^{H}, \boldsymbol{y}_n^{H}, u_n^{H})}{z_{T,R}^{MS}(\boldsymbol{x}_n^{*}, \eta_n^{*}, \boldsymbol{y}_n^{*}, u_n^{*})}\ge 1$. According to Theorem \ref{thm:ratio}, we have 
\begin{align*}
   1\le & \frac{z_{T,R}^{MS}(\boldsymbol{x}_n^{H}, \eta_n^{H}, \boldsymbol{y}_n^{H}, u_n^{H})}{z_{T,R}^{MS}(\boldsymbol{x}_n^{*}, \eta_n^{*}, \boldsymbol{y}_n^{*}, u_n^{*})} \\
{\le} &1+ {\frac{M\sum_{t=1}^Tf_{t,\rm max} }{M_{\rm min}\sum_{t=1}^Tf_{t,\rm min}+\sum_{t=1}^Tc_{t,\rm min}\min_{n\in\mathcal{T}_t}\{\sum_{j=1}^Nd_{n,j}\}}}\\
 =&1+\frac{\sum_{t=1}^TO(1)}{\sum_{t=1}^TO(1)+\sum_{t=1}^TO(t)}\\
 =&1+ \frac{O(T)}{O(T^2)}\to 1\ (\text{when }T\to\infty).
\end{align*} 
This completes the proof.
 
\endproof

}

{
\section{\color{black}An Example to Illustrate the Usefulness of the Derived Bounds}\label{e-companion:example}
\normalsize
In this appendix, we provide computational details of Example \ref{ex:tight} to illustrate the usefulness of lower bound ${\rm VMS_R^{LB}}$ and upper bound ${\rm VMS_R^{UB}}$ presented in Section \ref{sec:VMS}.
    Here, we consider one facility with $h_1=50$ and one customer site, where the customer demand $d_1=0$ in node 1, $d_2=50$ in node 2, and $d_3=150$ in node 3 (see Figure \ref{fig:eg}). In a two-stage model, to cover the demand in both nodes, one must invest $x_2=x_3=3$ for both nodes 2 and 3 up front. On the contrary, in a multistage model, the decision maker has the chance to wait and see the revealed demand and then make investment decisions, and as a result, the optimal solution is $x_2=1,\ x_3=3$, which obviously costs less than the two-stage counterpart.

 Because $d_1=0$ when $t=1$, all costs will occur at $t=2$ and thus we only analyze the cost at $t=2$. Denote $f$ as the unit maintenance cost for one stage and $c$ as the operational cost, respectively. The optimal resource-allocation decisions at the second stage for both two-stage and multistage models are $y^{*}_2 = 50,\ y^{*}_{3}=150$. Next, we calculate the optimal objective values of the risk-averse two-stage and multistage models while assuming $\alpha=0.5$ as follows:
    \begin{itemize}
        \item Two-stage model \eqref{model:tworiskaverse}:
        \begin{itemize}
            \item Node 2: $g_2={f}{x}^{*}_2+{c}{y}^{*}_2= 3f + 50c$ with probability 0.5;
            \item Node 3: $g_3={f}{x}^{*}_3+{c}{y}^{*}_3= 3f + 150c$ with probability 0.5.
        \end{itemize}
        Because $\alpha = 0.5$, we have $\eta^*={\rm VaR_{\alpha}}[g] = 3f + 50c,\ {\rm CVaR_{\alpha}}[g] = {\rm VaR_{\alpha}}[g] + \frac{1}{1-\alpha}\mathbb{E}[g - {\rm VaR_{\alpha}}[g]]_+ = 3f + 50c + \frac{1}{0.5}[0.5 \times 100c] = 3f + 150c$ and correspondingly $u^{*}_2=0,\ u^{*}_3=100c$. Therefore, $z_R^{TS} = \rho(g) = (1-\lambda)\mathbb{E}[g] + \lambda {\rm CVaR}(g) = (1-\lambda)\frac{6f+200c}{2} + \lambda(3f + 150c)=3f+100c+50\lambda c$.
        \item Multistage model \eqref{model:multiriskaverse}:
            \begin{itemize}
            \item Node 2: $g_2={f}{x}^{MS}_2+{c}{y}^{MS}_2 = f + 50c$ with probability 0.5;
            \item Node 3: $g_3={f}{x}^{MS}_3+{c}{y}^{MS}_3 = 3f + 150c$ with probability 0.5.
        \end{itemize}
        Because $\alpha = 0.5$, we have ${\rm VaR_{\alpha}}[g] = f + 50c$ and ${\rm CVaR_{\alpha}}[g] = {\rm VaR_{\alpha}}[g] + \frac{1}{1-\alpha}\mathbb{E}[g - {\rm VaR_{\alpha}}[g]]_+ = f + 50c + \frac{1}{0.5}[0.5 \times (2f + 100c)] = 3f + 150c$. Therefore, $z_R^{MS} = \rho(g) = (1-\lambda)\mathbb{E}[g] + \lambda {\rm CVaR}(g) = (1-\lambda)\frac{4f+200c}{2} + \lambda(3f + 150c)=2f+100c+\lambda f+50\lambda c$.

    \end{itemize}
    
    As a result, the gap between the optimal objective values of the two-stage and multistage models only exists in the maintenance cost, which is ${\rm VMS_R} = z_R^{TS} - z_R^{MS} = (1-\lambda)f$. Based on Proposition \ref{prop:substructure}, the constructed $\eta$-solutions for the root node are $\eta^{TS}=\max\{3f + 50c-0, 3f + 150c - 100c\} = 3f + 50c$ and $\eta^{MS}=\max\{f + 50c - 0, 3f + 150c - 100c\} = 3f + 50c = \eta^{TS}$ and thus the lower bound ${\rm VMS_R^{LB}}$ only depends on the {\color{black}variability} of units of resources $\boldsymbol{B} \boldsymbol{y}^{*}=\frac{y^*}{h}$ across different scenarios. We notice that the number of resources $\frac{y^*_2}{h} = 1$ at node 2, and $\frac{y^*_3}{h} = 3$ at node 3, which leads to ${\rm VMS_R^{LB}} = 0.5 \times (1-\lambda) \times f \times (3 - 1) =(1-\lambda)f = {\rm VMS_R}$. {\color{black}On the other hand, the optimal solution to the LP relaxation of the multistage model is $x_2^{MSLP}=1,\ x_3^{MSLP}=3,\ y_2^{MSLP}=50,\ y_3^{MSLP}=150,\ \eta^{MSLP}=f+50c,\ u_2^{MSLP}=0,\ u_3^{MSLP}=(3f+150c)-(f+50c)=2f+100c$. According to Theorem \ref{thm:risk-averse-upper-bound}, the constructed $\eta$-solutions for the root node are $\eta^{TS}=\max\{3f+50c-0,3f+150c-(2f+100c)\}=3f+50c$ and $\eta^{MS}=\max\{f+50c-0, 3f+150c-(2f+100c)\}=f+50c$. As a result, ${\rm VMS_R^{UB}}=0.5\times(1-\lambda)\times 2f+\lambda\times2f=(1+\lambda)f$.}

    This example also provides a special case where both ${\rm VMS_R^{LB}}$ and ${\rm VMS_R}$ are positively related to the maintenance cost and negatively impacted by the risk attitude $\lambda$.

}

{\color{black}
\section{Additional Computational Results}\label{e-companion:results}
This section presents additional computational results. We record the percentages of instances in different cases with SD scenario trees under varying risk attitude in Table \ref{tab:percentage-lambda}, under varying demand standard deviation in Table \ref{tab:percentage-sigma}, and the percentages with SD scenario trees under different parameter settings when $\delta_1 = 5\%,\ \delta_2=20\%$ in Tables \ref{tab:percentage-deltaC}--\ref{tab:percentage-deltadeviation}. Finally, we present computational time comparison with SD scenario trees and different parameter settings in Tables \ref{tab:time-CSAA}--\ref{tab:time-NSAA}, respectively.
% Table generated by Excel2LaTeX from sheet 'Sheet1'
\begin{table}[ht!]
% \OneAndAHalfSpacedXI
  \centering
  \caption{Percentage of instances in different cases with SD scenario trees and varying risk attitude $\lambda$ when $\delta_1 = 10\%,\ \delta_2=30\%$}
    \begin{tabular}{rrrr}
    \hline
    $\lambda$ & Case (i) & Case (ii) & Case (iii) \\
    \hline
    0     & 52/100 & 46/100 & 2/100 \\
    0.2   & 37/100 & 55/100 & 8/100 \\
    0.4   & 29/100 & 67/100 & 4/100 \\
    0.6   & 17/100 & 77/100 & 6/100 \\
    0.8   & 5/100 & 92/100 & 3/100 \\
    1     & 1/100 & 93/100 & 6/100 \\
    \hline
    \end{tabular}%
  \label{tab:percentage-lambda}%
\end{table}%

% Table generated by Excel2LaTeX from sheet 'Sheet1'
\begin{table}[ht!]
% \OneAndAHalfSpacedXI
  \centering
  \caption{Percentage of instances in different cases with SD scenario trees and varying demand standard deviation $\sigma$ when $\delta_1 = 10\%,\ \delta_2=30\%$}
    \begin{tabular}{rrrr}
    \hline
    $\sigma$ & Case (i) & Case (ii) & Case (iii) \\
    \hline
    0.2   & 0/100 & 100/100 & 0/100 \\
    0.4   & 3/100 & 97/100 & 0/100 \\
    0.6   & 13/100 & 84/100 & 3/100 \\
    0.8   & 20/100 & 73/100 & 7/100 \\
    \hline
    \end{tabular}%
  \label{tab:percentage-sigma}%
\end{table}%

\begin{table}[ht!]
% \OneAndAHalfSpacedXI
  \centering
  \caption{Percentage of instances in different cases with SD scenario trees and varying number of branches $C$ when $\delta_1 = 5\%,\ \delta_2=20\%$}
    \begin{tabular}{rrrr}
    \hline
    $C$ & Case (i) & Case (ii) & Case (iii) \\
    \hline
    2     & 66/100 & 17/100 & 17/100 \\
    3     & 96/100 & 1/100 & 3/100 \\
    4     & 100/100 & 0/100 & 0/100 \\
    5     & 100/100 & 0/100 & 0/100 \\
    \hline
    \end{tabular}%
  \label{tab:percentage-deltaC}%
\end{table}%

% Table generated by Excel2LaTeX from sheet 'Sheet1'
\begin{table}[ht!]
% \OneAndAHalfSpacedXI
  \centering
  \caption{Percentage of instances in different cases with SD scenario trees and varying number of stages $T$ when $\delta_1 = 5\%,\ \delta_2=20\%$}
    \begin{tabular}{rrrr}
     \hline
    $T$ & Case (i) & Case (ii) & Case (iii) \\
    \hline
    3     & 66/100 & 17/100 & 17/100 \\
    4     & 94/100 & 0/100 & 6/100 \\
    5     & 100/100 & 0/100 & 0/100 \\
    6     & 100/100 & 0/100 & 0/100 \\
    7     & 100/100 & 0/100 & 0/100 \\
    8     & 100/100 & 0/100 & 0/100 \\
    \hline
    \end{tabular}%
  \label{tab:percentage-deltaT}%
\end{table}%

% Table generated by Excel2LaTeX from sheet 'Sheet1'
\begin{table}[ht!]
% \OneAndAHalfSpacedXI
  \centering
  \caption{Percentage of instances in different cases with SD scenario trees and varying risk attitude $\lambda$ when $\delta_1 = 5\%,\ \delta_2=20\%$}
    \begin{tabular}{rrrr}
    \hline
    $\lambda$ & Case (i) & Case (ii) & Case (iii) \\
    \hline
    0     & 90/100 & 4/100 & 6/100 \\
    0.2   & 84/100 & 8/100 & 8/100 \\
    0.4   & 74/100 & 12/100 & 14/100 \\
    0.6   & 59/100 & 20/100 & 21/100 \\
    0.8   & 43/100 & 34/100 & 23/100 \\
    1     & 26/100 & 48/100 & 26/100 \\
    \hline
    \end{tabular}%
  \label{tab:percentage-deltarisk}%
\end{table}%

% Table generated by Excel2LaTeX from sheet 'Sheet1'
\begin{table}[ht!]
% \OneAndAHalfSpacedXI
  \centering
  \caption{Percentage of instances in different cases with SD scenario trees and varying demand standard deviation $\sigma$ when $\delta_1 = 5\%,\ \delta_2=20\%$}
    \begin{tabular}{rrrr}
    \hline
    $\sigma$ & Case (i) & Case (ii) & Case (iii) \\
    \hline
    0.2   & 8/100 & 91/100 & 1/100 \\
    0.4   & 46/100 & 47/100 & 7/100 \\
    0.6   & 67/100 & 19/100 & 14/100 \\
    0.8   & 66/100 & 17/100 & 17/100 \\
    \hline
    \end{tabular}%
  \label{tab:percentage-deltadeviation}%
\end{table}%

% Table generated by Excel2LaTeX from sheet 'Sheet3'
\begin{table}[ht!]
% \OneAndAHalfSpacedXII
  \centering
  \caption{Computational time comparison with SD scenario trees and different number of branches $C$}
  \resizebox{0.8\textwidth}{!}{
    \begin{tabular}{r|rr|rr|rr}
    \hline
          & \multicolumn{2}{c|}{$z_R^{TS}$ via Gurobi} & \multicolumn{2}{c|}{$z_R^{MS}$ via Gurobi} & \multicolumn{2}{c}{$z_R^{MS}$ via AA} \\
    $C$     & Time (sec.)  & Obj. (\$) & Time (sec.)  & Obj. (\$) & Time (sec.)  & Obj. (\$) \\
    \hline
    2     & 2.78  & 20,051K & 1800 (0.05\%) & 18,752K & 0.10  & 18,864K \\
    3     & 1.26  & 21,684K & 1800 (0.07\%) & 19,640K & 0.24  & 19,999K \\
    4     & 3.67  & 27,488K & 1800 (0.13\%) & 23,569K & 0.45  & 23,751K \\
    5     & 1800 (0.03\%) & 22,915K & 1800 (0.06\%) & 18,464K & 1.09  & 18,741K \\
    \hline
    \end{tabular}%
    }
  \label{tab:time-CSAA}%
\end{table}%

% Table generated by Excel2LaTeX from sheet 'Sheet3'
\begin{table}[ht!]
% \OneAndAHalfSpacedXII
  \centering
  \caption{Computational time comparison with SD scenario trees and different number of stages $T$}
  \resizebox{0.8\textwidth}{!}{
    \begin{tabular}{r|rr|rr|rr}
    \hline
          & \multicolumn{2}{c|}{$z_R^{TS}$ via Gurobi} & \multicolumn{2}{c|}{$z_R^{MS}$ via Gurobi} & \multicolumn{2}{c}{$z_R^{MS}$ via AA} \\
    $T$     & Time (sec.)  & Obj. (\$) & Time (sec.)  & Obj. (\$) & Time (sec.)  & Obj. (\$) \\
    \hline
    3     & 2.78  & 20,051K & 1800 (0.05\%) & 18,752K & 0.10  & 18,864K \\
    4     & 1.84  & 40,836K & 1800 (0.06\%) & 36,445K & 0.24  & 36,819K \\
    5     & 1800 (0.06\%) & 62,020K & 1800 (0.06\%) & 53,045K & 0.96  & 53,279K \\
    6     & 1800 (0.05\%) & 101,280K & 1800 (0.12\%) & 82,199K & 2.25  & 82,970K \\
    \hline
    \end{tabular}%
    }
  \label{tab:time-TSAA}%
\end{table}%

% Table generated by Excel2LaTeX from sheet 'Sheet3'
\begin{table}[ht!]
% \OneAndAHalfSpacedXII
  \centering
  \caption{Computational time comparison with SD scenario trees and different number of facilities $M$}
  \resizebox{0.8\textwidth}{!}{
    \begin{tabular}{r|rr|rr|rr}
    \hline
          & \multicolumn{2}{c|}{$z_R^{TS}$ via Gurobi} & \multicolumn{2}{c|}{$z_R^{MS}$ via Gurobi} & \multicolumn{2}{c}{$z_R^{MS}$ via AA} \\
    $M$     & Time (sec.)  & Obj. (\$) & Time (sec.)  & Obj. (\$) & Time (sec.)  & Obj. (\$) \\
    \hline
    5     & 2.78  & 20,051K & 1800 (0.05\%) & 18,752K & 0.10  & 18,864K \\
    10    & 1.66  & 18,401K & 1800 (0.01\%) & 16,675K & 0.17  & 17,087K \\
    15    & 3.58  & 20,186K & 1800 (0.04\%) & 17,912K & 0.51  & 18,325K \\
    20    & 3.11  & 21,632K & 1800 (0.06\%) & 19,481K & 0.35  & 20,215K \\
    \hline
    \end{tabular}%
    }
  \label{tab:time-MSAA}%
\end{table}%

% Table generated by Excel2LaTeX from sheet 'Sheet3'
\begin{table}[ht!]
% \OneAndAHalfSpacedXII
  \centering
  \caption{Computational time comparison with SD scenario trees and different number of customer sites $N$}
  \resizebox{0.8\textwidth}{!}{
    \begin{tabular}{r|rr|rr|rr}
    \hline
          & \multicolumn{2}{c|}{$z_R^{TS}$ via Gurobi} & \multicolumn{2}{c|}{$z_R^{MS}$ via Gurobi} & \multicolumn{2}{c}{$z_R^{MS}$ via AA} \\
    $N$     & Time (sec.)  & Obj. (\$) & Time (sec.)  & Obj. (\$) & Time (sec.)  & Obj. (\$) \\
    \hline
    10    & 2.78  & 20,051K & 1800 (0.05\%) & 18,752K & 0.10  & 18,864K \\
    20    & 14.33 & 35,191K & 1800 (0.03\%) & 34,521K & 0.12  & 34,873K \\
    30    & 1800 (0.03\%) & 60,628K & 1800 (0.02\%) & 58,157K & 0.23  & 58,427K \\
    40    & 1800 (0.05\%) & 81,626K & 1800 (0.07\%) & 74,741K & 0.28  & 74,936K \\
    \hline
    \end{tabular}%
    }
  \label{tab:time-NSAA}%
\end{table}%

}

%%%%%%%%%%%%%%%%%
\end{document}